\documentclass[copyedit]{informs1}  

\usepackage{multirow}
\usepackage{endnotes}
\usepackage{enumerate}
\usepackage{rotating}
\usepackage{subfigure}
\usepackage{graphicx}
\usepackage{multirow}
\usepackage{hhline}
\usepackage{amsmath}
\usepackage{epstopdf}
\usepackage{eufrak}

\usepackage{setspace,graphics,booktabs,color}
\usepackage{multirow}

\usepackage{float}
\floatstyle{ruled}
\newfloat{algorithm}{htp}{lop}
\floatname{algorithm}{Algorithm}

\usepackage{natbib}
 \bibpunct[, ]{(}{)}{,}{a}{}{,}%
 \def\bibfont{\small}%
\makeatletter

\newcommand{\Rmnum}[1]{\expandafter\@slowromancap\romannumeral #1@}
\makeatother

\TheoremsNumberedThrough     

\EquationsNumberedThrough    

\definecolor{DSgray}{cmyk}{0,1,0,0}

\begin{document}

%
%
%
%
%
%
%
%
%
%
%

\RUNAUTHOR{}
%
\RUNTITLE{}

\TITLE{\Large Comparison-Based Algorithms for One-Dimensional Stochastic Convex Optimization}

\ARTICLEAUTHORS{
	
	\AUTHOR {Xi Chen} \AFF{Stern School of Business, New York
		University, New York City, NY 10012, {xchen3@stern.nyu.edu}}
	
	\AUTHOR {Qihang Lin} \AFF{Tippie College of Business, University of
		Iowa, Iowa City, IA, 52245, {qihang-lin@uiowa.edu}}
	
	\AUTHOR {Zizhuo Wang} \AFF{Department of Industrial and Systems
		Engineering, University of Minnesota, Minneapolis, MN 55455,
	} 
	
	\AFF{Institute of Data and Decision Analytics, The Chinese University of Hong Kong, Shenzhen,
		{zwang@umn.edu}} 
}

\ABSTRACT{Stochastic optimization finds a wide range of applications
in operations research and management science. However, existing
stochastic optimization techniques usually require the information of
random samples (e.g., demands in the newsvendor problem) or the
objective values at the sampled points (e.g., the lost sales cost),
which might not be available in practice. In this paper, we consider
a new setup for stochastic optimization, in which the decision maker
has access to only comparative information between a random sample
and two chosen decision points in each iteration. We propose a
comparison-based algorithm (CBA) to solve such problems
 in one dimension with convex objective functions.
Particularly, the CBA properly chooses the two points in each
iteration and constructs an unbiased gradient estimate for the
original problem. We show that the CBA achieves the same convergence
rate as the optimal stochastic gradient methods (with the
samples observed).  We also consider extensions of our approach
to multi-dimensional quadratic problems as well as problems with
non-convex objective functions. Numerical experiments show that the
CBA performs well in test problems.}


\HISTORY{This version on \today}

\maketitle
%



\section{Introduction}
\label{sec:intro}  In this paper, we consider the following
stochastic optimization problem:
\begin{eqnarray}\label{formulation}
\min_{\ell\le x\le u} & H(x) = {\mathbb E}_{\xi}\left[h(x,
\xi)\right]
\end{eqnarray}
where $-\infty\le\ell\le u \le +\infty$ and $\xi$ is a random variable.\footnote{Throughout the
note, when $\ell = -\infty$ ($u = +\infty$, resp.), the notation
$\ell \le x$ ($x\le u$, resp.) will be interpreted as $x > -\infty$
($x < +\infty$, resp.).} This problem has many applications and is
fundamental for stochastic optimization. For example,
\begin{enumerate}
\item If $h(x,\xi) = (x - \xi)^2$ with $\ell = -\infty$ and $u = +\infty$, then the problem is to find the
expectation of $\xi$.
\item If $h(x,\xi) = {\mathfrak{h}}\cdot (x - \xi)^+ + {\mathfrak{b}} \cdot(\xi - x)^+$, then
the problem is the classical newsvendor problem with unit holding
cost ${\mathfrak{h}}$ and unit backorder cost ${\mathfrak{b}}$. It
can also be viewed as the problem of finding the
${\mathfrak{b}}/({\mathfrak{h}}+{\mathfrak{b}})$-th quantile of
$\xi$ when $\ell = -\infty$ and $u = +\infty$. Furthermore, one  can
consider a more general version of this problem in which
\begin{eqnarray*}
h(x, \xi) = \left\{\begin{array}{ll} \mathfrak{h}_+(x,\xi) & \mbox{ if } x
\ge \xi\\
\mathfrak{h}_{-}(x,\xi) & \mbox{ if } x < \xi.
\end{array}
\right.
\end{eqnarray*}
 This problem can be viewed as a newsvendor problem
with general holding and backorder costs (see, e.g.,
\citealt{halman} and references therein for discussions of this
problem, where $\mathfrak{h}_+(x,\xi)$ is a general holding cost
function and $\mathfrak{h}_{-}(x,\xi)$ is a general backorder cost
function). It can also be viewed as a single period appointment
scheduling problem with general waiting and overtime costs (see,
e.g., \citealt{gupta}) or a staffing problem with general underage
and overage costs (see, e.g., \citealt{kolker}).
\item If $h(x,\xi) = - x \cdot 1(\xi \ge x)$, then the problem can
be viewed as an optimal pricing problem where $x$ is the price set
by the seller, $\xi$ is the valuation of each customer, and
$h(x,\xi)$ is the negative of the revenue obtained
from the customer (a customer purchases at price $x$ if and only if his/her valuation $\xi$ is greater than or equal to $x$).\\
\end{enumerate}

 In many practical situations, the distribution of $\xi$
(whose c.d.f. will be denoted by $F(\cdot)$) is unknown a priori.
Existing stochastic optimization techniques for \eqref{formulation}
usually require sampling from the distribution of $\xi$ and use
random samples to update the decision $x$ toward optimality. In the
existing methods, it is assumed that either the random samples of
$\xi$ are fully observed or the objective value $h(x,\xi)$ under a
decision $x$ and random samples $\xi$ can be observed. However, such
information may not always be available in practice (see Examples
\ref{example:mean}-\ref{example:revenue} below).

 In this paper, we investigate whether having full sample
information is always critical for solving one-dimensional
stochastic convex optimization problems. We realize that, in some
cases where full samples are not observed, comparative relation
between a chosen decision variable and a random sample may still be
accessible. This motivates us to study stochastic optimization with
only the presence of comparative information.
Specifically, given a decision $x$, a sample $\xi$ is drawn from the
underlying distribution, and we assume that we only have information
about whether $\xi$ is greater than or less than (or equals to) $x$.
In addition, after knowing the comparative relation between $x$ and
$\xi$, we further assume that we can choose another point $z$ and
obtain information about whether $\xi$ is greater than or less than
(or equals to) $z$. Such a $z$ is not as a decision
variable but a randomly sampled point.
We show that, in fact, having the comparative information in this way can sometimes be sufficient
for solving \eqref{formulation}. In the following, we list several
scenarios in which such situations may arise:
\begin{example}\label{example:mean}
Suppose $x$ represents a certain feature of a product (e.g., size or
taste, etc) and $\xi$ is the preference of each customer about that
feature, and the firm selling this product would like to find out
the average preference of the customers (or equivalently, to find
the optimal offering to minimize the expected customer
dissatisfaction, measured by $h(x,\xi) = (x - \xi)^2$). Such
a firm faces a stochastic optimization problem described in the
first example above. In many cases, it is hard for a customer to
give an exact value for his/her preference (i.e., the exact value of
his/her $\xi$). However, it is quite plausible that the customer can
report comparative relation between his/her preferred value
of the feature and the actual value of the feature of the product
presented to him/her (e.g., whether the product should be larger or
whether the taste should be saltier). Furthermore, it is
possible to ask one customer to compare his/her preferred value with
two different values of the same feature of the product, for
example, by giving the customer two different samples. Moreover,
the second sample may be given in a customer satisfaction survey,
and the customer will not count the second sample toward its
(dis)satisfaction value. Therefore, such a scenario fits the setting
described above.
\end{example}
\begin{example}\label{example:inventory}
In a newsvendor problem, it is sometimes hard to observe the exact
demand in each period due to demand censorship. In such situations,
one does not have direct access to the sample point (the demand) nor
does one have access to the cost in the corresponding period (the
lost sales cost).\footnote{There is a vast literature on
newsvendor/inventory problems with censored demand. For some recent
references, we refer the readers to \cite{Ding}, \citet{Bensoussan},
and \cite{Besbes_newvendor}.} However, the seller usually has
comparative information between the realized demand and the
chosen inventory level  (e.g., by observing if there is a stock out
or a leftover). Moreover, by allowing the seller to make a one-time
additional ordering in each time period (this ability is sometimes
called the quick response ability for the seller, see e.g.,
\citealt{Cachon}), it is possible that one can obtain such
information at two points. In such cases, the firm will face a
newsvendor problem as described in the second example above, and
thus it will correspond to the setting in our problem.
\end{example}
\begin{example}\label{example:revenue}
In a revenue management problem, by offering a price to each
customer, the seller can observe whether the customer purchased the
product, and the seller faces a stochastic optimization
problem described in the third example above. In practice, it is
hard to ask the customer to report a true valuation of the
product. However, it is possible to ask the customer in a
market survey whether he or she will purchase the product at a different
price. Such an example can also be extended to a divisible product
case in which a customer can buy a continuous amount of a product with
a maximum of $1$. In this case, the $h$ function can be
redefined as $h(x,\xi) = -x\min\{1, g(x,\xi)\}$ where $x$ is the
offered price, $g(x,\xi)$ is the unconstrained purchase amount of
the customer, and $\xi$ is the maximum price this customer is
willing to buy the full amount of this product (i.e., $g(x, \xi)$ is
decreasing in $x$ with $g(x, \xi)
> 1$ when $x < \xi$ and $g(x, \xi) <1$ when $x > \xi$). Such a
purchase behavior can be explained by a quadratic utility function
of the customer, which is often used in the literature (see e.g.,
\citealt{candogan}). For the seller, by observing whether the
customer buys the full amount of the product, he or she can infer whether
an offered price is greater than or less than the $\xi$ value of
this customer.\footnote{There is abundant recent literature
that studies the setting in which the seller can observe the full
information of $\xi$ for each customer. In particular, it has been
shown that in this case, the seller can obtain asymptotic optimal
revenue as the selling horizon grows. For a review of this
literature, we refer the readers to \cite{denboer}.}
\end{example}

In this paper, we propose an efficient algorithm to solve the
above-described stochastic optimization problem. More precisely, we
propose a stochastic approximation algorithm that only utilizes
comparative information between each sample point and two chosen decision
points in each iteration. We show that by properly choosing the two
points (one point has to be chosen randomly according to a
specifically designed distribution), we can obtain unbiased gradient
estimates for the original problem.\footnote{
If $h(x,\xi)$ is piecewise linear with two pieces, e.g., $h(x,\xi) = {\mathfrak{h}}\cdot (x - \xi)^+ + {\mathfrak{b}} \cdot(\xi - x)^+$,
only comparing $x$ and $\xi$ may be sufficient to compute the stochastic gradient $h_x'(x,\xi)$ (that equals $\mathfrak{h}$ or $-\mathfrak{b}$). }
The unbiased gradient estimates
will in turn give rise to efficient algorithms based on a standard
stochastic gradient method (we will review the related
literature shortly). Under some mild conditions, we show that if the
original problem is convex, then our algorithm will achieve a
convergence rate of $O(1/\sqrt{T})$ for the objective value
 (where $T$ is the number of iterations); if the
original problem is strongly convex, then the convergence rate can
be improved to $O(1/T)$. Moreover, the information at two points is
necessary in this setting as we show that only knowing
comparative information between the sample and one point in
each iteration is insufficient for any algorithm to converge to the
optimal solution (see Example \ref{ex:onecomp}). We also perform
several numerical experiments using our algorithm. The
experimental results show that our algorithms are indeed
efficient, with convergence speed in the same order compared to the
case when one has direct observations of the samples.  We also
extend our algorithm to a multi-dimensional setting with quadratic
objective function, a setting with non-convex objective function and
a setting in which
multiple comparisons can be conducted in each iteration.\\

{\bf\noindent Literature Review.} Broadly speaking, our work falls
into the area of stochastic optimization, a subject on which there is vast literature. 
There has been a vast
literature on stochastic optimization. For a comprehensive review of
this literature, we refer the readers to \citet{shapiro}. In
particular, in this literature, it is usually assumed that one has
access to random samples (or alternatively, the objective values at
the sampled points). Two main types of algorithms have been
proposed, namely   the sample average approximation (SAA) method (see, e.g., \citealt{shapiro})
and the stochastic approximation (SA)
methods (see, e.g., \citealt{robbins}, \citealt{kiefer}). A
typical SAA method collects a number of samples from the underlying
distribution and uses the sample averaged objective function to
approximate the expected value function. This method has been widely
used in many operations management problems (see, e.g.,
\citealt{levi2015data, rudin2014big}). In contrast, the SA approach
(e.g., the stochastic gradient descent) is usually an iterative
algorithm. In each iteration, a new sample (or a small batch of new
samples) is drawn, and a new iterate is computed using the new
sample(s). Our work belongs to the category of SA. In the following,
we shall focus our literature review on the stochastic approximation
methods.


If the objective function is convex, then various stochastic
gradient methods can guarantee, under slightly different
assumptions, that the objective value of the iterates converges to
the optimal value in a rate of $O(1/\sqrt{T})$ after $T$ iterations
(see, e.g., \citealt{Nemirovski:09, Duchi:09, Hu:09, Xiao:10,
Lin:11, Chen:12, Lan:10a, Lan:10b, Rakhlin:12, Lan:10c, Shamir:13,
Hazan:14}). Furthermore, this convergence rate is known to be
optimal \citep{Nemirovski:83}. When the objective function is
strongly convex, some stochastic gradient methods can obtain an
improved convergence rate of $O(\log T/T)$ \citep{Duchi:09,Xiao:10}.
More recently, several papers have further improved the convergence
rate to $O(1/T)$. Among those papers, there are three different
methods used: (a) accelerated stochastic gradient method with
auxiliary iterates besides the main iterate \citep{Hu:09, Lin:11,
Lan:10a, Lan:10b, Chen:12, Lan:10c}; (b) averaging the historical
solutions \citep{Rakhlin:12,Shamir:13}; and (c) multi-stage
stochastic gradient method that periodically restarts
\citep{Hazan:14}. Again, the convergence rate of $O(1/T)$ has been shown to be optimal by
\citet{Nemirovski:83} for strongly
convex problems.

Apart from the setting of minimizing a single objective function,
stochastic gradient methods can also be applied to online learning
problems (for a comprehensive review, see \citealt{Shalev-Shwartz:2012}) 
where a sequence of functions is presented to a decision
maker who needs to provide a solution sequentially to each function
with the goal of minimizing the total regret. It is known that the
stochastic gradient methods can obtain a regret of $O(\sqrt{T})$
after $T$ decisions if the functions presented are convex. Moreover,
this regret has been shown to be optimal \citep{Bianchi:book06}.
When the functions are strongly convex, \citet{Zinkevish:03},
\citet{Duchi:09}, \citet{Xiao:10} and \citet{Duchi:Ada} show that
the regret can be further improved to $O(\log{T})$.

To distinguish our work from the above, we note that all of the
above works have assumed that either the sample is directly
accessible (one can observe the value of each sample) or the
objective value corresponding to the decision variable and the
current sample is accessible. In either case, it is easy to obtain
an estimate of the gradient of the objective function. In contrast,
in our case, we do not have access to the sample or the objective
value. Instead, we only have comparative information between each
sample and two chosen decision points. Indeed as we shall discuss in the next
section, one of the main challenges in our problem is to use this
very limited information to construct an unbiased gradient for the
original problem and then further use it to find the optimal
solution. Our contribution is to show that the same order of
convergence rate can still be achieved under this setting with less
information.


\section{Main Results}
\label{sec:results}

In this paper, we make the following assumptions:

\begin{assumption}\label{assumption:onedimensional}
\begin{enumerate}
\item []
\item[] (A1) The random variable $\xi$ follows a continuous
distribution.
\item[] (A2) For each $\xi$, $h(x,\xi)$ is continuously differentiable with respect to $x$ on $[\ell, \xi)$ and
$(\xi, u]$ with the derivative denoted by $h'_x(x,\xi)$.
Furthermore, for any $x\in [\ell, u]$, $h'_{-}(x) :=
\underset{z\rightarrow x-}{\lim} h'_x(x, z)$ and $h'_{+}(x) :=
\underset{z\rightarrow x+}{\lim} h'_x(x, z)$ exist and are finite.
\item[] (A3) For any $x\in[\ell, u]$, $x\neq \xi$, $h_{x,\xi}''(x, \xi)= \frac{\partial^2 h(x,\xi)}{\partial \xi \partial x}$ exists.
\item[] (A4) $H(x)$ in \eqref{formulation} is differentiable and $\mu$-convex on $[\ell,u]$ for some $\mu \ge 0$, namely,
\begin{eqnarray}
\label{eq:strconv} H(x_2) \ge H(x_1) + H'(x_1)(x_2-x_1) +
\frac{\mu}{2} (x_2-x_1)^2, \quad \forall x_1, x_2 \in [\ell, u].
\end{eqnarray}
Moreover, $\mathbb{E}_{\xi}(h'_x(x, \xi))=H'(x)$ for all
$x\in[\ell,u]$.
\item[] (A5) Either of the following statements is true:
\begin{enumerate}[(a)]
\item[a.] There exists a constant $K_1$ such that $\mathbb{E}_{\xi}(h_{x}'(x,\xi))^2\leq K_1^2$ for
any $x\in [\ell, u]$;
\item[b.] $H'(x)$ is $L$-Lipschitz continuous on $[\ell, u]$. Furthermore, there exists a constant $K_2$ such that
\begin{eqnarray*}
\label{eq:boundg} \mathbb{E}_{\xi}(h_x'(x,\xi)-H'(x))^2\leq K_2^2,
\quad \forall x\in[\ell, u].
\end{eqnarray*}
\end{enumerate}
\end{enumerate}
\end{assumption}

Now we make several comments on the above assumptions. The first
assumption that $\xi$ is continuously distributed is mainly for the
ease of discussion. In fact, all of our results will continue to
hold as long as with probability $1$, for all iterates $x$ in our
algorithm, ${\mathbb P}(\xi = x) = 0$. We shall revisit this
assumption in Section \ref{sec:discussions}. Assumptions A2-A4 are
some regularity assumptions on the functions $h$ and $H$.
Particularly, the last point of Assumption A4 is satisfied under
many cases, for example, when $h'_x(x,\xi)$ is continuous in $\xi$
and $\xi$ is supported on a finite set (\citealt{widder}), or when
$h(x,\xi)$ is convex in $x$ for each $\xi$ (by monotone convergence
theorem). When \eqref{eq:strconv} holds with $\mu>0$, we call $H$ a $\mu$-\emph{strongly convex} function.
The last assumption states that the partial derivative
$h'_x(x,\xi)$ has uniformly bounded second-order moment or variance.
This is used to guarantee that the step in each iteration in our
algorithm has bounded variance, which is a common assumption in
stochastic approximation literature (see, e.g.,
\citealt{Nemirovski:09}, \citealt{Duchi:09}, \citealt{Lan:10c}).

In addition, Assumption \ref{assumption:onedimensional} is not hard
to satisfy in our examples mentioned earlier. Specifically, for
Example \ref{example:mean}, it satisfies Assumption
\ref{assumption:onedimensional} when $\xi$ is continuously
distributed and has finite variance. For Example
\ref{example:inventory}, it satisfies Assumption
\ref{assumption:onedimensional} when $\xi$ is continuously
distributed and the cost functions are linear. When the cost
functions are nonlinear ($\mathfrak{h}_+$ and $\mathfrak{h}_{-}$
respectively), it satisfies Assumption
\ref{assumption:onedimensional} if both $\mathfrak{h}_{+}$ and
$\mathfrak{h}_{-}$ are second-order continuously differentiable on
their respective domains, have bounded first-order derivatives (for
example, when $x$ and $\xi$ are restricted to finite intervals), and
$h$ is convex in $x$. For Example \ref{example:revenue}, it
satisfies Assumption \ref{assumption:onedimensional} under the
divisible case when $\xi$ is a continuous random variable and the
expected revenue function $x{\mathbb E}_{\xi}\min\{g(x,\xi), 1\}$ is
concave on $[\ell, u]$, which holds, for example, when $g(x,\xi)$ is
a piecewise linear function and when the range $[\ell, u]$ is
small.\footnote{The indivisible product case with $h(x,\xi) = - x
\cdot 1(\xi \ge x)$ does not satisfy Assumption (A4). In particular,
it does not satisfy $\mathbb{E}_{\xi}(h'_x(x, \xi))=H'(x)$ (it does
satisfy all the other assumptions under mild conditions though). In
order to satisfy $\mathbb{E}_{\xi}(h'_x(x, \xi))=H'(x)$, it is
sufficient that $h(x,\xi)$ is continuous in $x$, which holds in the
divisible product case.}



In the following, we propose a comparison-based algorithm (CBA) to
solve (\ref{formulation}). Let $\Xi$ denote the support of $\xi$ with
$-\infty \le \underline{s} :=\inf\{\Xi\} \le \bar{s} := \sup\{\Xi\} \le
+\infty$. The algorithm requires specification of two functions,
$f_{-}(x, z)$ and $f_{+}(x, z)$, which need to satisfy the following
conditions.
\begin{itemize}
\item (C1) $f_{-}(x,z) = 0$ for all $z \geq x$ and $f_{-}(x,z) > 0$ for all $\underline{s} \le z < x$. In addition,  for all $x$, we have $\int_{-\infty}^{x-}
f_{-}(x,z) dz = 1$.
\item (C2) $f_{+}(x,z) =0$ for all $z \leq x$ and $f_{+}(x,z) > 0$ for all $\bar{s} \ge z > x$. In addition, for all $x$, we have $\int_{x+}^{\infty}
f_{+}(x,z) dz = 1$.
\item (C3) There exists a constant $K_3$ such that $\int_{\underline{s}}^{x-}\frac{F(z)(h''_{x,z}(x,z))^2}{f_{-}(x,z
)}dz\leq K_3$ and
$\int_{x+}^{\bar{s}}\frac{(1-F(z))(h''_{x,z}(x,z))^2}{f_{+}(x,z)}dz\leq
K_3$ for all $x\in[\ell, u]$, where $F(\cdot)$ is the c.d.f. of $\xi$.
\end{itemize}

Note that, for any given $x\in[\ell, u]$,
$f_{-}(x,z)$ and $f_{+}(x,z)$ essentially define two
density functions of $z$ on $(-\infty,x]$ and $[x,+\infty)$. (We will discuss how to choose $f_{-}(\cdot,
\cdot)$ and $f_{+}(\cdot, \cdot)$ in Section \ref{sec:choice of f}.)
Also, we note that $\underline{s}$ and $\bar{s}$ need not to be
known in advance. If one is unsure about $\underline{s}$ ($\bar{s}$,
resp.), then one can choose a sufficiently small (large, resp.)
value, or just choose $\underline{s} = -\infty$ ($\bar{s} = +\infty$,
resp.). Condition (C3) is a technical condition and
may not be straightforward to verify at the first glance. However,
in Section \ref{sec:choice of f}, we show that under mild conditions
(e.g., $\xi$ has a light tail and $h''_{x,z}(x,z)$ is uniformly
bounded), it is not hard to choose the functions $f_{-}(x,z)$ and
$f_{+}(x,z)$ such that condition (C3) is satisfied (we will leave
the detailed discussions in Section \ref{sec:choice of f}). Next,
in Algorithm \ref{alg:cba1}, we describe the detailed procedure of
the CBA.

\begin{algorithm}[h] \caption{Comparison-Based Algorithm
(CBA):\label{alg:cba1}}
\begin{enumerate}
\item {\bf Initialization.} Set $t = 1$, $x_1 \in[\ell, u]$. Define $\eta_t$
for all $t\ge 1$. Set the maximum number of iterations $T$.
Choose functions $f_{-}(x,z)$ and $f_{+}(x,z)$ that satisfy
(C1)-(C3).
\item {\bf Main iteration.} Sample $\xi_t$ from the distribution of $\xi$. If $\xi_t = x_t$, then resample
$\xi_t$ until it does not equal $x_t$. (This step will always
terminate in a finite number of steps as long as $\xi$ is not
deterministic.)
\begin{enumerate}
\item If $\xi_t < x_t$, then generate $z_t$ from a
distribution on $(-\infty, x_t] $ with p.d.f. $f_{-}(x_t,z_t)$. Set
\begin{eqnarray}\label{gradient_case1}
g(x_t, \xi_t, z_t) = \left\{\begin{array}{ll} h'_{-}(x_t), & \mbox{
if } z_t < \xi_t,
\\
h'_{-}(x_t) - \frac{h''_{x,z}(x_t,z_t)}{f_{-}(x_t,z_t)}, & \mbox{ if
} z_t \ge \xi_t. \end{array}\right.
\end{eqnarray}
\item If $\xi_t > x_t$, then generate $z_t$ from a
distribution on $[x_t, +\infty)$ with p.d.f. $f_{+}(x_t, z_t)$. Set
\begin{eqnarray}\label{gradient_case2}
g(x_t, \xi_t, z_t) = \left\{\begin{array}{ll} h'_{+}(x_t), & \mbox{
if } z_t > \xi_t,
\\
h'_{+}(x_t) + \frac{h''_{x,z}(x_t,z_t)}{f_{+}(x_t,z_t)}, & \mbox{ if
} z_t \le \xi_t.\end{array} \right.
\end{eqnarray}
\end{enumerate}
Let
\begin{eqnarray} \label{eq:proj}
x_{t+1}&=&\text{Proj}_{[\ell, u]}\left(x_t-\eta_tg(x_t, \xi_t,
z_t)\right) = \max\left(\ell, \min\left(u, x_t - \eta_tg(x_t, \xi_t,
z_t)\right)\right).
\end{eqnarray}
\item {\bf Termination.} Stop when $t\ge T$. Otherwise, let $t\leftarrow t+ 1$ and go
back to Step 2.
\item {\bf Output.} $\textbf{CBA}(x_1,T,\{\eta_t\}_{t=1}^T) = \bar x_T=\frac{1}{T}\sum_{t=1}^Tx_t$.
\end{enumerate}
\end{algorithm}


In iteration $t$ of CBA, the solution $x_t$ will be updated based on
two random samples. First, a sample $\xi_t$ is drawn from the
distribution of $\xi$. In contrast to existing
stochastic gradient methods, CBA does not require exactly observing
$\xi_t$ but only needs to know whether $\xi_t<x_t$ or $\xi_t>x_t$.
Based on the result of the comparison between  $\xi_t$ and $x_t$, a
second sample $z_t$ is drawn from the density function
$f_{+}(x_t,z)$ or $f_{-}(x_t,z)$. An unbiased stochastic gradient,
$g(x_t, \xi_t, z_t)$, of $H(x_t)$ is then constructed and used to
update $x_t$ with the standard gradient descent step.

We note that the output of Algorithm \ref{alg:cba1} is the average
of the historical solutions $\bar{x}_T =
\frac{1}{T}\sum_{t=1}^Tx_t$. This is because the convergence of
objective value is established based on $\bar x_T$. However,
Algorithm \ref{alg:cba1} can be applied to the online learning
setting where one can use the solution $x_t$ as the decision in each
stage $t$ and obtain the desired expected total regret (see Proposition
\ref{thm:CBA1}-\ref{thm:restart}). We have the following proposition
about the stochastic gradient $g(x_t, \xi_t , z_t )$ in CBA.

\begin{proposition}\label{prop:unbiased}
Suppose $f_{-}(x,z)$ and $f_{+}(x,z)$ satisfy (C1)-(C3) and
Assumption \ref{assumption:onedimensional} holds. Then
\begin{enumerate}
\item
${\mathbb E}_{z }g(x , \xi , z ) = h'_x(x ,\xi )$, for all $x\in
[\ell, u]$, $x \neq \xi $.
\item ${\mathbb E}_{z , \xi } g(x , \xi , z ) =
H'(x )$, for all $x \in [\ell, u]$.
\item If Assumption A5(a) holds, then $\mathbb{E}_{z,\xi} (g(x, \xi, z))^2\leq G^2:= K_1^2+2K_3$. If Assumption A5(b) holds, then ${\mathbb E}_{z,\xi} (g(x, \xi, z)-H'(x))^2\leq \sigma^2:= K_2^2+2K_3$.
\end{enumerate}
\end{proposition}

{\bf\noindent Proof of Proposition \ref{prop:unbiased}.} First, we
consider the case when $\xi  < x $. We have
\begin{eqnarray*}
{\mathbb E}_{z } g(x , \xi , z ) = h_{-}'(x ) - \int_{\xi }^{x -}
h''_{x,z}(x , z)dz = h'_x(x , \xi).
\end{eqnarray*}
Similarly, when $\xi > x$,
\begin{eqnarray*} {\mathbb
E}_{z } g(x , \xi , z ) = h_{+}'(x ) + \int_{x +}^{\xi } h''_{x,z}(x
, z)dz =  h'_x(x , \xi).
\end{eqnarray*}
Thus the first conclusion of the proposition is proved. The second
conclusion of the proposition follows from Assumption A1 (which
ensures $\xi=x$ is a zero-measure event) and Assumption A4.

Next, we show the first part of the third conclusion when Assumption
A5(a) is true. If $\xi< x$, then we have
\begin{eqnarray}
\label{eq:A1} {\mathbb E}_{z} (g(x, \xi, z))^2&=&\int_{-\infty}^{x-
}(h'_{-}(x))^2f_{-}(x ,z )dz + \int_{\xi }^{x -} \left(-
2h'_{-}(x)\frac{h''_{x,z}(x,z)}{f_{-}(x ,z
)}+\left(\frac{h''_{x,z}(x,z)}{f_{-}(x ,z
)}\right)^2\right)f_{-}(x ,z )dz \nonumber \\
&=&(h'_{-}(x))^2-2h'_{-}(x)(h'_{-}(x)-h'_x(x , \xi
))+ \int_{\xi }^{x -}\frac{(h''_{x,z}(x,z))^2}{f_{-}(x ,z )}dz\\
&\leq&(h'_x(x , \xi))^2+\int_{\xi }^{x
-}\frac{(h''_{x,z}(x,z))^2}{f_{-}(x ,z )}dz.\nonumber
\end{eqnarray}
where the last inequality is because $a^2 + b^2 \ge 2ab$ for any
$a$, $b$. By similar arguments, if $\xi> x$, then
\begin{eqnarray*}
\label{eq:A2} {\mathbb E}_{z} (g(x, \xi, z))^2 &\leq&(h'_x(x ,
\xi))^2+ \int_{x + }^{\xi}\frac{(h''_{x,z}(x,z))^2}{f_{+}(x ,z )}dz.
\end{eqnarray*}
These two inequalities and Assumption A5(a) further imply
\begin{eqnarray}
{\mathbb E}_{z,\xi} (g(x, \xi,
z))^2&\leq&K_1^2+\int_{\underline{s}}^{x-}\left(\int_{\xi }^{x
-}\frac{(h''_{x,z}(x,z))^2}{f_{-}(x ,z )}dz\right)dF(\xi)+\int_{x +
}^{\bar{s}}\left(\int_{x + }^{\xi}\frac{(h''_{x,z}(x,z))^2}{f_{+}(x
,z )}dz\right)dF(\xi)\nonumber \\
&=&K_1^2+\int_{\underline{s}}^{x-}\frac{F(z)(h''_{x,z}(x,z))^2}{f_{-}(x
,z )}dz+\int_{x +
}^{\bar{s}}\frac{(1-F(z))(h''_{x,z}(x,z))^2}{f_{+}(x ,z
)}dz\label{eq:A1A2} \\ &\leq&K_1^2+2K_3, \nonumber
\end{eqnarray}
where the interchanging of integrals in the equality is justified by
Tonelli's theorem and the last inequality is due to (C3).

Next, we show the second part of the third conclusion when
Assumption A5(b) is true. If $\xi< x$, then following the similar
analysis as in \eqref{eq:A1}, we have
\begin{eqnarray*}
{\mathbb E}_{z} (g(x, \xi, z)- h'_x(x ,\xi ))^2 ={\mathbb E}_{z}
(g(x, \xi, z))^2- (h'_x(x ,\xi ))^2 \leq \int_{\xi }^{x
-}\frac{(h''_{x,z}(x,z))^2}{f_{-}(x ,z )}dz.
\end{eqnarray*}
Similarly, if $\xi> x$, then
\begin{eqnarray*}
{\mathbb E}_{z} (g(x, \xi, z)- h'_x(x ,\xi ))^2&\leq& \int_{x +
}^{\xi}\frac{(h''_{x,z}(x,z))^2}{f_{+}(x ,z )}dz.
\end{eqnarray*}
By using the same argument as in \eqref{eq:A1A2}, we have
\begin{eqnarray*}
{\mathbb E}_{z, \xi} (g(x, \xi, z)- h'_x(x ,\xi ))^2 \le 2K_3.
\end{eqnarray*}
Finally, we note that,
\begin{eqnarray*}
{\mathbb E}_{z,\xi} (g(x, \xi, z)-H'(x))^2&=&{\mathbb E}_{z,\xi}
(g(x, \xi, z)- h'_x(x ,\xi ))^2+
\mathbb{E}_{\xi}(h_x'(x,\xi)-H'(x))^2.
\end{eqnarray*}
Therefore, when Assumption A5(b) holds, we have ${\mathbb E}_{z,\xi}
(g(x, \xi, z)-H'(x))^2 \le K_2^2 + 2K_3$. Thus the proposition
holds. $\hfill\Box$\\

Proposition \ref{prop:unbiased} shows that in the CBA, the gradient
estimate $g(x, \xi, z)$ is an unbiased estimate of the true gradient
at $x$ and can be utilized as a stochastic gradient of $H(x)$. Note
that such an unbiased gradient is generated without accessing the
sample $\xi$ itself nor the value of the objective function $h(x,\xi)$ at the sampled
point. The only information used in generating the unbiased gradient
is comparative information between the sample and two points.


 Using the unbiased stochastic gradient, we can
characterize the convergence results of the CBA in the following
propositions.
\begin{proposition}
    \label{thm:CBA1}
    Suppose $\mu=0$.
    Let $G^2$ and $\sigma^2$ be defined
    as in Proposition~\ref{prop:unbiased} and $x^*$ be any optimal
    solution to (\ref{formulation}).
    \begin{itemize}
        \item If Assumption A5(a) holds, then by choosing $\eta_t=\frac{1}{\sqrt{T}}$, the CBA ensures that
        \begin{eqnarray*}
            \mathbb{E}(H(\bar
            x_T)-H(x^*))\leq\frac{(x_1-x^*)^2}{2\sqrt{T}}+\frac{G^2}{2\sqrt{T}} \qquad\mbox{and}\qquad \sum_{t=1}^{T}\mathbb{E}(H(x_t)-H(x^*))\leq
            \frac{\sqrt{T}}{2}(x_1-x^*)^2+\frac{\sqrt{T}G^2}{2}.
        \end{eqnarray*}
        If, in addition, $u$ and $\ell$ are finite, then by choosing $\eta_t=\frac{1}{\sqrt{t}}$, the CBA ensures that
        \begin{eqnarray*}
            \mathbb{E}(H(\bar
            x_T)-H(x^*))\leq\frac{(u-\ell)^2}{2\sqrt{T}}+\frac{G^2}{\sqrt{T}} \qquad \mbox{and}\qquad \sum_{t=1}^{T}\mathbb{E}(H(x_t)-H(x^*))\leq
            \frac{\sqrt{T}}{2}(u-\ell)^2+\sqrt{T}G^2.
        \end{eqnarray*}
        \item If Assumption A5(b) holds, then by choosing $\eta_t=\frac{1}{L+\sqrt{T}}$, the CBA ensures that
        \begin{eqnarray*}
            \mathbb{E}(H(\bar
            x_T)-H(x^*))&\leq&\frac{L+\sqrt{T}}{2T}(x_1-x^*)^2+\frac{H(x_1)-H(x^*)}{T}+\frac{\sigma^2}{2\sqrt{T}}, \mbox{   and}\\
            \sum_{t=1}^{T}\mathbb{E}(H(x_t)-H(x^*))&\leq&
            \frac{L+\sqrt{T}}{2}(x_1-x^*)^2+H(x_1)-H(x^*)+\frac{\sqrt{T}\sigma^2}{2}.
        \end{eqnarray*}
 If, in addition, $u$ and $\ell$ are finite, then by choosing $\eta_t=\frac{1}{L+\sqrt{t}}$, the CBA ensures that
        \begin{eqnarray*}
            \mathbb{E}(H(\bar
            x_T)-H(x^*))&\leq&\frac{(u-\ell)^2}{2\sqrt{T}}+\frac{H(x_1)-H(x^*)}{T}+\frac{L(x_1-x^*)^2}{2T}+\frac{\sigma^2}{\sqrt{T}}, \mbox{   and}\\
            \sum_{t=1}^{T}\mathbb{E}(H(x_t)-H(x^*))&\leq&
            \frac{\sqrt{T}}{2}(u-\ell)^2+H(x_1)-H(x^*)+\frac{L}{2}(x_1-x^*)^2+\sqrt{T}\sigma^2.
        \end{eqnarray*}
     \end{itemize}
\end{proposition}

Proposition \ref{thm:CBA1} gives the performance of the CBA when
$\mu=0$. When $\mu=0$ and $u$ and/or $\ell$ are infinite, the
stepsize $\eta_t$ is chosen to be a constant depending on the total
number of iterations $T$ in order to achieve the optimal convergence
rate (i.e., $\eta_t=\frac{1}{\sqrt{T}}$ when Assumption A5(a) holds
and $\eta_t=\frac{1}{L+\sqrt{T}}$ when Assumption A5(b) holds).
Therefore, one needs to determine the total number of iterations $T$
before running the optimization algorithm for the computation of the
stepsize. When $\mu=0$ and $u$ and $\ell$ are finite, the stepsize
$\eta_t$ can be chosen as a decreasing sequence in $t$ (i.e.,
$\eta_t=\frac{1}{\sqrt{t}}$ when Assumption A5(a) holds and
$\eta_t=\frac{1}{L+\sqrt{t}}$ when Assumption A5(b) holds). In such
a case, one does not need to pre-specify the total number of
iterations $T$. The convergence result when Assumption A5(a) holds
has a proof similar to \citet{Nemirovski:09} and \citet{Duchi:09}
except using our comparison-based stochastic gradient. Also, the
convergence result when Assumption A5(b) holds is largely built upon
the results in \cite{Lan:10a}. The detailed proof of the proposition are
given in Appendix \ref{appendix:proof_of_theorem1-3}.

Next, we have a further result when $\mu > 0$.
\begin{proposition}
    \label{thm:CBA1sc}
    Suppose $\mu>0$.
    Let $G^2$ and $\sigma^2$ be defined
    as in Proposition~\ref{prop:unbiased} and $x^*$ be any optimal
    solution to (\ref{formulation}).
    \begin{itemize}
        \item If Assumption A5(a) holds, then by choosing $\eta_t=\frac{1}{\mu t}$, the CBA ensures that
        \begin{eqnarray*}
            \mathbb{E}(H(\bar x_T)-H(x^*))\leq\frac{G^2}{2\mu}\frac{\log T +1}{T}
            \qquad\mbox{and}\qquad\sum_{t=1}^{T}\mathbb{E}(H(x_t)-H(x^*))\leq\frac{G^2}{2\mu}(\log
            T +1).
        \end{eqnarray*}
        \item If Assumption A5(b) holds, then by choosing $\eta_t=\frac{1}{\mu t+L}$, the CBA ensures that
        \begin{eqnarray*} \mathbb{E}(H(\bar
            x_T)-H(x^*))&\leq&\frac{\sigma^2}{2\mu}\frac{\log T
                +1}{T}+\frac{H(x_1)-H(x^*)}{T}+\frac{L(x_1-x^*)^2}{2T}, \mbox{   and}\\
            \sum_{t=1}^{T}\mathbb{E}(H(x_t)-H(x^*))&\leq&\frac{\sigma^2}{2\mu
            }(\log T +1)+H(x_1)-H(x^*)+\frac{L(x_1-x^*)^2}{2}.
        \end{eqnarray*}
    \end{itemize}
\end{proposition}

Proposition \ref{thm:CBA1sc} gives the performance of the CBA
when $\mu>0$. The convergence rate of $O\left(\log T/T\right)$ when
$\mu>0$ is better than the convergence rate of
$O\left(1/\sqrt{T}\right)$ when $\mu=0$. In such case, one does not
need to know $T$ in advance and can always choose the stepsize as a
decreasing sequence in $t$. Again, the detailed proof of the proposition
is given in Appendix \ref{appendix:proof_of_theorem1-3}.

According to Proposition~\ref{thm:CBA1sc}, when $\mu > 0$, the
convergence rate of CBA is $O\left(\log T/T\right)$. In the
following, we improve the convergence rate when $\mu > 0$ to
$O(1/T)$ using a restarting method first proposed by
\citet{Hazan:14}. In \citet{Hazan:14}, the authors show that the
restarting method works when Assumption A5(a) holds. In this paper,
we extend the result by showing that the restarting method can also
obtain the $O(1/T)$ rate if Assumption A5(b) holds. According to
\citet{Nemirovski:83}, no algorithm can achieve convergence rate
better than $O\left(1/T\right)$, thus we have obtained the best
possible convergence rates in those settings. We now describe the
restarting method in Algorithm \ref{alg:rsa}, which we will later
refer to as the multi-stage comparison-based algorithm (MCBA).

\begin{algorithm} \caption{Multi-stage Comparison-Based Algorithm (MCBA)}\label{alg:rsa}
\begin{enumerate}
\item Initialize the number of stages $K\geq 1$, the starting solution $\hat x^1$. Set $k=1$.
\item Let $T_k$ be the number of iterations in stage $k$ and $\eta_t^k$ be the step length in iteration $t$ of CBA  in stage $k$ for $t=1,2,\dots,T_k$.
\item Let $\hat x^{k+1}=\textbf{CBA}(\hat
x^k,T_k,\{\eta_t^k\}_{t=1}^{T_k})$.
\item Stop when $k \geq K$. Otherwise, let $k \leftarrow k + 1$ and go back to step 2.
\item Output $\hat x^{K+1}$.
\end{enumerate}
\end{algorithm}

We have the following proposition about the performance of Algorithm
\ref{alg:rsa}. The proof of the proposition is given in Appendix
\ref{appendix:proof_of_theorem1-3}.

\begin{proposition}
    \label{thm:restart} Suppose $\mu > 0$. Let $G^2$ and $\sigma^2$ be defined
    as in Proposition~\ref{prop:unbiased}, $x^*$ be any optimal
    solution to (\ref{formulation}), and $T=\sum_{k=1}^KT_k$ with $T_k$ defined in MCBA.
    \begin{itemize}
        \item If Assumption A5(a) holds, then by choosing $\eta_t^k=\frac{1}{2^{k+1}\mu}$ and $T_k=2^{k+3}$, the MCBA
        ensures that
        \begin{eqnarray*}
            \mathbb{E}(H(\hat x_{K+1})-H(x^*))
            \leq\frac{16(H(\hat{x}_1) - H(x^*)+G^2/\mu)}{T}.
        \end{eqnarray*}
        \item If Assumption A5(b) holds, then by choosing $\eta_t^k=\frac{1}{2^{k+1}\mu+L}$ and $T_k=2^{k+3}+4$,
        the MCBA ensures that
        \begin{eqnarray*}
            \mathbb{E}(H(\hat x_{K+1})-H(x^*)) \leq\frac{32(H(\hat{x}_1) - H(x^*)+L(\hat{x}_1-x^*)^2/2+\sigma^2/\mu)}{T}.
        \end{eqnarray*}
    \end{itemize}
\end{proposition}

 Both Proposition \ref{thm:CBA1sc} and \ref{thm:restart}
give the convergence rates of CBA and MCBA for strongly convex
problems (i.e., $\mu>0$). The convergence rate of MCBA is $O(1/T)$
which improves the $O\left(\log T/T\right)$ convergence rate of CBA
by a factor of $\log T$. The intuition for this difference is that
the output $\bar x_T$ of CBA is the average of all historical
solutions so its quality is reduced by the earlier solutions (i.e.,
$x_t$ with a small $t$) which are far from the optimal solution. On
the contrary, MCBA restarts CBA periodically with a better initial
solution   (i.e., $\hat x_k$) for each restart. As a result, the
output of MCBA is the average of historical solutions only in the
last ($K$th) call of CBA which does not involve the earlier
solutions and thus has a higher quality. This is the main reason for
MCBA to have a better solution after the same number of iterations,
or equivalently, a better convergence rate than CBA. However, CBA is
easier to implement as it does not require periodic restart as
needed in MCBA. Moreover, the theoretical convergence of MCBA
requires strong convexity in the problem while CBA converges without
strong convexity requirement (see Proposition
\ref{thm:CBA1}).\footnote{ Note that, in
Proposition~\ref{thm:restart}, we only focus on the case when $\mu>0$
since MCBA is mainly designed to improve the convergence rate
$O\left(\log T/T\right)$ of CBA when the problem is strongly convex.
When $\mu=0$, the convergence rate of MCBA is still
$O(1/\sqrt{T})$.}

By Proposition \ref{thm:CBA1}-\ref{thm:restart}, we have shown that
under some mild assumptions (Assumption
\ref{assumption:onedimensional}), if one has access to comparative
information between each sample $\xi_t$ and \emph{two} points, then
one can still find the optimal solution to (\ref{formulation}), and
the convergence speed is in the same order as when one can observe
the actual value of the sample (or the objective value at the
sampled point). One natural question is whether the same convergence
result can be achieved by only having comparative information
between each sample $\xi_t$ and \emph {one} point. The next example
gives a negative answer to this question, showing that it is
impossible to always find the optimal solution in this case, even if
one allows the algorithm to be a randomized one. Thus it verifies
the necessity of having comparative information at two points in
each iteration (for each sample).

\begin{example}\label{ex:onecomp}
Let $h(x, \xi) = (x - \xi)^2$, $\ell = -1$ and $u = 1$. In this
case, the optimization problem (\ref{formulation}) is to find the
projection of the expected value of $\xi$ onto the interval $[-1,
1]$. Suppose there are two underlying distributions for $\xi$. In
the first case, $\xi$ follows a uniform distribution on $[-3, -2]$
or $[2, 3]$, each with probability $0.5$. In the second case, $\xi$
follows a uniform distribution on $[-3, -2]$ or $[3, 4]$, each with
probability $0.5$. In the following, we denote the distributions
corresponding to the first and second cases by $F_1(\cdot)$ and
$F_2(\cdot)$, respectively. It is easy to verify that in the first
case, the optimal solution to (\ref{formulation}) is $x^* = 0$,
while in the second case, the optimal solution to
(\ref{formulation}) is $x^* = 0.25$. And it is also easy to verify
that the above settings (both cases) satisfy Assumption
\ref{assumption:onedimensional}.

Now we consider any algorithm that only utilizes the comparative
information between $\xi_t$ and one decision point $x_t$ in each
iteration (however, the point has to be chosen between $[\ell, u]$
since we can modify $h(x, \xi)$ such that it is undefined or
$\infty$ on $x\notin [\ell, u]$). Suppose the algorithm maps $x_t$
and the comparative information between $\xi_t$ and the chosen decision point
to a distribution of $x_{t+1}$ (thus we allow randomized algorithm).
Note that for any $x\in [\ell, u]$, $F_1(x) = F_2(x)=0.5$. In other
words, there is a $0.5$ probability that $\xi > x$ and a $0.5$
probability that $\xi <x$ no matter whether $\xi$ is drawn from
$F_1(\cdot)$ or $F_2(\cdot)$. For any algorithm, the distribution of
each $x_{t}$ will be the same under either case. Thus, no algorithm
can return the optimal solution in both cases. In other words, no
algorithm can guarantee to solve (\ref{formulation}) with comparative
information between each sample and only one point, even under
Assumption \ref{assumption:onedimensional}. $\hfill\Box$
\end{example}

\section{Choice of $f_{-}$ and $f_{+}$}
\label{sec:choice of f}

In the CBA, one important step is the specification of the two sets
of density functions $f_{-}(x, z)$ and $f_{+}(x, z)$. In the last
section, we only said that $f_{-}$ and $f_{+}$ need to satisfy
conditions (C1)-(C3) but did not give any specific examples. Nor did
we discuss what are good choices of $f_{-}$ and $f_{+}$. In this
section, we address this issue by first showing several examples of
$f_{-}$ and $f_{+}$ which could be useful in practice and then
discussing the effect of choices of $f_{-}$ and $f_{+}$ on the
efficiency of the algorithms. We start with the following examples
of choices of $f_{-}$ and $f_{+}$.

\begin{example}[Uniform Sampling Distribution]
\label{example:funiform} Suppose the support of $\xi$ is known to be
contained in a finite interval $[\underline{s}, \bar{s}]$ and the
optimal decision $x^*$ is known to be within a finite interval
$[\ell, u]$. (Without loss of generality, we assume $[\underline{s},
\bar{s}] \subseteq [\ell, u]$. Otherwise, we can expand $[\ell, u]$
to contain $[\underline{s}, \bar{s}]$.) And we assume
$h''_{x,z}(x,z)$ is uniformly bounded on $[\ell,
u]\times[\underline{s}, \bar{s}]$, $x\neq z$. Then we can set both
$f_{-}(x,z)$ and $f_{+}(x,z)$ to be uniformly distributed, i.e., for
$x \in (\ell, u)$,
\begin{equation*}
f_{-}(x, z) = \begin{cases}
    \frac{1}{x-\ell}  & \ell \leq z < x\\
    0    &  \text{otherwise}
  \end{cases}
  \qquad\mbox{and}\qquad
  f_{+}(x, z) = \begin{cases}
    \frac{1}{u-x} & x < z \leq u  \\
    0             & \text{otherwise}
  \end{cases}.
\end{equation*}
When $x = \ell$, we can set $f_{-}$ to be a uniform distribution on
$[\ell - 1, \ell]$; and when $x = u$, we can set $f_{+}$ to be a
uniform distribution on $[u, u+1]$. It is not hard to verify that
this set of choices satisfy conditions (C1)-(C3).\\



\begin{example}[Exponential Sampling Distribution]
\label{example:fexponential} Suppose the support of $\xi$ is
${\mathbb R}$ or unknown, and $\xi$ follows a light tail
distribution (more precisely, there exists a constant $\bar{\lambda}
> 0$ such that $\lim_{t \rightarrow\infty} e^{\bar{\lambda}{t}}
{\mathbb P} (|\xi| > t) = 0$). Moreover, $x$ is constrained on a
finite interval $[\ell, u]$ and we assume $h''_{x,z}(x,z)$ is
uniformly bounded on $[\ell, u]\times \mathbb{R}$, $x\neq z$. Then we can
choose $f_{-}$ and $f_{+}$ to be exponential distributions. More
precisely, we can choose
\begin{equation*}
f_{-}(x, z) = \begin{cases}
   0                                        & z \ge x \\
    \lambda_{-} \exp(-\lambda_{-} (x-z))    & z < x
  \end{cases} \qquad\mbox{and}\qquad
  f_{+}(x, z) = \begin{cases}
    \lambda_{+} \exp(-\lambda_{+} (z-x)) & z > x  \\
    0                                    & z \le x
  \end{cases}
\end{equation*}
where $0 < \lambda_{-}, \lambda_{+} < \bar\lambda$ are two parameters
one can adjust. Apparently under this choice of $f_{-}$ and $f_{+}$,
conditions (C1)-(C2) are satisfied. For (C3), we note that by the
light tail assumption, there exists a constant $C$ such that
$e^{-\bar{\lambda}t} F(t) \le C$ and $e^{\bar{\lambda}t}(1-F(t)) \le
C$ for all $t$. Therefore, we have
\begin{eqnarray*}
\int_{-\infty}^{x-} \frac{F(z)}{f_{-}(x,z)}dz = \frac{
e^{\lambda_{-}x}}{\lambda_{-}}\int_{-\infty}^{x} e^{-\lambda_{-}z}
F(z) dz \le \frac{Ce^{\lambda_{-}x}}{\lambda_{-}} \int_{-\infty}^x
e^{(\bar{\lambda} -\lambda_{-})z} dz \le
\frac{Ce^{\bar{\lambda}x}}{(\bar{\lambda} -
\lambda_{-})\lambda_{-}},\\
 \int_{x+}^{\infty}\frac{1 - F(z)}{f_{+}(x,z)}dz = \frac{
e^{-\lambda_{+}x}}{\lambda_{+}}\int_{x}^{\infty} e^{\lambda_{+}z} (1
- F(z)) dz \le \frac{Ce^{-\lambda_{+}x}}{\lambda_{+}}
\int_{x}^\infty e^{-(\bar{\lambda}- \lambda_{+})z}dz \le
\frac{Ce^{-\bar{\lambda}x}}{(\bar{\lambda} -
\lambda_{+})\lambda_{+}}.
\end{eqnarray*}
Combined with the uniform boundedness of $h''_{x,z}(x,z)$, condition (C3) also holds in this case.\\
\end{example}

Next we discuss the effect of choosing different $f_{-}$ and $f_{+}$
on the efficiency of the algorithm and the optimal choices of
$f_{-}$ and $f_{+}$. First we note that by Proposition
\ref{thm:CBA1}-\ref{thm:restart}, the choice of $f_{-}$ and $f_{+}$
does not affect the asymptotic convergence rate of the algorithms as
long as they satisfy conditions (C1)-(C3). All what they affect is
the constant in the convergence results, which depends on ${\mathbb
E}_{z, \xi}(g(x,\xi, z))^2$ or ${\mathbb E}_{z,\xi}(g(x,\xi,z) -
H'(x))^2$ (depending on whether Assumption A5(a) or A5(b) holds).
However, by Proposition \ref{prop:unbiased}, ${\mathbb
E}_{z,\xi}(g(x,\xi,z) - H'(x))^2 = {\mathbb E}_{z, \xi}(g(x, \xi,
z))^2 - (H'(x))^2$, i.e., the two terms only differ by a constant
which does not depend on the choice of $f_{-}$ and $f_{+}$.
Therefore, in what follows, we focus on choosing $f_{-}$ and $f_{+}$
to minimize ${\mathbb E}_{z, \xi}(g(x,\xi, z))^2$.

By (\ref{eq:A1}), for any $\xi < x$, we have
\begin{eqnarray*}
{\mathbb E}_{z}(g(x, \xi, z))^2 = 2h'_{-}(x)h'_{x}(x,\xi) -
(h'_{-}(x))^2 +
\int_{\xi}^{x-}\frac{(h''_{x,z}(x,z))^2}{f_{-}(x,z)}dz.
\end{eqnarray*}
By a similar argument, for $\xi > x$, we have
\begin{eqnarray*}
{\mathbb E}_{z}(g(x, \xi, z))^2 = 2h'_{-}(x)h'_{x}(x,\xi) -
(h'_{-}(x))^2 +
\int_{x+}^{\xi}\frac{(h''_{x,z}(x,z))^2}{f_{+}(x,z)}dz.
\end{eqnarray*}

Further taking expectation over $\xi$, we have
\begin{eqnarray*}
 && {\mathbb E}_{z, \xi} (g(x, \xi, z))^2 \\
 &= &  2h'_{-}(x)H'(x) - (h'_{-}(x))^2 + \int_{-\infty}^{x-}\left(
\int_{\xi}^{x-}\frac{(h''_{x,z}(x,z))^2}{f_{-}(x,z)}dz\right)
dF(\xi) + \int_{x+}^{\infty}\left(
\int_{x+}^{\xi}\frac{(h''_{x,z}(x,z))^2}{f_{+}(x,z)}dz\right) dF(\xi)\\
 &= & 2h'_{-}(x)H'(x) - (h'_{-}(x))^2 +
\int_{-\infty}^{x-}\frac{F(z)(h''_{x,z}(x,z))^2}{f_{-}(x,z)} dz +
\int_{x+}^{\infty}\frac{(1 - F(z))(h''_{x,z}(x,z))^2}{f_{+}(x,z)}dz.
\end{eqnarray*}
By Cauchy-Schwarz inequality and condition (C1), we have
\begin{eqnarray*}
\int_{-\infty}^{x-}\frac{F(z)(h''_{x,z}(x,z))^2}{f_{-}(x,z)} dz =
 \int_{-\infty}^{x-}\frac{F(z)(h''_{x,z}(x,z))^2}{f_{-}(x,z)} dz
\cdot \int_{-\infty}^{x-} f_{-}(x,z)dz  \ge \left(\int_{-\infty}^{x}
\sqrt{F(z)}|h''_{x,z}(x,z)|dz\right)^2.
\end{eqnarray*}
And the equality holds only if $f_{-}(x, z) = C_{-}
\sqrt{F(z)}|h''_{x,z}(x,z)|$, for all $z < x$ for some $C_{-} > 0$.
Therefore, if $\sqrt{F(z)}|h''_{x,z}(x,z)|$ is integrable on
$(-\infty, x]$, then the optimal choice for $f_{-}(x,z)$ is
\begin{eqnarray}
\label{eq:optfneg}
f_{-}(x, z) = \frac{\sqrt{F(z)}|h''_{x,z}(x,z)|}{\int_{-\infty}^x
\sqrt{F(z)}|h''_{x,z}(x,z)|dz}, \quad\forall z < x.
\end{eqnarray}
Similarly, if $\sqrt{1 - F(z)}|h''_{x,z}(x,z)|$ is integrable on
$[x, \infty)$, then the optimal choice for $f_{+}(x,z)$ is
\begin{eqnarray}
\label{eq:optfpos}
f_{+}(x, z) = \frac{\sqrt{1 - F(z)}|h''_{x,z}(x,z)|}{\int_x^\infty
\sqrt{1 - F(z)}|h''_{x,z}(x,z)|dz}, \quad\forall z > x.
\end{eqnarray}

Now we use an example to illustrate the above results. Suppose
$h(x,\xi) = (x-\xi)^2$ and $F(z)$ is a uniform distribution on $[a,
b]$. Then the optimal choice of $f_{-}$ and $f_{+}$ are
\begin{eqnarray*}
f_{-}(x, z) =  \frac{3(z-a)^{1/2}}{2(x-a)^{3/2}}, \; \forall a \le z
< x \leq b \quad \mbox{and} \quad f_{+}(x,z)
=\frac{3(b-z)^{1/2}}{2(b-x)^{3/2}}, \; \forall a\leq x < z \le b.
\end{eqnarray*}
Similarly, when $h(x,\xi) = (x-\xi)^2$ and $F(z)$ is a normal
distribution $\mathcal{N}(a, b^2)$, the optimal choice of $f_{-}$
and $f_{+}$ are (it is easy to show that the integrals on the bottom
are finite)
\begin{eqnarray*}
f_{-}(x,z) =
\frac{\sqrt{\Phi\left(\frac{z-a}{b}\right)}}{\int_{-\infty}^{x}\sqrt{\Phi\left(\frac{z-a}{b}\right)}dz}\quad
\forall z < x \quad \mbox{and}\quad f_{+}(x,z) = \frac{\sqrt{1 -
\Phi\left(\frac{z-a}{b}\right)}}{\int_{x}^{\infty}\sqrt{1 -
\Phi\left(\frac{z-a}{b}\right)}dz}\quad\forall z >x
\end{eqnarray*}
where $\Phi(\cdot)$ is the c.d.f. of a standard normal distribution.

However, in many cases, the optimal choice of $f_{-}$ and $f_{+}$
may not exist due to either 1) $\sqrt{F(z)}|h''_{x,z}(x,z)|$ or
$\sqrt{1 - F(z)}|h''_{x,z}(x,z)|$ is not integrable, or 2) the
integration is $0$ (e.g., in the case when $h(x,\xi)$ is piecewise
linear in $x$ and $\xi$). In those cases, either the optimal $f_{-}$
and $f_{+}$ are not attainable, or the choice of $f_{-}$ and $f_{+}$
does not matter (e.g., in the piecewise linear case\footnote{In the
piecewise linear case, the comparative information between
$x$ and $\xi$ will imply the knowledge of the stochastic gradient at
the current sample point, and the gradient does not depend on the
choice of $f_{-}$ and $f_{+}$ functions.}). Moreover, finding the
optimal $f_{-}$ and $f_{+}$ essentially needs the knowledge of the
distribution of $\xi$, which is not known in advance. Therefore, one
can only use an approximate (or prior) distribution of $\xi$ to
calculate the distribution.\footnote{Such issues are common in
variance reduction problems in simulation. The optimal choice to
reduce variance relies on the knowledge of the underlying
distribution. See, e.g., \citet{glynn}.} In addition, sampling from
the distributions described above usually involves much more
computational efforts than sampling from a uniform distribution or
an exponential distribution, the overhead of which may well
overshadow the improvement of the convergence speed. Therefore, in
practice, choosing a heuristic sampling distribution $f_{-}$ and
$f_{+}$ may be more preferable, such as the uniform or exponential
distributions described earlier in this section. Indeed, as we will
see in later in Section \ref{sec:numerical}, using uniform or
exponential distributions lead to efficient solutions in our test
cases.


\end{example}

\section{Extensions}
\label{sec:exts} In this section, we discuss a few extensions of our
comparison-based algorithms. In particular, in Section \ref{sec:QP},
we extend our discussions to multi-dimensional problems with
quadratic objective functions. In Section \ref{sec:NonConvex}, we
consider the case in which the objective function is not a convex
function. In Section \ref{sec:MB}, we consider the case in which
multiple samples can be drawn and multiple comparisons can be
conducted in each iteration.  We also consider a case
in which categorical results (depending on the difference between
the decision point and the sampled point) instead of binary results
can be obtained from each comparison in Appendix
\ref{appendix:multilevelCBA}. Overall, we show that our proposed
ideas can still be applied in those settings (with some variations).

\subsection{Multi-Dimensional Convex Quadratic Problem}\label{sec:QP}
In this section, we extend our model to a multi-dimensional convex
quadratic problem and propose a stochastic optimization algorithm
for such a setting based on comparative information.
Specifically, we consider the following stochastic convex
optimization problem
\begin{eqnarray}
\label{eq:QPobj} \min_{x\in {\mathcal X}} H(x) = {\mathbb E}_{\xi}
\left[h(x, \xi):=\frac{1}{2}(x-\xi)^\top Q(x-\xi)\right],
\end{eqnarray}
where $x\in\mathbb{R}^d$, $\xi\in\mathbb{R}^d$ is a random variable,
$Q$ is a positive definite matrix and $\mathcal X$ is a closed
convex set in $\mathbb{R}^d$.

    We denote the gradient of $H$ by $\nabla H(x):=\mathbb{E}Q(x-\xi)$,
    the directional derivative of $h(x ,\xi )$ along a direction
    $u\in\mathbb{R}^d$ by $\nabla_u h(x ,\xi ):=u^{T}Q(x-\xi)$, and the gradient of $h(x,\xi)$ with respect to $x$ by $\nabla h(x ,\xi ):=Q(x-\xi)$.
    The following assumption is made in this section:
    \begin{assumption}\label{assumption:qp}
        There exists a constant $K_4$ such that
        $\mathbb{E}\|\xi-\mathbb{E}\xi\|_2^2\leq K_4^2$.
    \end{assumption}
    This assumption simply requires that $\xi$ has a finite variance,
    which is not hard to satisfy in practice.

    Suppose we can generate a random vector $u$ in $\mathbb{R}^d$ that
satisfies $\mathbb{E}(uu^\top)=I_{d}$, where $I_{d}$ is the $d\times d$ identity
matrix.  We will have $\mathbb{E}_{u}u \nabla_u h(x ,\xi
)=\mathbb{E}_{u}uu^\top Q(x-\xi)=\nabla h(x ,\xi )$. Hence, to construct an unbiased stochastic gradient for $H(x)$, we only need to
    construct an unbiased stochastic estimation for $u^\top Q(x-\xi)$ and multiply it to $u$. In the following,
    we show that this can be done by first comparing $x+zu$ and $x-zu$ (in the value of $h(\cdot,\xi)$) with a random positive number $z$
    and then comparing the better one between $x+zu$ and $x-zu$ with $x$.  In other words, we can
    still construct a stochastic gradient for $H(x)$ in \eqref{eq:QPobj} using two comparisons.

    Similar as Example~\ref{example:mean}, this problem may still represent a problem of finding the average features of a group of customers. In particular, in such problems, $x$ may represent $d$ features (e.g., size, taste, etc.) of a product
    while $\xi$ is the preference of a random customer on those features. The firm would like to find
    the optimal features to minimize the expected customer
    dissatisfaction which is measured by $h(x,\xi)$ in \eqref{eq:QPobj}. Suppose the current product is $x$. To implement the two comparisons above,
    the firm can generate a random change $u$ and a random level $z$ and then ask a customer for his/her preference between two new products $x+zu$ and $x-zu$ and his/her preference between the better new product and the current product.

    We call this method comparison-based algorithm for quadratic problem (CBA-QP) and provide its details in  Algorithm~\ref{alg:cbaqp}.
    In CBA-QP, we need to specify a density function $f(z)$ on
$[0,+\infty)$, which needs to satisfy the following condition.

\begin{itemize}
    \item (C4) There exists a constant $K_5$ such that for any $x\in {\mathcal X}$,
    $$\int_{0}^{+\infty}\frac{\text{Prob}\left(\frac{\lambda_{\max}(Q)\|x-\xi\|_2}{\lambda_{\min}(Q)\sqrt{d}}\geq z\right)}{f(z)} dz=\int\int_{0}^{\frac{\lambda_{\max}(Q)\|x-\xi\|_2}{\lambda_{\min}(Q)\sqrt{d}}}\frac{1}{f(z)} dzdF(\xi) \leq K_5$$
    where $\lambda_{\max}(Q)$ and $\lambda_{\min}(Q)$ are the largest
    and smallest eigenvalues of $Q$, respectively.
\end{itemize}

Below we give two examples of choices of $f(z)$ such that it
satisfies condition (C4).

\begin{example}[Uniform Sampling Distribution]
\label{example:funiform2} Suppose $\Xi$ (the support of $\xi$) and the feasible set $\mathcal{X}$ are both bounded. The quantity
$R:=\frac{\lambda_{\max}(Q)}{\lambda_{\min}(Q)\sqrt{d}}\max\limits_{x\in\mathcal{X},\xi\in\Xi}\|x-\xi\|_2$
is finite. We can set  $f(z)$ to be the density function of
a uniform distribution on $[0, R]$, i.e.,
        \begin{equation*}
        f(z) = \begin{cases}
       1/R                                   & 0 \leq z \leq R\\
        0    &  \text{otherwise.}
        \end{cases}
        \end{equation*}
        With this choice, we have
        $$
        \int\int_{0}^{\frac{\lambda_{\max}(Q)\|x-\xi\|_2}{\lambda_{\min}(Q)\sqrt{d}}}
        \frac{1}{f(z)} dzdF(\xi) \leq R^2,
        $$
        thus (C4) is satisfied with $K_5=R^2$.\\
    \end{example}
    \begin{example}[Exponential Sampling Distribution]
        \label{example:fexponential2} Suppose $\xi$ follows a light tail
        distribution with $\mathbb{E}\exp(\|\xi\|_2)\leq
        \bar{\sigma}$ for some $\bar{\sigma} > 0$.
Moreover, suppose the
        feasible set $\mathcal{X}$ is bounded. Then we can
        choose $f(z)$ to be an exponential distribution, i.e.,
        \begin{equation}
        \label{eq:fzexpdist}
        f(z) = \begin{cases}
        0                                        & z < 0 \\
        \lambda \exp(-\lambda z)    & z \geq 0,
        \end{cases}
        \end{equation}
        where $\lambda$ can be any positive constant less than $c:=\frac{\sqrt{d}\lambda_{\min}(Q)}{\lambda_{\max}(Q)}$.
        With this choice, we have
        \small
        $$
        \int_{0}^{\frac{\lambda_{\max}(Q)\|x-\xi\|_2}{\lambda_{\min}(Q)\sqrt{d}}}
        \frac{1}{f(z)} dz
        \leq
        \frac{1}{\lambda^2} \exp\left( \frac{\lambda_{\max}(Q)\|x-\xi\|_2\lambda}{\lambda_{\min}(Q)\sqrt{d}}\right)
        \leq \frac{1}{\lambda^2} \exp\left( \frac{\lambda\max\limits_{x\in\mathcal{X}}\|x\|_2}{c}\right)
        \exp\left( \frac{\lambda\|\xi\|_2}{c}\right).
        $$
        \normalsize
        By the light tail assumption and the concavity of $(\cdot)^{\lambda/c}$, we can use Jensen's inequality to show that $\mathbb{E}\exp\left( \frac{\lambda\|\xi\|_2}{c}\right)    \leq \left[\mathbb{E}\exp\left( \|\xi\|_2\right)\right]^\frac{\lambda}{c}\leq \bar\sigma^\frac{\lambda}{c}  $. Using this inequality, we
        further have
        \small$$
        \int\int_{0}^{\frac{\lambda_{\max}(Q)\|x-\xi\|_2}{\lambda_{\min}(Q)\sqrt{d}}}
        \frac{1}{f(z)} dzdF(\xi)\leq  \frac{\bar\sigma^\frac{\lambda}{c}}{\lambda^2} \exp\left(
        \frac{\lambda\max\limits_{x\in\mathcal{X}}\|x\|_2}{c}\right),
        $$
        \normalsize
        thus (C4) is satisfied with $K_5= \frac{\bar\sigma^\frac{\lambda}{c}}{\lambda^2} \exp\left(
        \frac{\lambda}{c}\max\limits_{x\in\mathcal{X}}\|x\|_2\right)$.\\
        \end{example}

    \begin{algorithm}[h] \caption{CBA for Quadratic Problem (CBA-QP):\label{alg:cbaqp}}
        \begin{enumerate}
            \item {\bf Initialization.} Set $t = 1$, $x_1 \in {\mathcal X}\subset\mathbb{R}^d$. Define $\eta_t$
            for all $t\ge 1$. Set the maximum number of iterations $T$.
            Choose a density function $f(z)$ on $[0,+\infty)$.  Let $\mathcal{Q}$ be the uniform distribution on a sphere in $\mathbb{R}^d$ with radius
            $\sqrt{d}$. Note that $\mathbb{E}uu^T=I_d$ for $u$ following distribution $\mathcal{Q}$.
            \item {\bf Main iteration.} Sample $\xi_t$ from the distribution of $\xi$. Sample $u_t$ from $\mathcal{Q}$.
            Sample $z_t$ from $f(z)$
            \begin{enumerate}
                \item If $h(x_t+z_tu_t, \xi_t) < h(x_t-z_tu_t, \xi_t)$, set
                \begin{eqnarray*}\label{gradient_case1qp}
                    g(x_t, \xi_t, u_t,z_t) = \left\{\begin{array}{ll} 0, & \mbox{
                            if } h(x_t+z_tu_t, \xi_t) > h(x_t, \xi_t),
                        \\
                        -\frac{1}{2}\frac{u_t^TQu_t}{f(z_t)}u_t, & \mbox{ if
                        } h(x_t+z_tu_t, \xi_t) \leq h(x_t, \xi_t). \end{array}\right.
                \end{eqnarray*}
                \item If $h(x_t+z_tu_t, \xi_t) \geq h(x_t-z_tu_t, \xi_t)$, set
                \begin{eqnarray*}\label{gradient_case2qp}
                    g(x_t, \xi_t, u_t,z_t) = \left\{\begin{array}{ll} 0, & \mbox{
                            if } h(x_t-z_tu_t, \xi_t) > h(x_t, \xi_t),
                        \\
                        \frac{1}{2}\frac{u_t^TQu_t}{f(z_t)}u_t, & \mbox{ if
                        } h(x_t-z_tu_t, \xi_t) \leq h(x_t, \xi_t). \end{array}\right.
                \end{eqnarray*}
            \end{enumerate}
            Let
            \begin{eqnarray} \label{eq:projqp}
            x_{t+1}&=&\text{Proj}_{\mathcal X}\left(x_t-\eta_tg(x_t, \xi_t, u_t
            z_t)\right) .
            \end{eqnarray}
            \item {\bf Termination.} Stop when $t\ge T$. Otherwise, let $t\leftarrow t+ 1$ and go
            back to Step 2.
            \item {\bf Output.} $\textbf{CBA-QP}(x_1,T,\{\eta_t\}_{t=1}^T) = \bar x_T=\frac{1}{T}\sum_{t=1}^Tx_t$.
        \end{enumerate}
    \end{algorithm}

    The proposition below gives some properties of the stochastic gradient $g(x_t, \xi_t, u_t,z_t)$ in CBA-QP (see Algorithm \ref{alg:cbaqp}):
    \begin{proposition}\label{prop:unbiased_multi_dimensional}
       Let $z$, $u$, and $g$ be defined in Algorithm \ref{alg:cbaqp}. Then the following properties hold.
        \begin{enumerate}
            \item
            ${\mathbb E}_{z,u}g(x , \xi , u, z ) = \nabla h(x ,\xi )$, for all $x\in
            {\mathcal X}$.
            \item ${\mathbb E}_{z ,u, \xi } g(x , \xi , u, z ) =
            \nabla H(x )$, for all $x \in  {\mathcal X}$.
            \item ${\mathbb E}_{z,u,\xi} \|g(x , \xi , u, z )-\nabla H(x)\|_2^2\leq \sigma^2:=\lambda_{\max}(Q)^2 K_4^2+\frac{\lambda_{\max}(Q)^2d^3K_5}{4}$.
        \end{enumerate}
    \end{proposition}

    {\bf\noindent Proof of Proposition
        \ref{prop:unbiased_multi_dimensional}.} First, we consider the case
    when $h(x+zu, \xi) < h(x-zu, \xi)$, which is equivalent to
    $u^TQ(x-\xi)<0$ by the definition of $h$ and the non-negativity of
    $z$. Note that $h(x+zu, \xi) \leq h(x, \xi)$ if and only if
    $zu^TQ(x-\xi)+\frac{z^2}{2}u^TQu\leq0$, or equivalently, $0\leq z\leq
    -\frac{2u^TQ(x-\xi)}{u^TQu}$. It then follows from the definition of
    $g$ that
    \begin{eqnarray*}
        {\mathbb E}_{z } g(x , \xi , u, z ) = \int_{0}^{-\frac{2u^TQ(x-\xi)}{u^TQu}}
        -\frac{1}{2}u^TQu \cdot udz = u^TQ(x-\xi)u=\nabla_u h(x ,\xi )u.
    \end{eqnarray*}
    Similarly, the second case $h(x+zu, \xi) \geq h(x-zu, \xi)$ occurs when $u^TQ(x-\xi)\geq 0$.  The inequality $h(x-zu, \xi) \leq h(x, \xi)$ holds if and only if $-zu^TQ(x-\xi)+\frac{z^2}{2}u^TQu\leq0$, or equivalently,  $0\leq z\leq \frac{2u^TQ(x-\xi)}{u^TQu}$.
    It then follows again from the definition of $g$ that
    \begin{eqnarray*}
        {\mathbb E}_{z } g(x , \xi , u, z ) = \int_{0}^{\frac{2u^TQ(x-\xi)}{u^TQu}}
        \frac{1}{2}u^TQu \cdot udz = u^TQ(x-\xi)u=\nabla_u h(x ,\xi )u.
    \end{eqnarray*}
    Therefore, in both cases, we have ${\mathbb E}_{z } g(x , \xi , u, z
    )=\nabla_u h(x ,\xi )u$.  Since $\mathbb{E}(uu^T)=I_{d}$, further
    taking expectation over $u$ on both sides of this equality gives the
    first conclusion of the proposition.

    The second conclusion can be obtained by taking expectation over
    $\xi$ on both sides of the first conclusion, namely, ${\mathbb
        E}_{z,u, \xi } g(x , \xi , u, z ) ={\mathbb E}_{ \xi }{\mathbb
        E}_{z,u} g(x , \xi , u, z) ={\mathbb E}_{\xi }\nabla h(x ,\xi )=
    \nabla H(x)$. (This holds because $h$ is a convex function, thus one
    can apply the monotone convergence theorem.)

    Next, we prove the third conclusion. The first conclusion implies that
        \begin{eqnarray*}
        {\mathbb E}_{u,z} \|g(x, \xi, u,z)-\nabla h(x ,\xi )\|^2 &=&{\mathbb E}_{u,z}
        \|g(x, \xi, u,z)\|^2- \|\nabla h(x ,\xi )\|^2.
    \end{eqnarray*}

    If $h(x+zu, \xi) < h(x-zu, \xi)$, by the definition of $g$, we can show that
\begin{eqnarray*}
      {\mathbb E}_{z}
        \|g(x, \xi, u,z)\|^2  &=&
        \int_{0}^{-\frac{2u^TQ(x-\xi)}{u^TQu}}
        \frac{(u^TQu)^2}{4f(z)} \|u\|_2^2dz\\
        &\leq&\int_{0}^{\frac{\lambda_{\max}(Q)\|u\|_2\|x-\xi\|_2}{\lambda_{\min}(Q)\|u\|_2^2}}
        \frac{\lambda_{\max}(Q)^2\|u\|_2^4}{4f(z)} \|u\|_2^2dz\\
        &=&\frac{\lambda_{\max}(Q)^2d^3}{4}\int_{0}^{\frac{\lambda_{\max}(Q)\|x-\xi\|_2}{\lambda_{\min}(Q)\sqrt{d}}}
        \frac{1}{f(z)} dz,
    \end{eqnarray*}
where the inequality is by Cauchy-Schwarz inequality and the second equality is because $\|u\|=\sqrt{d}$.
 Similarly, if $h(x+zu, \xi) \geq h(x-zu, \xi)$, then we can also show that
    \begin{eqnarray*}
        {\mathbb E}_{z} \|g(x, \xi, u,z)-\nabla h(x ,\xi )\|^2
        \leq
        \frac{\lambda_{\max}(Q)^2d^3}{4}\int_{0}^{\frac{\lambda_{\max}(Q)\|x-\xi\|_2}{\lambda_{\min}(Q)\sqrt{d}}}
        \frac{1}{f(z)} dz.
    \end{eqnarray*}
As a result, we have $${\mathbb E}_{u,z} \|g(x, \xi, u,z)-\nabla h(x ,\xi )\|^2
\le \frac{\lambda_{\max}(Q)^2d^3}{4}\int_{0}^{\frac{\lambda_{\max}(Q)\|x-\xi\|_2}{\lambda_{\min}(Q)\sqrt{d}}}
\frac{1}{f(z)} dz.$$

    According to this inequality and condition (C4), we have
    \begin{eqnarray*}
        {\mathbb E}_{z, u, \xi} \|g(x, \xi, u, z)- \nabla h(x ,\xi )\|^2 &\le& \frac{\lambda_{\max}(Q)^2d^3}{4}\int\int_{0}^{\frac{\lambda_{\max}(Q)\|x-\xi\|_2}{\lambda_{\min}(Q)\sqrt{d}}}
        \frac{1}{f(z)} dzdF(\xi) \leq \frac{\lambda_{\max}(Q)^2d^3K_5}{4}.
    \end{eqnarray*}
    In addition, by Assumption~\ref{assumption:qp}, we have
    $$
    \mathbb{E}_{\xi}\|\nabla h(x,\xi)-\nabla H(x)\|_2^2= \mathbb{E}_{\xi}\|Q(x-\xi)-Q(x-\mathbb{E}\xi)\|_2^2\leq \lambda_{\max}(Q)^2\mathbb{E}_{\xi}\|\xi-\mathbb{E}\xi\|_2^2\leq \lambda_{\max}(Q)^2K_4^2.
    $$
    Finally, we note that
    \begin{eqnarray*}
        {\mathbb E}_{z,u,\xi} \|g(x, \xi, u, z)-\nabla H(x)\|_2^2&=&{\mathbb E}_{z,u,\xi}
        \|g(x, \xi,u, z)- \nabla h(x ,\xi )\|_2^2+
        \mathbb{E}_{\xi}\|\nabla h(x,\xi)-\nabla H(x)\|_2^2.
    \end{eqnarray*}
    Therefore,
    ${\mathbb E}_{z,u,\xi} \|g(x , \xi , u, z )-\nabla H(x)\|_2^2\leq
    \sigma^2:= \lambda_{\max}(Q)^2 K_4^2+\frac{\lambda_{\max}(Q)^2d^3K_5}{4}$. Thus the proposition holds.
    $\hfill\Box$\\

    Based on Proposition \ref{prop:unbiased_multi_dimensional}, we have
    the following proposition about the performance of CBA-QP.
    \begin{proposition}\label{thm:CBA1-QP}
        Let $\sigma^2$ be defined as in
        Proposition~\ref{prop:unbiased_multi_dimensional}. Let $x^*$ be any
        optimal solution to \eqref{eq:QPobj} and $\mu > 0 $ and $L > 0$ be the smallest
        and largest eigenvalues of $Q$, respectively. Then by choosing $\eta_t=\frac{1}{\mu t+L}$, the CBA-QP ensures that
            \begin{eqnarray*} \mathbb{E}(H(\bar
                x_T)-H(x^*))&\leq&\frac{\sigma^2}{2\mu}\frac{\log T
                    +1}{T}+\frac{H(x_1)-H(x^*)}{T}+\frac{L\|x_1-x^*\|_2^2}{2T}, \mbox{   and}\\
                \sum_{t=1}^{T}\mathbb{E}(H(x_t)-H(x^*))&\leq&\frac{\sigma^2}{2\mu
                }(\log T +1)+H(x_1)-H(x^*)+\frac{L\|x_1-x^*\|_2^2}{2}.
            \end{eqnarray*}
    \end{proposition}

    Similar as to CBA, we can further improve the
    theoretical performance of CBA-QP using a
    restarting strategy as described in MCBA. We denote the resulting restarted algorithm by MCBA-QP.

    \begin{algorithm} \caption{Multi-stage Comparison-Based Algorithm for Quadratic Problem (MCBA-QP)}\label{alg:rsaqp}
        \begin{enumerate}
            \item Initialize the number of stages $K\geq 1$, the starting solution $\hat x^1$. Set $k=1$.
            \item Let $T_k$ be the number of iterations in stage $k$ and $\eta_t^k$ be the step length in iteration $t$ of CBA-QP  in stage $k$ for $t=1,2,\dots,T_k$.
            \item Let $\hat x^{k+1}=\textbf{CBA-QP}(\hat
            x^k,T_k,\{\eta_t^k\}_{t=1}^{T_k})$.
            \item Stop when $k \geq K$. Otherwise, let $k \leftarrow k + 1$ and go back to step 2.
            \item Output $\hat x^{K+1}$.
        \end{enumerate}
    \end{algorithm}

    \begin{proposition}
        \label{thm:restart-QP} Let $\sigma^2$ be defined
        as in Proposition~\ref{prop:unbiased_multi_dimensional}, $\mu > 0$ and $L > 0$ be defined as in Proposition~\ref{thm:CBA1-QP}, and
        $T=\sum_{k=1}^KT_k$ with $T_k$ defined in MCBA-QP.
        By choosing $\eta_t^k=\frac{1}{2^{k+1}\mu+L}$ and $T_k=2^{k+3}+4$,
        the MCBA-QP ensures that
        \begin{eqnarray*}
            \mathbb{E}(H(\hat x_{K+1})-H(x^*)) \leq\frac{32(H(\hat{x}_1) - H(x^*)+L\|\hat{x}_1-x^*\|_2^2/2+\sigma^2/\mu)}{T}.
        \end{eqnarray*}
    \end{proposition}

  The proofs of Proposition \ref{thm:CBA1-QP} and
    \ref{thm:restart-QP} are very similar to those of
    the second part of Proposition \ref{thm:CBA1sc} and the second part of Proposition \ref{thm:restart}, respectively, and are provided in Appendix \ref{appendix:proofoftheorem4-5}.

Although we mainly focus on the quadratic problem for the
multi-dimensional case, it is easy to see that our approach can also
be applied to the multi-dimensional cases where the objective
function is separable with respect to each dimension while the
constraint may not be separable. In particular, our method can be
applied to the following problem
$$
\min\limits_{x\in\mathcal{X}}H(x) = {\mathbb E}_{\xi} \left[h(x, \xi):=\sum_{i=1}^d h_i(x_i,\xi_i)\right]
$$
where $x=(x_1,\dots,x_d)^\top\in\mathbb{R}^d$, $\xi=(\xi_1,\dots,\xi_d)^\top\in\mathbb{R}^d$,  $\mathcal{X}$ is a convex closed set in $\mathbb{R}^d$, and
$h_i(x_i,\xi_i)$ is a function satisfying Assumption~\ref{assumption:onedimensional}.
In such cases, one can apply the idea of CBA to construct an unbiased stochastic gradient of each $h_i(x_i,\xi_i)$ based on comparisons and then concatenate them into an unbiased stochastic gradient of $h(x,\xi)$ in order to apply the projected gradient step~\eqref{eq:projqp}. The resulting algorithm will have the same convergence rates as CBA under each setting in Proposition~\ref{thm:CBA1}-\ref{thm:restart}. Note that, we do not assume $\mathcal{X}$ is separable so we still cannot solve the problem above as $d$ independent problems.

\subsection{Non-Convex Problem}
\label{sec:NonConvex}The focus of our study in this
    paper is to construct an unbiased stochastic gradient based on
    comparative information. Once the construction is done, one can also
    apply the stochastic gradient to non-convex stochastic optimization
    problems. Although the stochastic
    gradient method can no longer guarantee to reach a global optimal solution for a non-convex problem,
    the recent result by \cite{davis2018stochastic} shows that the
    iterative solution generated by stochastic gradient method still
    converges to a nearly stationary point under some conditions.

    In this section, we still consider one-dimensional problem (\ref{formulation}) but
    $H(x)$ is no longer convex. More specifically, we still assume Assumption~\ref{assumption:onedimensional}
    holds except that Assumption (A4) is replaced by
    \begin{itemize}
        \item[]
    (A4') $H(x)$ in \eqref{formulation} is differentiable and $\rho$-\emph{weakly convex} on $[\ell,u]$ for some $\rho > 0$, namely,
    \begin{eqnarray}
    \label{eq:wkconv} H(x_2) \ge H(x_1) + H'(x_1)(x_2-x_1) -
    \frac{\rho}{2} (x_2-x_1)^2, \quad \forall x_1, x_2 \in [\ell, u].
    \end{eqnarray}
    Moreover, $\mathbb{E}_{\xi}(h'_x(x, \xi))=H'(x)$ for all
    $x\in[\ell,u]$.
    \end{itemize}

    For any $\lambda>0$, the \emph{Moreau
        envelope} for \eqref{formulation} is defined as a function
    \begin{eqnarray}
    \label{eq:moreau}
    H_\lambda(x):=\min_{\ell\le y\le u}\left\{H(y)+\frac{1}{2\lambda}(x-y)^2\right\}.
    \end{eqnarray}
    By definition, as long as $\lambda<\frac{1}{\rho}$, the minimization
    problem \eqref{eq:moreau} has a strongly convex objective function
    and has a unique solution, denoted by $\hat x$. Moreover, the
    function $H_\lambda(x)$ is continuously differentiable with the
    gradient given by $ H'_\lambda(x):=\frac{1}{\lambda}(x-\hat x)$.
    According to \cite{davis2018stochastic}, the value of $|H'_\lambda(x)|$ measures the near stationarity of a solution
    $x\in[\ell, u]$ because
    $$
    |x-\hat x|=\lambda|H'_\lambda(x)|\quad\text{  and  }\quad\text{dist}(0;\partial H(\hat x))\leq |H'_\lambda(x)|
    $$
    where $\partial H(\hat x)$ the subdifferential of $H$ at $x$, i.e., the set of all $v$ satisfying $H(y)\geq H(x)+v(y-x)+o(\|y-x\|)$ as $y\rightarrow x$ and $\text{dist}(x; A)$ is the nearest distance between point $x$ and set $A$ defined as $\text{dist}(x; A) = \inf_{y\in A}||x-y||$.
    This means that if
    $|H'_\lambda(x)|=\frac{1}{\lambda}|x-\hat x|$ is small, then the solution $x$ is closed
    to a point $\hat x$ (because of small $|x-\hat x|$) which is nearly
    stationary (because of small $\text{dist}(0;\partial H(\hat x))$).

    According to \cite{davis2018stochastic}, in order to find a solution
    $x$ with a small $|H'_\lambda(x)|$, we can still use
    CBA except that the output will be $x_{t^*}$
    where $t^*$ is a random index such that ${\mathbb
        P}(t^*=t)=\frac{\eta_t}{\sum_{s=1}^T\eta_s}$ for $t=1,2,\dots,T$.
    Then, we have $\mathbb{E}| H'_{\lambda}(x_{t^*})|$ with $\lambda=\frac{1}{2\rho}$
    converges to zero in a rate of $O(\frac{1}{\sqrt{T}})$. We state
    this result formally as follows.
    \begin{proposition}[Corollary 2.2 by \citealt{davis2018stochastic}]
        Suppose Assumption (A4) is replaced by Assumption (A4') and
         the optimal value of
        \eqref{formulation}  is finite. Let $G^2$ be defined as in
        Proposition~\ref{prop:unbiased} and $t^*$ be a random index such
        that ${\mathbb P}(t^*=t)=\frac{\eta_t}{\sum_{s=1}^T\eta_s}$ for
        $t=1,2,\dots,T$. By choosing $\eta_t=\frac{1}{\sqrt{T}}$, the CBA
        ensures that
        \begin{eqnarray*}
            \mathbb{E}|H'_{\frac{1}{2\rho}}(x_{t^*})|^2\leq 2\frac{(
                H_{\frac{1}{2\rho}}(x_1)-\min_{\ell\le x\le u}  H(x) )+\rho
                G^2}{\sqrt{T}}.
        \end{eqnarray*}
    \end{proposition}
    Therefore, the comparison-based algorithm can be partly applied to
    non-convex problems.

\begin{remark}
Recently, there has been growing interest on the convergence of
first-order methods when the objective function is non-convex but
satisfies the Kurdyka-{\L}ojasiewicz (KL) property
~\citep{bolte2007lojasiewicz,bolte2014proximal,attouch2010proximal,karimi2016linear,noll2014convergence,xu2017globally,li2017convergence}.
A function $H(x)$ satisfies the KL property at a point $x^*$ if
there exists $\eta>0$, a neighborhood $U$ of $x^*$, and a concave
function $\kappa:[0,\eta]\rightarrow [0,+\infty)$ such that: (i)
$\kappa(0)=0$, (ii) $\kappa$ is continuously differentiable on
$(0,\eta)$, (iii) $\kappa'(\cdot)>0$ on $(0,\eta)$, and (iv)
    $$
    \kappa'(H(x)-H(x^*))\cdot \text{dist}(0;\partial H(x^*))\geq 1
    $$
for any $x\in U$ that satisfies $H(x^*)<H(x)<H(x^*)+\eta$. Existing
results show that, if $H(x)$ is non-convex but satisfies the KL
property, each bounded sequence generated by a first-order method
converges to a stationary point $x^*$ of $H(x)$, and the local convergence
rate to $x^*$ can be characterized by the geometry of  $\kappa$ near
$x^*$ (see, e.g.,  Proposition 3 in \citealt{attouch2010proximal}).

However, almost all existing studies on first-order methods under
the KL property focus on deterministic cases. Since our approaches
construct stochastic gradients, most of the existing results
utilizing KL property cannot be directly applied to our problem. The
only result we are aware of about stochastic gradient descent under
the KL property is given by \citet{karimi2016linear}. However, they
make stronger assumptions than the KL property on the problem. In
particular, they require that the problem be unconstrained, that the
aforementioned neighborhood $U$ be the entire space $\mathbb{R}$,
and that $\kappa$ be in the form of $\kappa(t)=c\sqrt{t}$ for some
$c>0$. If these conditions are satisfied by $H(x)$,  the CBA can
also find an $\epsilon$-optimal solution within $O(1/\epsilon)$
iterations according to Theorem 4 in \citet{karimi2016linear}.
This is because CBA is essentially a
stochastic gradient descent method except that the stochastic
gradient is constructed using comparative information.
\end{remark}

\subsection{Mini-Batch Method with Additional Comparisons}
\label{sec:MB}  In this section, we consider the
scenario where multiple comparisons can be conducted in each
iteration. In such cases, a mini-batch technique can be implemented
in CBA and MCBA to reduce the noise
in stochastic gradient and improve the performance of the
algorithms. In particular, we still compare $\xi_t$ and $x_t$ in
iteration $t$ in Algorithm~\ref{alg:cba1}. In case (a) where
$\xi_t<x_t$, we generate $S$ independent samples, denoted by $z_t^s$
for $s=1,2,\dots,S$, from a distribution on $(-\infty, x_t] $ with
p.d.f. $f_{-}(x_t,z)$ and construct a stochastic gradient $g(x_t,
\xi_t,z_t^s)$ as in \eqref{gradient_case1} with $z_t$ replaced by
$z_t^s$. Similarly, in case (b) where $\xi_t>x_t$,  we generate
$z_t^s$ from a distribution on $[x_t,+\infty) $ with p.d.f.
$f_{+}(x_t,z)$ and construct $g(x_t, \xi_t,z_t^s)$ as in
\eqref{gradient_case2} with $z_t$ replaced by $z_t^s$. After
obtaining $g(x_t, \xi_t,z_t^s)$ in either case, we replace
\eqref{eq:proj} in CBA with the following two
steps
\begin{eqnarray*} \label{eq:projbatch}
\bar g_t&=&\frac{1}{S}\sum_{s=1}^Sg(x_t, \xi_t,z_t^s)\\
x_{t+1}&=&\text{Proj}_{[\ell, u]}\left(x_t-\eta_t\bar g_t\right) = \max\left(\ell, \min\left(u, x_t - \eta_t\bar g_t\right)\right).
\end{eqnarray*}
Here, $\bar g_t$ is the average gradient constructed by a mini-batch, which satisfies $\mathbb{E}_{z_t^s,s=1,\dots,S,\xi_t}\bar g_t=H'(x_t)$ by conclusion 2 in Proposition~\ref{prop:unbiased} and has a smaller noise than  $g(x_t, \xi_t, z_t)$. MCBA can also benefit from this technique by calling CBA after the aforementioned modification.

This mini-batch technique will not improve the asymptotic
convergence rate of CBA and MCBA because it will not completely
eliminate  the noise in the stochastic gradient. However, by
reducing the noise, it will improve the algorithm's performance in
practice as we demonstrate in
Section~\ref{subsec:numerical_minibatch}.\footnote{Note that in the
mini-batch method, the multiple samples are drawn simultaneously. It
is worth noting that it is also possible to draw multiple samples
sequentially, each one depending on the results of all previous
ones. However, that mechanism will be quite complex. Moreover, the
asymptotic performance will not improve because it will not surpass
the asymptotic performance when the samples can be directly observed
(which is already achieved by our current algorithm with two
comparisons). Therefore, we here only choose to present the
mini-batch approach.}


\section{Numerical Tests}
\label{sec:numerical}

In this section, we conduct numerical experiments to show that
although much less information is used, the proposed algorithms
based only on comparative information converge at the same
rate as the stochastic gradient methods. We will also investigate
the impact of choices of $f_{-}$ and $f_{+}$  and the impact of
mini-batch on the performances of the proposed methods.  We
implemented all algorithms in MATLAB running on a 64-bit Microsoft
Windows 10 machine with a 2.70 Ghz Intel Core i7-6820HQ CPU and 8GB
of memory.

\subsection{Convergence of Objective Value}
\label{subsec:numerical:convergence}

In this section, we conduct numerical experiments to test the
performance of the CBA and the MCBA. We consider two objective
functions:
\begin{eqnarray*}
\mbox{ 1)}\quad h_1(x,\xi) = (x - \xi)^2 \qquad\mbox{   and   }
\qquad\mbox{ 2)}\quad h_2(x,\xi) = \left\{\begin{array}{ll} (x - \xi)^2 + (x-\xi) & \mbox{if  } \xi < x\\
2(x-\xi)^2 + 2(\xi- x) & \mbox{if  } \xi\ge x\end{array} \right..
\end{eqnarray*}
Here the two objective functions correspond to the two examples we
described in the beginning with $h_1$ being smooth but not $h_2$.
For each choice of $h(x,\xi)$, we consider two distributions of
$\xi$, a uniform distribution $\mathcal{U}[50,150]$ and a normal
distribution $\mathcal{N}(100,100)$. Thus we have four cases in
total. It is easy to see that for the first objective function, the
optimal solution under either underlying distribution is $x^* =
100$. For the second objective function, the optimal solution under
the uniform distribution is $x^* = 108.66$ while the optimal
solution under the normal distribution is $x^* = 102.82$. In all
experiments, we choose the feasible set to be $[\ell,u]=[50,150]$.
For the cases where $\xi\sim\mathcal{U}[50,150]$, we choose $f_{-}$
and $f_{+}$ to be uniform distributions as in Example
\ref{example:funiform} with $\ell = 50$ and $u = 150$. For the cases
where $\xi\sim\mathcal{N}(100,100)$, we choose $f_{-}$ and $f_{+}$
to be exponential distributions as in Example
\ref{example:fexponential} with
$\lambda_{+}=\lambda_{-}=2^{-4}=0.0625$. (Later in Section
\ref{subsec:numerical_choice_f} we will test the effect of choosing
different $f_{-}$ and $f_{+}$ on the convergence speed of the
algorithms.)

In each of the four settings above, Assumption
\ref{assumption:onedimensional} holds with both A5(a) and A5(b)
satisfied, and $H(x)$ is strongly convex. To compare the performance
of the CBA with or without utilizing the strong convexity of the
problem, we test both step lengths $\eta_t=1/\sqrt{t}$ and $\eta_t=
1/(\mu t)$ for the CBA as suggested by Proposition \ref{thm:CBA1} and
Proposition \ref{thm:CBA1sc}. For the MCBA, we use the step sizes
$\eta_t^k= 1/(2^{k+1}\mu)$ and $T_k=2^{k+3}$ as suggested by Proposition
\ref{thm:restart}. We choose $\mu=0.5$ in all settings.

\begin{figure*}[!h]
    \centering
    \subfigure{
        \centering
        \includegraphics[width=0.45\columnwidth]{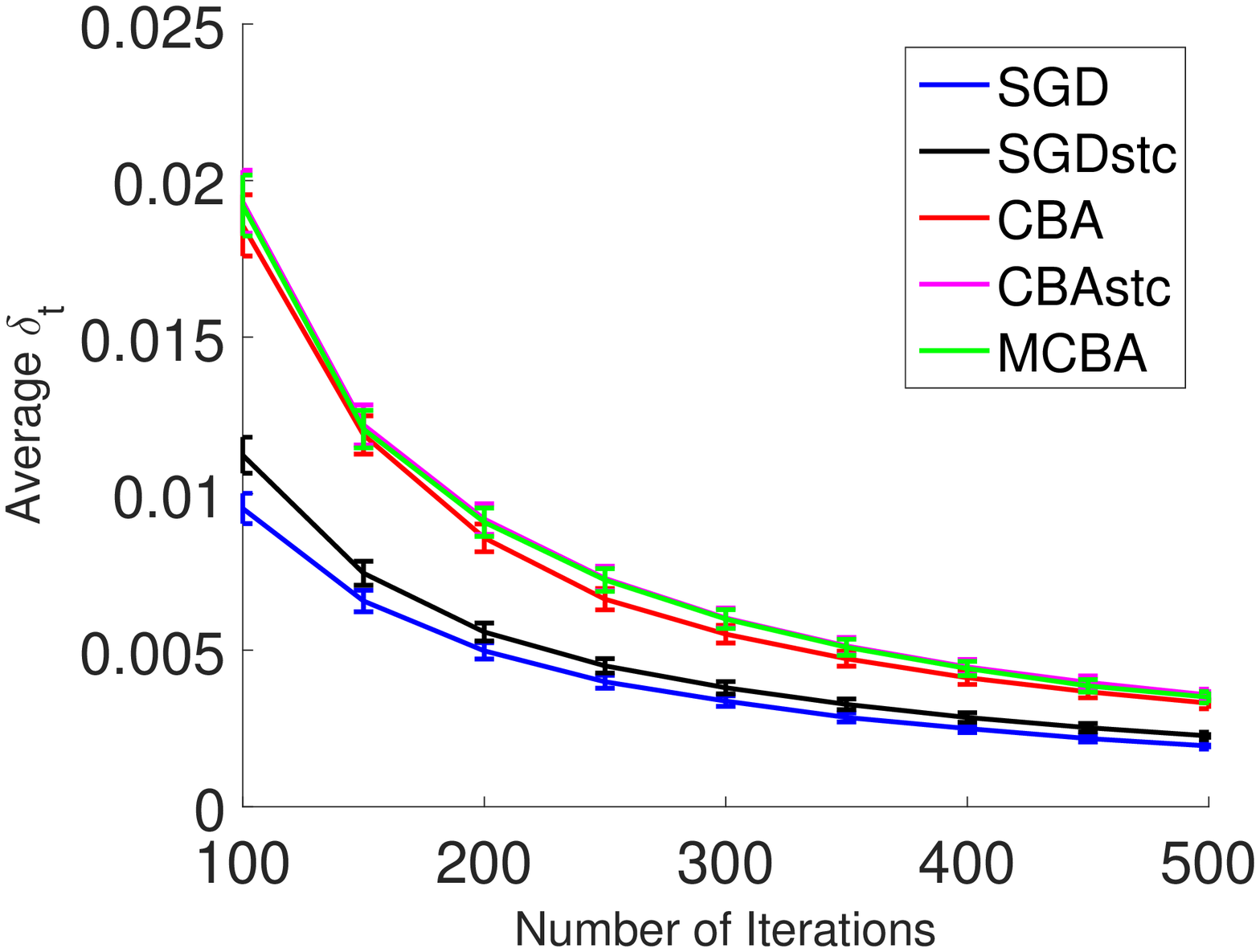}
        \label{fig:uniform_quad_uniform}
    }
    \subfigure{
        \centering
        \includegraphics[width=0.45\columnwidth]{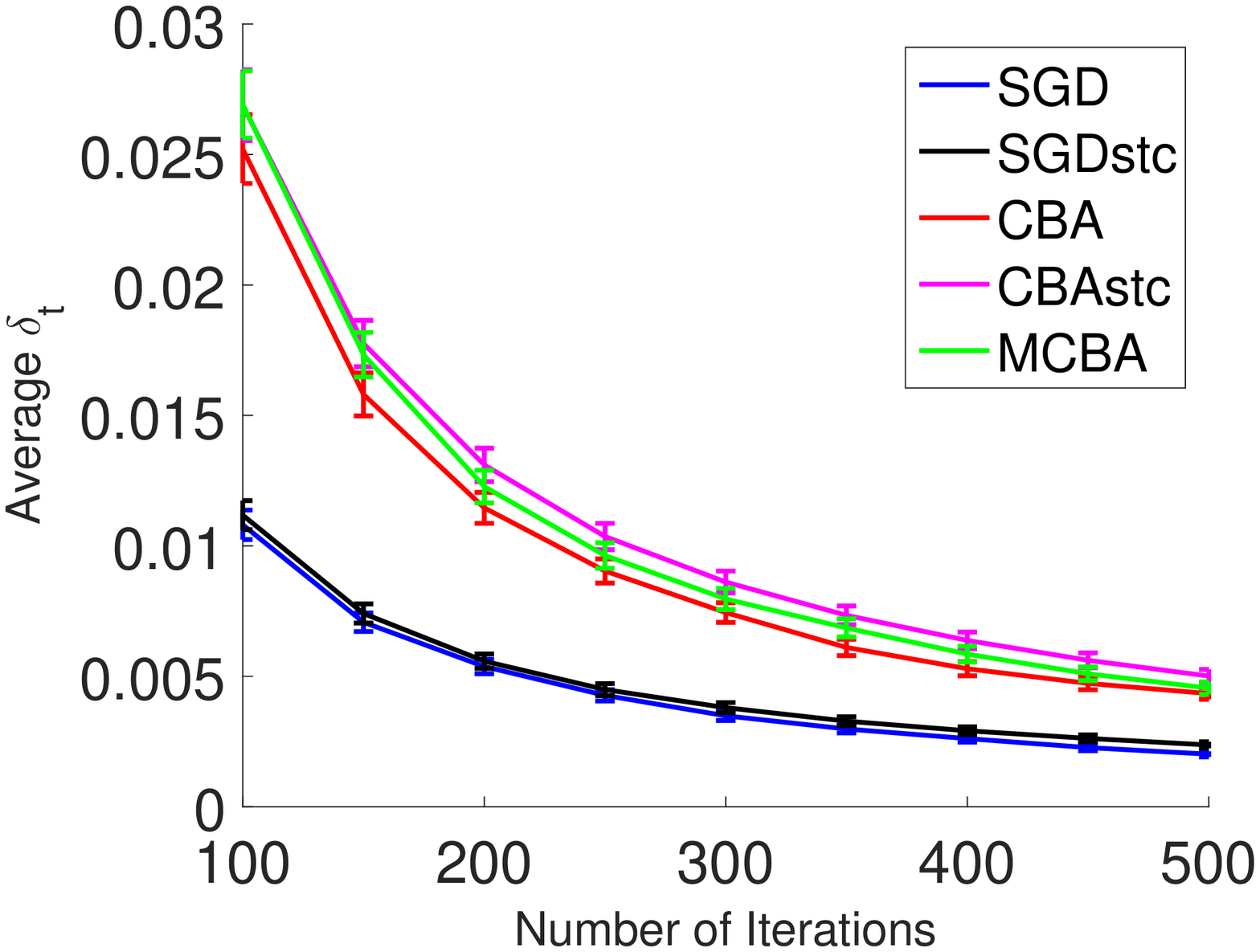}
        \label{fig:normal_quad_exp}
    }
    \subfigure{
        \centering
        \includegraphics[width=0.45\columnwidth]{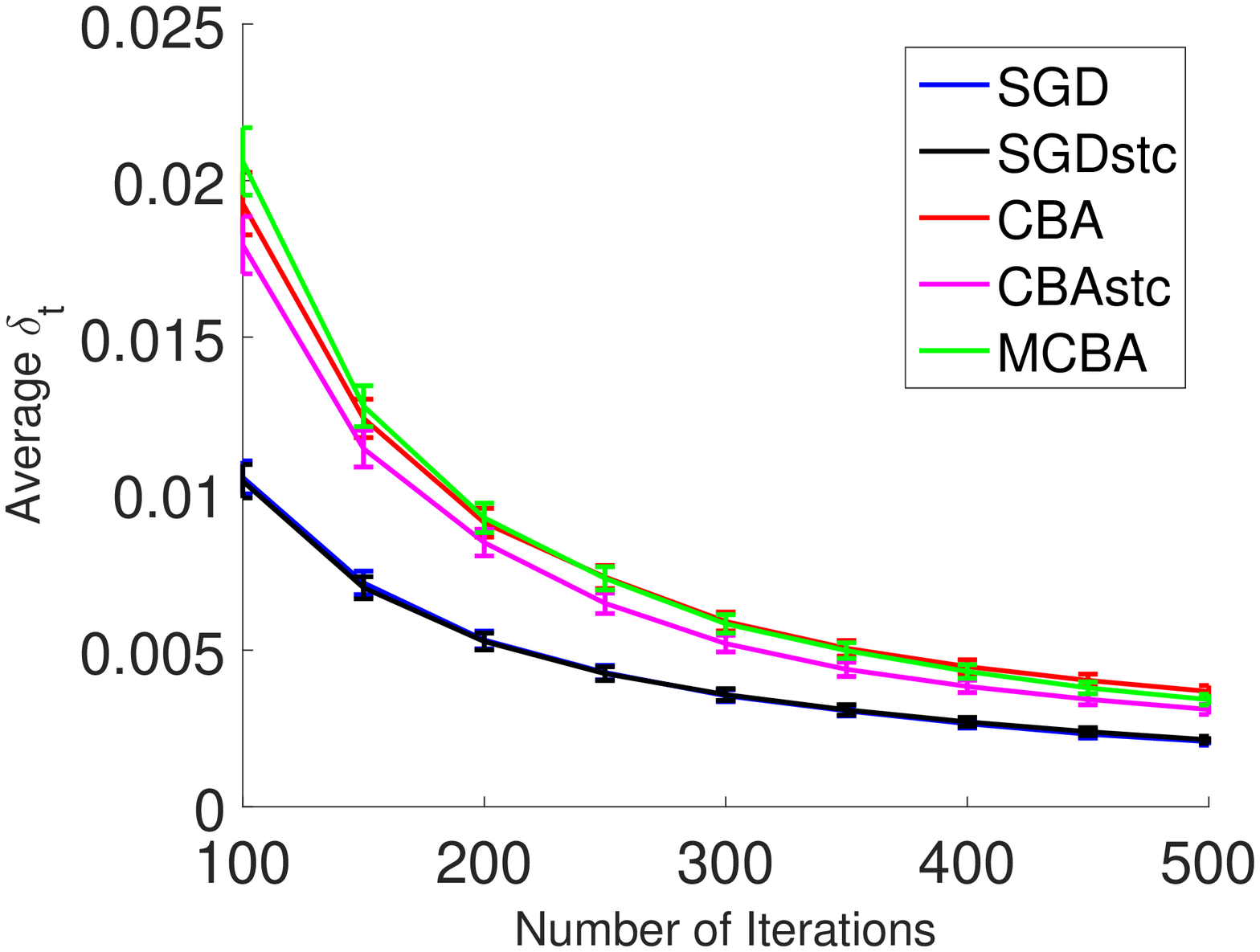}
        \label{fig:uniform_pquad_uniform}
    }
    \subfigure{
        \centering
        \includegraphics[width=0.45\columnwidth]{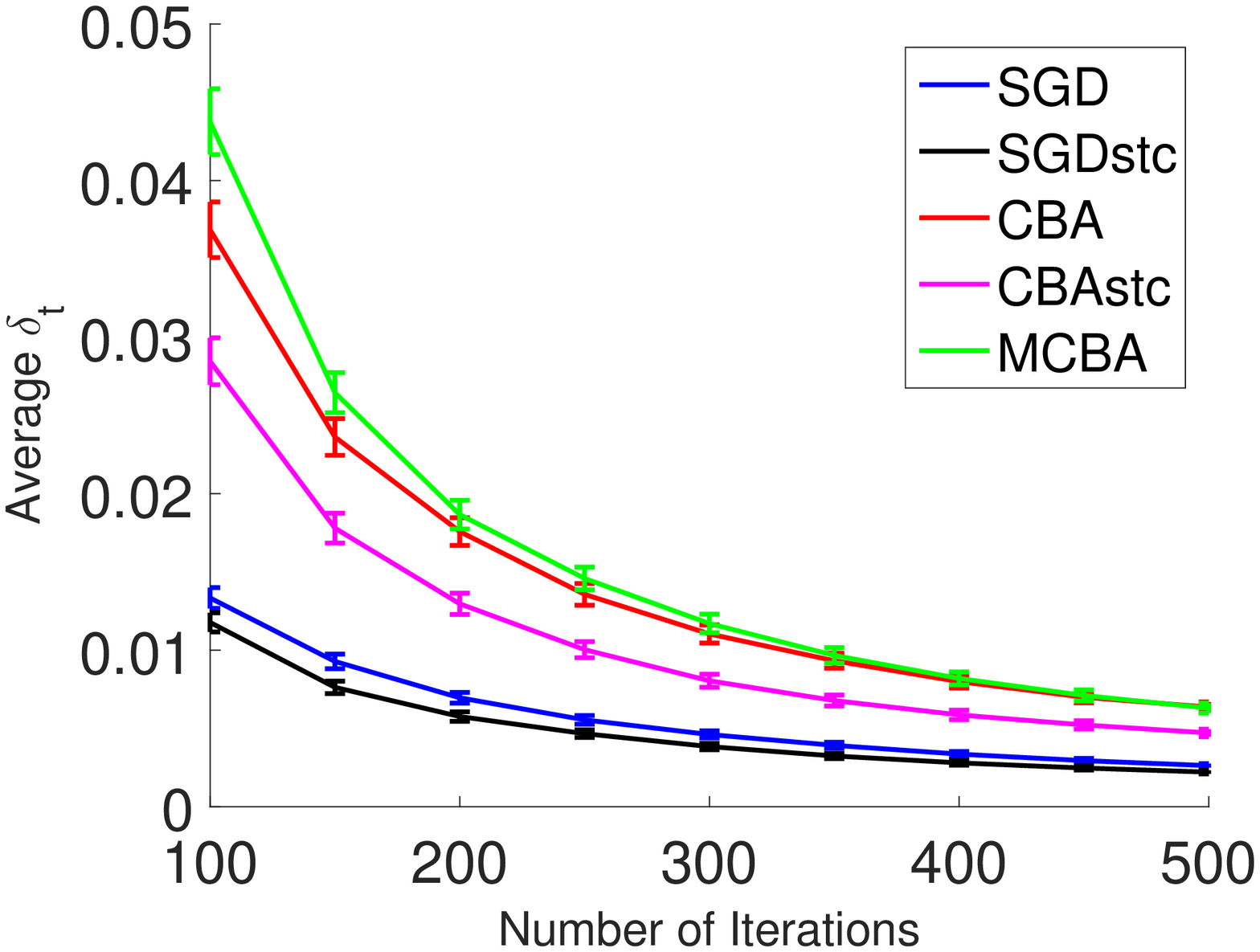}
        \label{fig:normal_pquad_exp}
    }
    \caption{The convergence of the average relative
        optimality gap for different algorithms for the four instances.
        First row: $h_1(x,\xi)$; Second row: $h_2(x,\xi)$. First
column: $\xi\sim\mathcal{U}[50,150]$;  Second column:
$\xi\sim\mathcal{N}(100,100)$.}
    \label{fig:synthetic}
    \vspace{-0.3cm}
\end{figure*}

In the following, we compare the CBA and the MCBA with the standard
stochastic gradient descent (SGD) method (see
\citealt{Nemirovski:09,Duchi:09}). In the standard SGD method, it is
assumed that $\xi_t$ can be observed directly and the stochastic
gradient update step $x_{t+1}=\text{Proj}_{[\ell,
u]}\left(x_t-\eta_th_x'(x_t, \xi_t)\right)$ is performed in each
iteration. In the experiments, we apply the same step lengths in the
SGD method as in the CBA. In all tests, we start from a random
initial point $x_1 \sim\mathcal{U}[50,150]$ and run each algorithm
for $T=500$ iterations and we report the average relative optimality
gap $\delta_t = \frac{H(\bar x_t)-H(x^*)}{H(x^*)}$ over $2000$
independent trails for each $t$. The results are reported in Figure
\ref{fig:synthetic}.

In Figure \ref{fig:synthetic}, the $x$-axis represents the number of
iterations and the $y$-axis represents the average relative
optimality gap for each algorithm. The curves \textbf{CBA} and
\textbf{SGD} represent the results of CBA and SGD with $\eta_t=
1/\sqrt{t}$ respectively while the curves \textbf{CBAstc} and
\textbf{SGDstc} represent the results of CBA and SGD with $\eta_t=
1/(\mu t)$ respectively. The curve \textbf{MCBA} represents the
results of the MCBA algorithm.  In each of the curves,
    the bar at each point represents the standard error of the
    corresponding $\delta_t$. As one can see, the standard errors are fairly
    small. Thus the test results are quite stable in these numerical
    experiments. We also present the average computation time of all algorithms in Table~\ref{table1}.

\begin{table}[!h]
    \vspace{-0.3cm}
    \begin{center}
        \begin{tabular}{|c|c|c | c | c|}
            \hline
            &\multicolumn{2}{c|}{$h_1(x,\xi)$}&\multicolumn{2}{c|}{$h_2(x,\xi)$}\\\hline
            Algorithm&$\xi\sim\mathcal{U}[50,150]$& $\xi\sim\mathcal{N}(100,100)$&$\xi\sim\mathcal{U}[50,150]$& $\xi\sim\mathcal{N}(100,100)$\\\hline
            SGD&0.008&0.007&0.010&0.040\\
            SGDstc&0.007&0.007&0.009&0.038\\
            CBA&0.011&0.018&0.013&0.049\\
            CBAstc&0.011&0.017&0.013&0.045\\
            MCBA&0.013&0.021&0.017&0.055\\
            \hline
        \end{tabular}
        \caption{The computation time (in seconds) of different algorithms for 500 iterations for the four instances.
            \label{table1}
        }
        \vspace{-0.5cm}
    \end{center}
\end{table}

From the results shown in Figure \ref{fig:synthetic}, we can see
that all of CBA, CBAstc and MCBA converge quite fast in these
problems. Even though they use much less information than the SGD
and the SGDstc methods, it takes only about twice as many iterations
to get the same accuracy. As shown in
Table~\ref{table1}, the computation time for 500 iterations is less
than $0.1$ seconds in each algorithm and the time is not much
different across different algorithms. (This short runtime is
because of the low dimensionality of the problems.) Moreover, in
our tests, CBA, CBAstc and MCBA have quite similar performance
despite their different theoretical
guarantees.\footnote{ Note that in Figure
\ref{fig:synthetic}, the CBA and the CBAstc sometimes
        perform even better than the MCBA despite the worse theoretical
        guarantee. In fact, this is common in convex optimization literature
        for such types of (restarting) methods. For example, \citet{Chen:12}
        proposed a stochastic gradient method called MORDA, which improves
        ORDA in the same paper in theoretical convergence rate using a similar restarting
        technique to our MCBA method. However, in the fourth column of Table
        2 in \citet{Chen:12}, the objective value in MORDA is higher
        than that in
        ORDA. Similarly, \citet{Lan:10b} developed a stochastic gradient
        method called Multistage AC-SA, which improves AC-SA in theoretical
        convergence rate but not necessarily in numerical performance (see
        Table 4.3 in \citealt{Lan:10b}). }



\subsection{Choices of $f_{-}$ and $f_{+}$}
\label{subsec:numerical_choice_f}

\begin{figure*}[!h]
    \centering
    \subfigure{
        \centering
        \includegraphics[width=0.45\columnwidth]{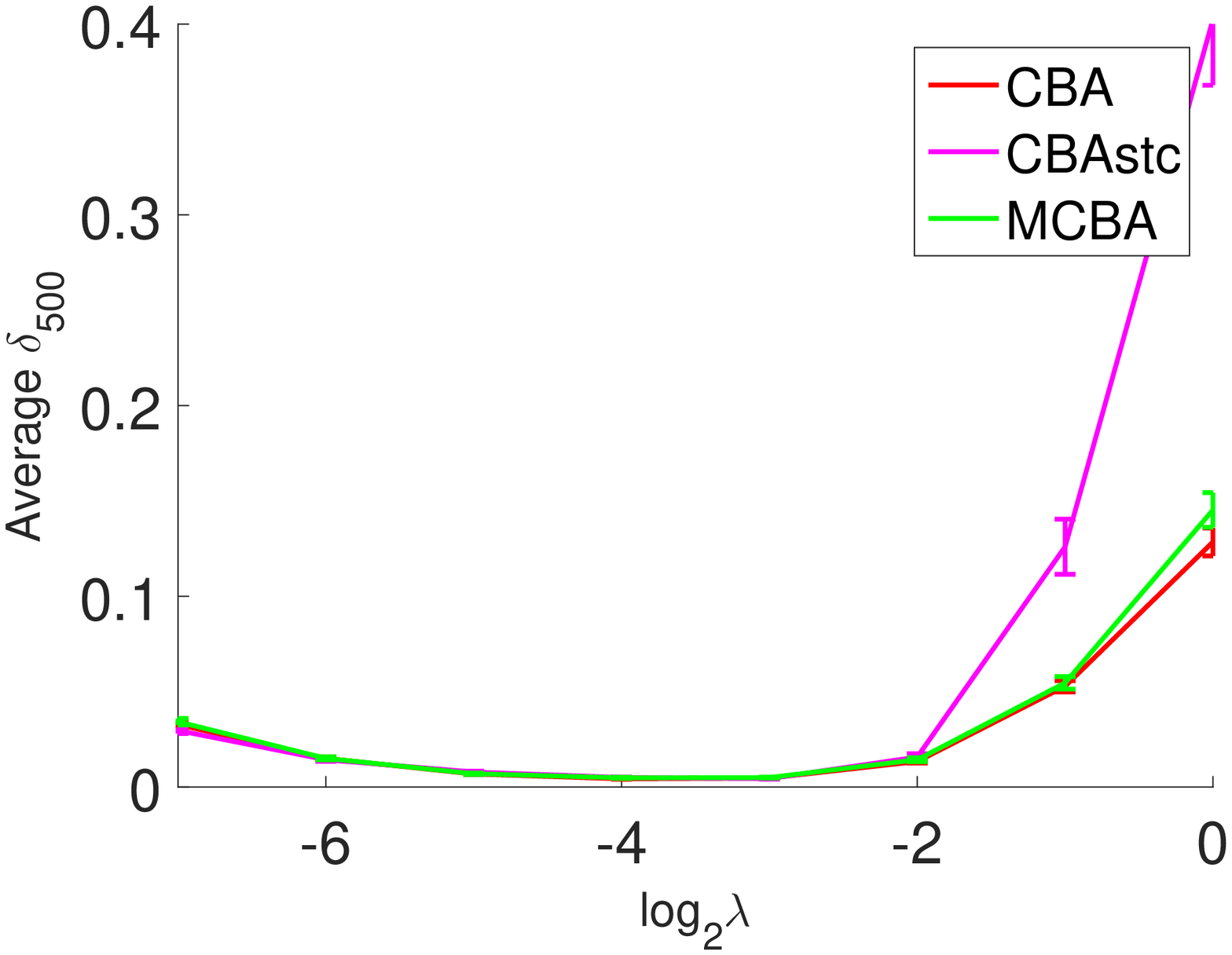}
        \label{fig:lambda_range1}
    }
    \subfigure{
        \centering
        \includegraphics[width=0.45\columnwidth]{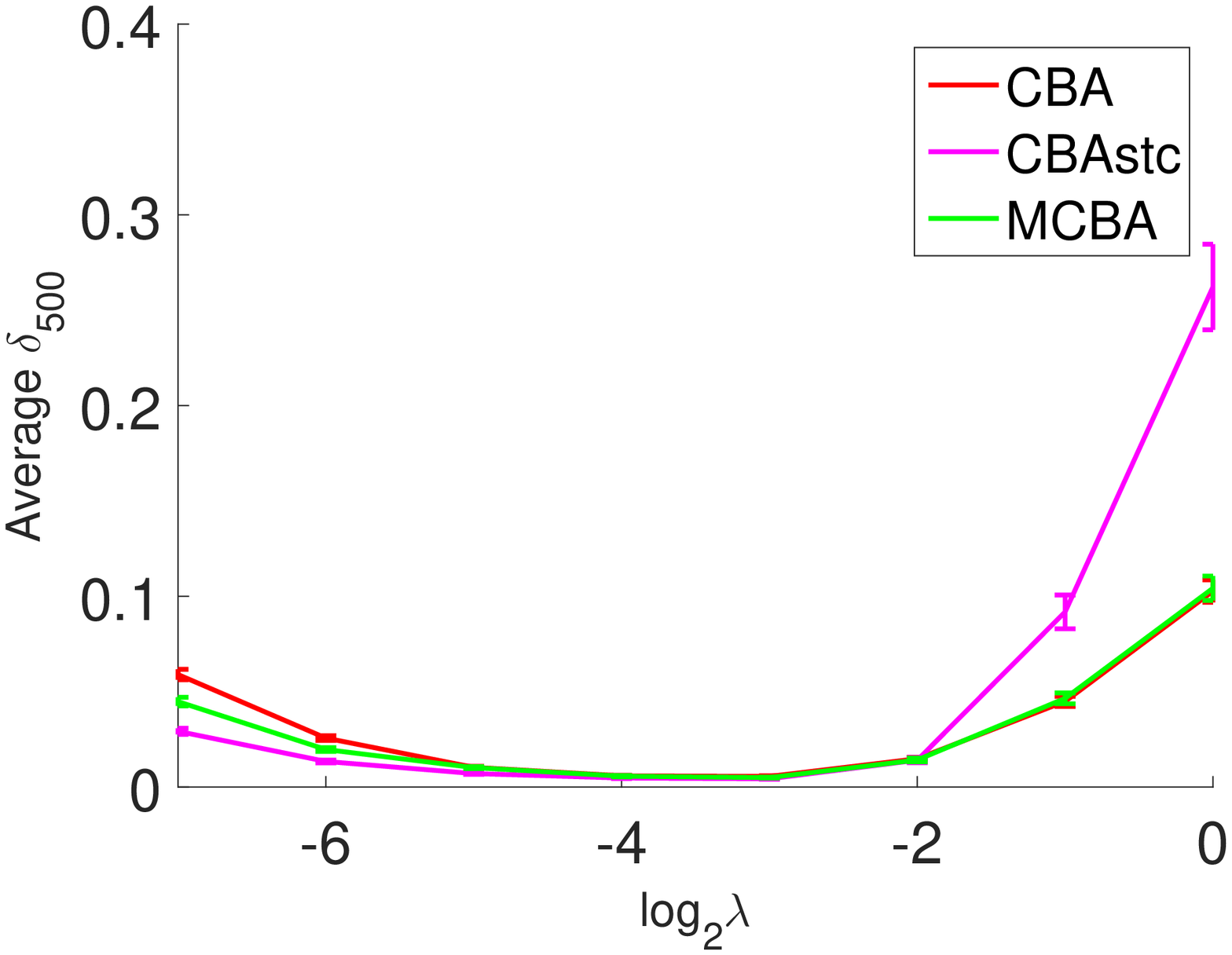}
        \label{fig:lambda_range2}
    }
    \caption{The impact of $\lambda$ in $f_+$ and $f_-$ to different algorithms, measured by
        the average relative optimality gap after 500 iterations.
        Left: $h_1(x,\xi)$; Right: $h_2(x,\xi)$. In both
        cases,  $\xi\sim\mathcal{N}(100,100)$.}
    \label{fig:synthetic_variance}
\end{figure*}

In this section, we perform some additional tests to study the
impact of different choices of $f_{-}$ and $f_{+}$ on the
performance of the CBA. We still consider the two objective
functions considered in the last section, but we focus on the case
in which $\xi\sim\mathcal{N}(100,100)$. First we keep $f_{-}$ and
$f_{+}$ to be exponential distributions (as in Example
\ref{example:fexponential}) and see how the performance is affected
by different values of $\lambda_{-}=\lambda_{+}=\lambda$. The
results are shown in Figure \ref{fig:synthetic_variance}.

In Figure \ref{fig:synthetic_variance}, the $x$-axis represents the
value of $\log_2\lambda$ while the $y$-axis represents the relative
optimality gap $\delta_T$ after $500$ iterations evaluated as the
average value of $2000$ independent trials  (again the
    bars show the standard error in these trials). The influences of
different values of $\lambda$ in $f_{-}$ and $f_{+}$ are presented
for the CBA, the CBAstc and the MCBA algorithms. We can see that the
value of $\lambda$ does influence the convergence speed of the
algorithms. Particularly, in both figures in Figure
\ref{fig:synthetic_variance}, the optimality gap after $T=500$
iterations decreases first as $\lambda$ increases but starts to
increase when $\lambda$ is large. Moreover, the influence is
relatively small when $\lambda$ is small but is large when $\lambda$
is large. And the influence is more pronounced for the CBAstc
algorithm. In our setting, the best performance is obtained around
$\lambda=2^{-4}=0.0625$ for all algorithms.

\begin{figure*}[!h]
    \centering
    \subfigure{
        \centering
        \includegraphics[width=0.45\columnwidth]{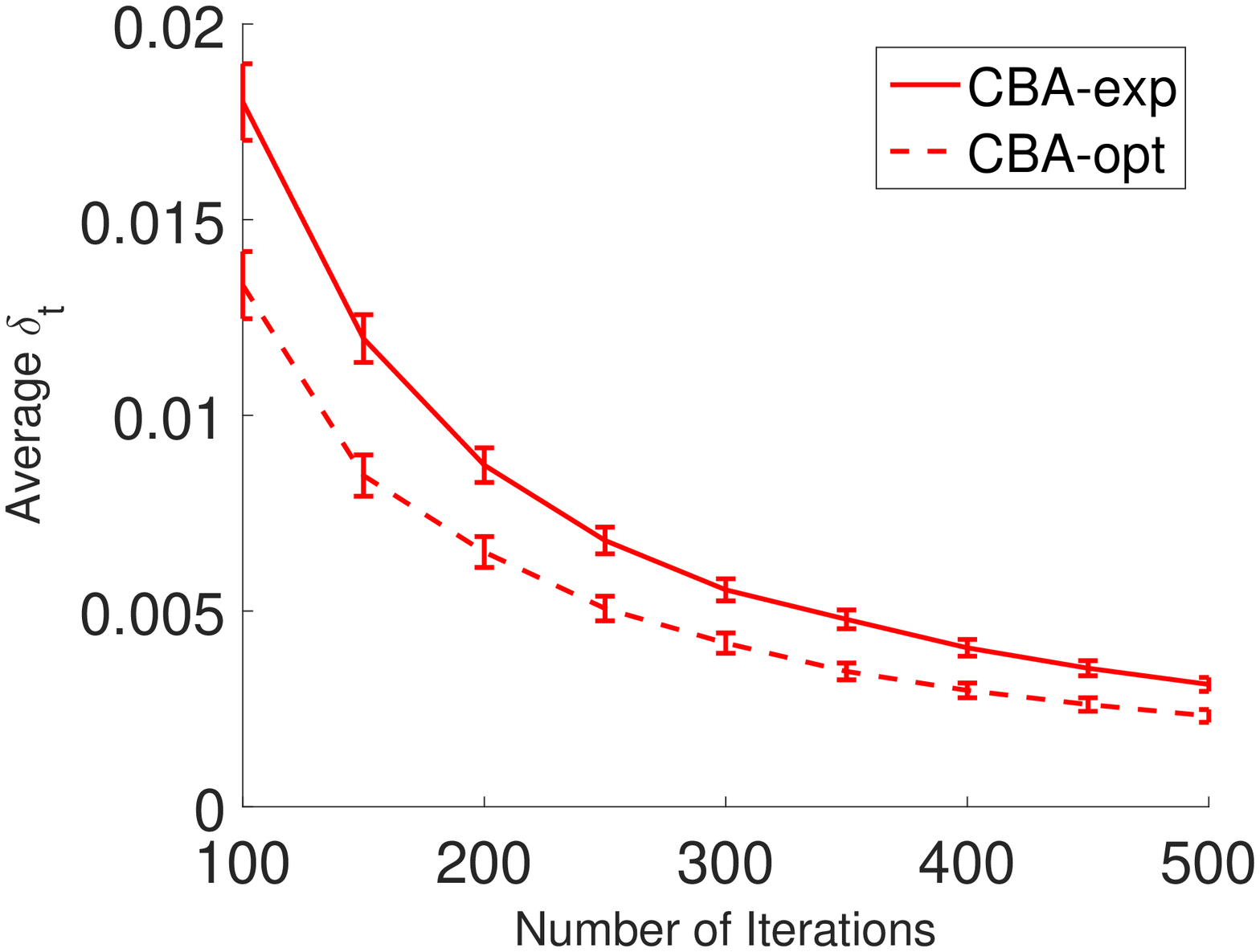}
        \label{fig:optf1}
    }
    \subfigure{
        \centering
        \includegraphics[width=0.45\columnwidth]{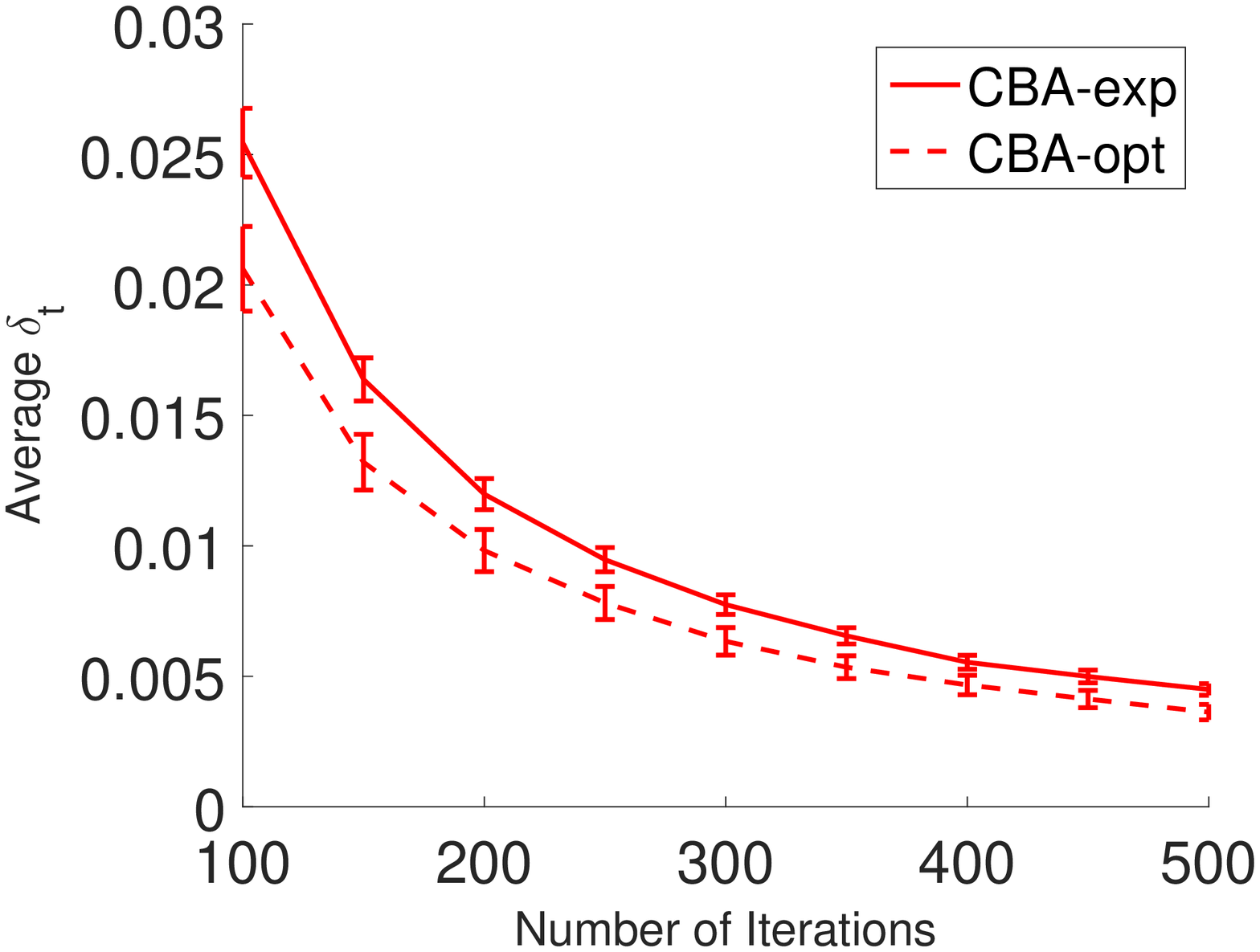}
        \label{fig:optf1}
    }
    \caption{The convergence of the average relative
        optimality gap in CBA when choosing $f_+$ and $f_-$ to be the exponential distribution and the optimal distribution in \eqref{eq:optfneg} and \eqref{eq:optfpos}.
        Left: $\xi\sim\mathcal{U}[50,150]$; Right: $\xi\sim\mathcal{N}(100,100)$. In both cases, $h(x,\xi)=(x-\xi)^2$.}
    \label{fig:optf}
\end{figure*}

In Section \ref{sec:choice of f}, we provided the optimal choice of
$f_{-}$ and $f_{+}$ in \eqref{eq:optfneg} and \eqref{eq:optfpos}. In
Figure~\ref{fig:optf}, we present the difference in the performances
of the CBA when $f_{-}$ and $f_{+}$ are chosen optimally versus when
they are chosen as exponential distributions with
$\lambda_{-}=\lambda_{+}=0.0625$. In this experiment, we choose
$h(x,\xi)=h_1(x,\xi)$ and run the CBA for 500 iterations and compute
the average relative optimality gap $\delta_{500}$ over $2000$
independent trials. The results are plotted in Figure
\ref{fig:optf}.

In Figure \ref{fig:optf}, we can see that the optimal choice of
$f_{-}$ and $f_{+}$ does improve the performance of the CBA, which
confirms our analysis in Section \ref{sec:choice of f}. However, the
improvement is not essential yet generating samples from the optimal
distribution is much more time consuming. For example, when $\xi\sim\mathcal{N}(100,100)$,
the computational time is less than $0.05$ seconds  when using exponential distribution
but about $1$ second when using the optimal distribution. Therefore, as we discussed in the end of Section
\ref{sec:choice of f}, one can just choose a simple distribution in
practice.

\subsection{Mini-Batch Method with Additional Comparisons}
\label{subsec:numerical_minibatch} 
\begin{figure*}[!h]
        \centering
        \subfigure{
            \centering
            \includegraphics[width=0.45\columnwidth]{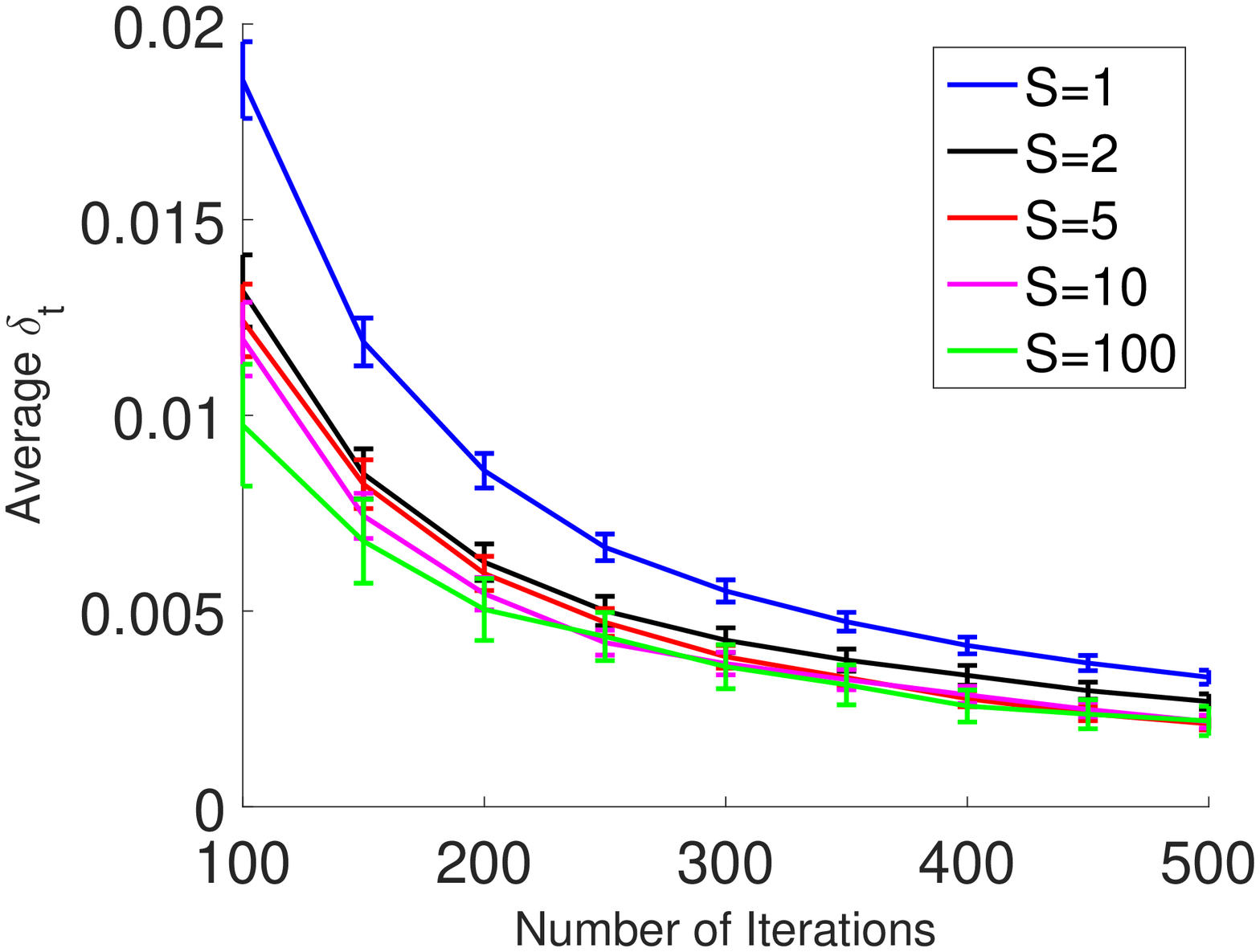}
            \label{fig:uniform_quad_uniform}
        }
        \subfigure{
            \centering
            \includegraphics[width=0.45\columnwidth]{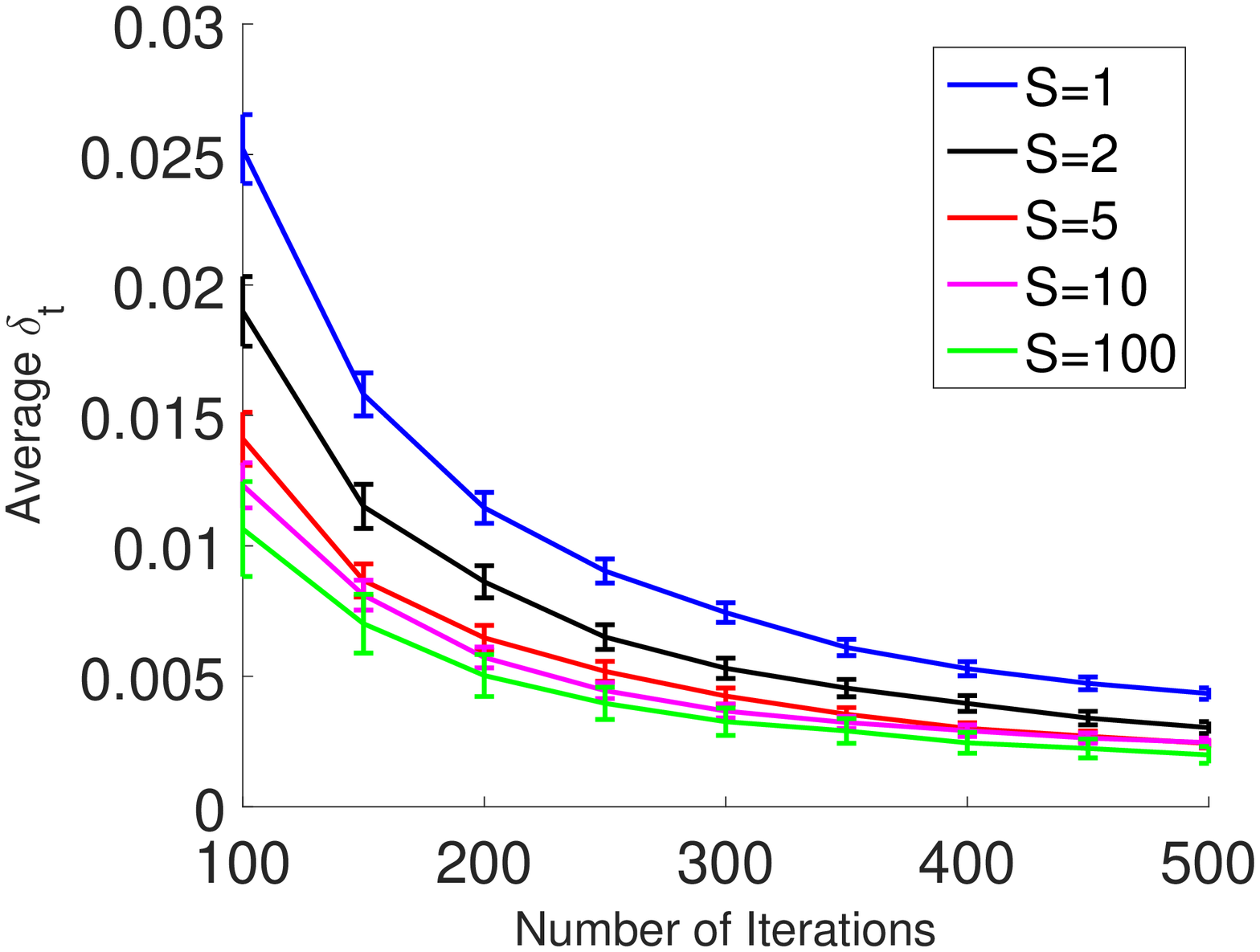}
            \label{fig:normal_quad_exp}
        }
        \subfigure{
            \centering
            \includegraphics[width=0.45\columnwidth]{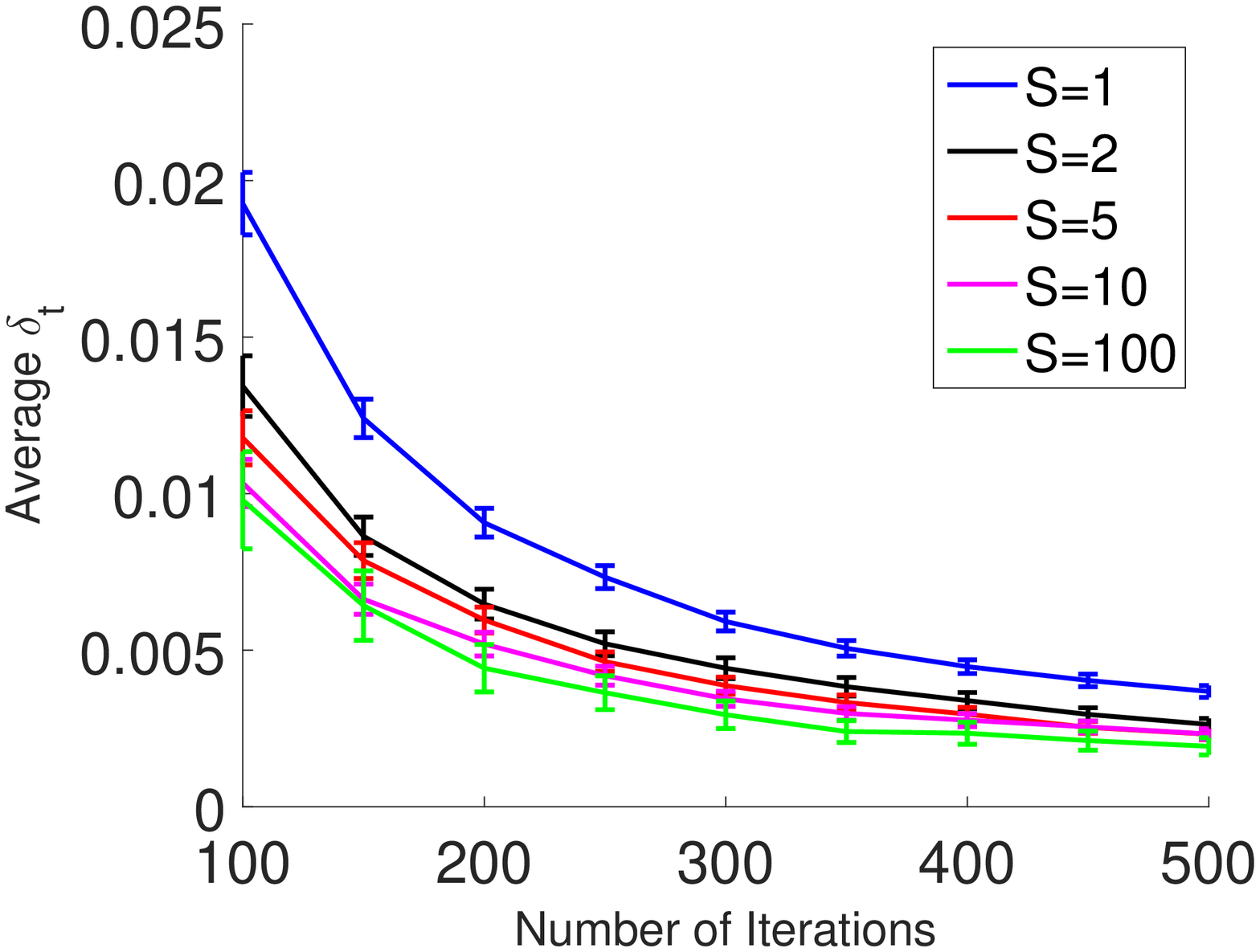}
            \label{fig:uniform_pquad_uniform}
        }
        \subfigure{
            \centering
            \includegraphics[width=0.45\columnwidth]{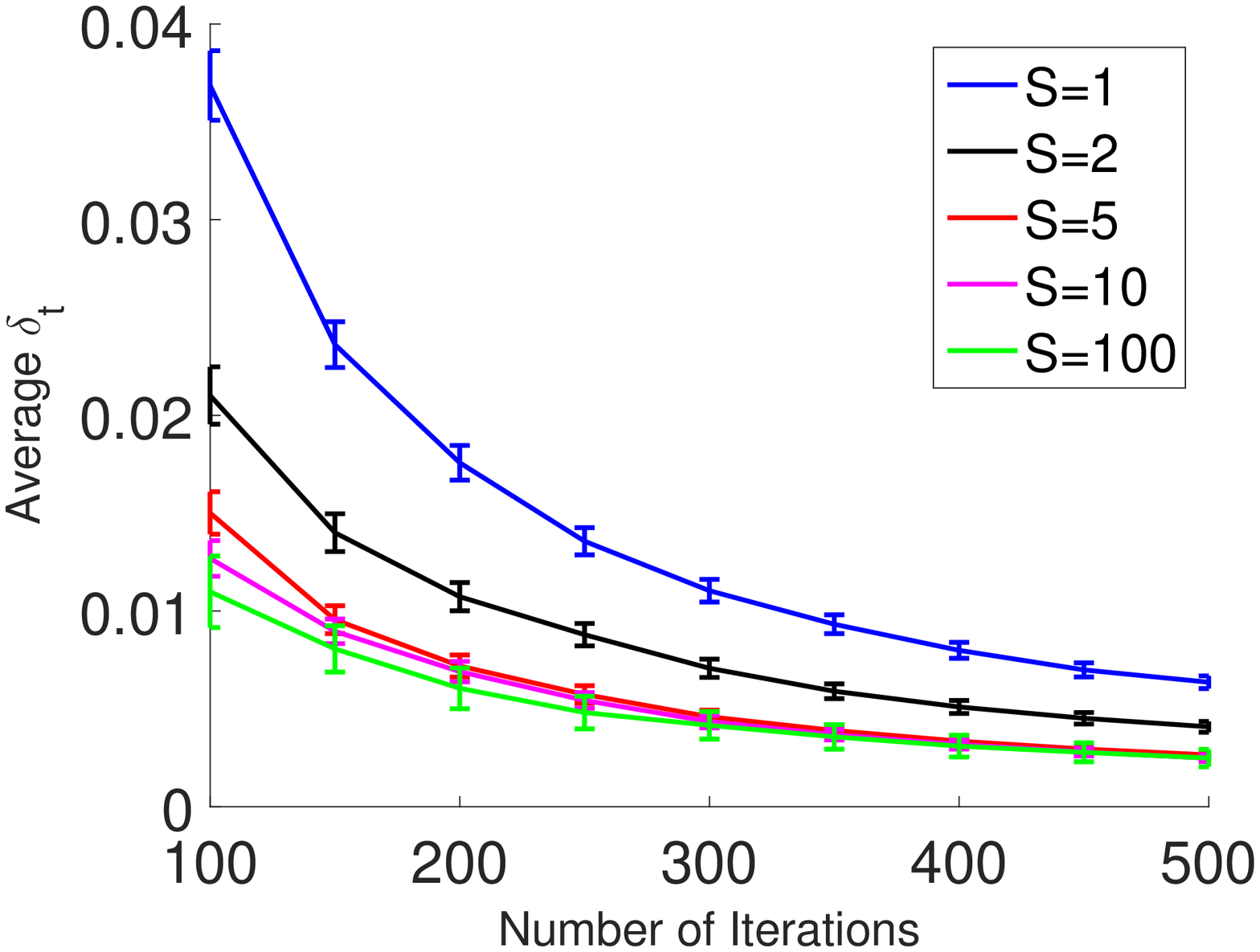}
            \label{fig:normal_pquad_exp}
        }
        \caption{The convergence of the average relative
            optimality gap in CBA when using different numbers of comparisons ($S$).
            First row: $h_1(x,\xi)$; Second row: $h_2(x,\xi)$. First
            column: $\xi\sim\mathcal{U}[50,150]$;  Second column:
            $\xi\sim\mathcal{N}(100,100)$.}
        \label{fig:synthetic_MB}
    \end{figure*}

In this section, we numerically
test how the performance of CBA depends on the sample size $S$ in
the mini-batch technique described in Section~\ref{sec:MB}. The
instances and the choice of parameters are all identical to
Section~\ref{subsec:numerical:convergence}. We present the
convergence of CBA with $S=1,2,5,10,100$ in
Figure~\ref{fig:synthetic_MB} and the associated runtimes in
Table~\ref{table2}. According to the figures, one additional
comparison (i.e. $S=2$) with $z$ can improve the convergence of
Algorithm~\ref{alg:cba1} in all four instances. Although increasing
$S$ can still improve the performance further, the effect diminishes
quickly. This is because the noise in the stochastic gradient is
generated from the sample noise of both $\xi$ and $z$. The
mini-batch technique for sampling $z$ does not help reduce the noise
due to $\xi$, which eventually dominates the noise due to $z$ when
$S$ is large enough.\footnote{In Appendix
\ref{appendix:additional_numerical}, we also present numerical
results of convergence with mini-batch method with respect to
runtime of the algorithm. The results show that choosing a small
batch size can usually lead to fastest convergence (in terms of
time). This is because of the tradeoff between the reduced number of
iterations required in the mini-batch method and the increased
computational efforts required in each iteration.} Although a
similar mini-batch technique can be applied to $\xi$, it might not
be practical when applying our algorithm in practice. For example,
when $\xi$ represents the ideal product of a customer, creating a
mini-batch for $\xi$ means asking different customers' preference
without updating the solution, which is not consistent with the
setting of online optimization where the solution is updated after
the visit of each customer.

    \begin{table}[!h]
        \begin{center}
            \begin{tabular}{|c|c|c | c | c|}
                \hline
                &\multicolumn{2}{c|}{$h_1(x,\xi)$}&\multicolumn{2}{c|}{$h_2(x,\xi)$}\\\hline
                $S$&$\xi\sim\mathcal{U}[50,150]$& $\xi\sim\mathcal{N}(100,100)$&$\xi\sim\mathcal{U}[50,150]$& $\xi\sim\mathcal{N}(100,100)$\\\hline
                1&0.011&0.018&0.013&0.049\\
                2&0.015&0.022&0.020&0.053\\
                5&0.021&0.025&0.022&0.057\\
                10&0.026&0.033&0.029&0.064\\
                100&0.117&0.158&0.125&0.175\\
                \hline
            \end{tabular}
            \caption{The computation time (in seconds) of CBA for 500 iterations when using different numbers of comparisons ($S$) per iteration.
                \label{table2}
            }
        \end{center}
    \end{table}

\subsection{Multi-Dimensional Problems}
\label{subsec:numerical:convergence_muti} In this section, we
conduct numerical experiments to test the performance of the CBA-QP
and the MCBA-QP on the multi-dimensional stochastic quadratic
program~\eqref{eq:QPobj}. To generate a testing instance
of~\eqref{eq:QPobj}, we first generate a $d\times d$ matrix $Q'$
where each entry is sampled from an i.i.d. standard normal
distribution $\mathcal{N}(0,1)$ and then set $Q=(Q')^\top Q'/d+I_d$
in~\eqref{eq:QPobj}, where $I_d$ is a $d\times d$ identity matrix.
We choose the random variable $\xi$ in \eqref{eq:QPobj} with a
multivariate normal distribution
$\mathcal{N}(100\mathbf{1}_d,50^2I_d)$, where $\mathbf{1}_d$ is the
all-one vector in $\mathbb{R}^d$. It is easy to show that, with this
choice, Assumption~\ref{assumption:qp} holds and we have
$\mathbb{E}\exp(\|\xi\|_2)\leq\bar{\sigma}$ for some $\bar{\sigma} >
0$ so that we can choose $f$ in Algorithms~\ref{thm:CBA1-QP}
and~\ref{thm:restart-QP} to be an exponential distribution in order
to guarantee (C4) according to Example~\ref{example:fexponential2}.
More specifically, we choose $f$ to be the exponential distribution
\eqref{eq:fzexpdist} with $\lambda=2^{-4}=0.0625$ (same as
$\lambda_+$ and $\lambda_-$ in
Section~\ref{subsec:numerical:convergence}). Finally, we choose the
feasible set of~\eqref{eq:QPobj} to be the box
$\mathcal{X}=\prod_{1}^d[50,150]$.

We compare the performance of the SGD, CBA-QP and MCBA-QP methods.  In the SGD method, each iteration computes $x_{t+1}=\text{Proj}_{\mathcal{X}}\left(x_t-\eta_tQ(x_t-\xi_t)\right)$ where $\xi_t$ is sampled from the distribution of $\xi$ and $Q(x_t-\xi_t)$ is the gradient of $h(x,\xi)$ with respect to $x$. As suggested by Proposition \ref{thm:CBA1-QP}, we choose $\eta_t=
1/(\mu t+L)$ in CBA-QP where $\mu$ and $L$ are the smallest and largest eigenvalues of $Q$, respectively. For MCBA, we choose $\eta_t^k=\frac{1}{2^{k+1}\mu+L}$ and $T_k=2^{k+3}+4$ as suggested by Proposition
\ref{thm:restart-QP}. In the experiments, we apply the same step lengths in the
SGD method as in the CBA-QP method.

\begin{figure*}[!h]
    \centering
    \subfigure{
        \centering
        \includegraphics[width=0.45\columnwidth]{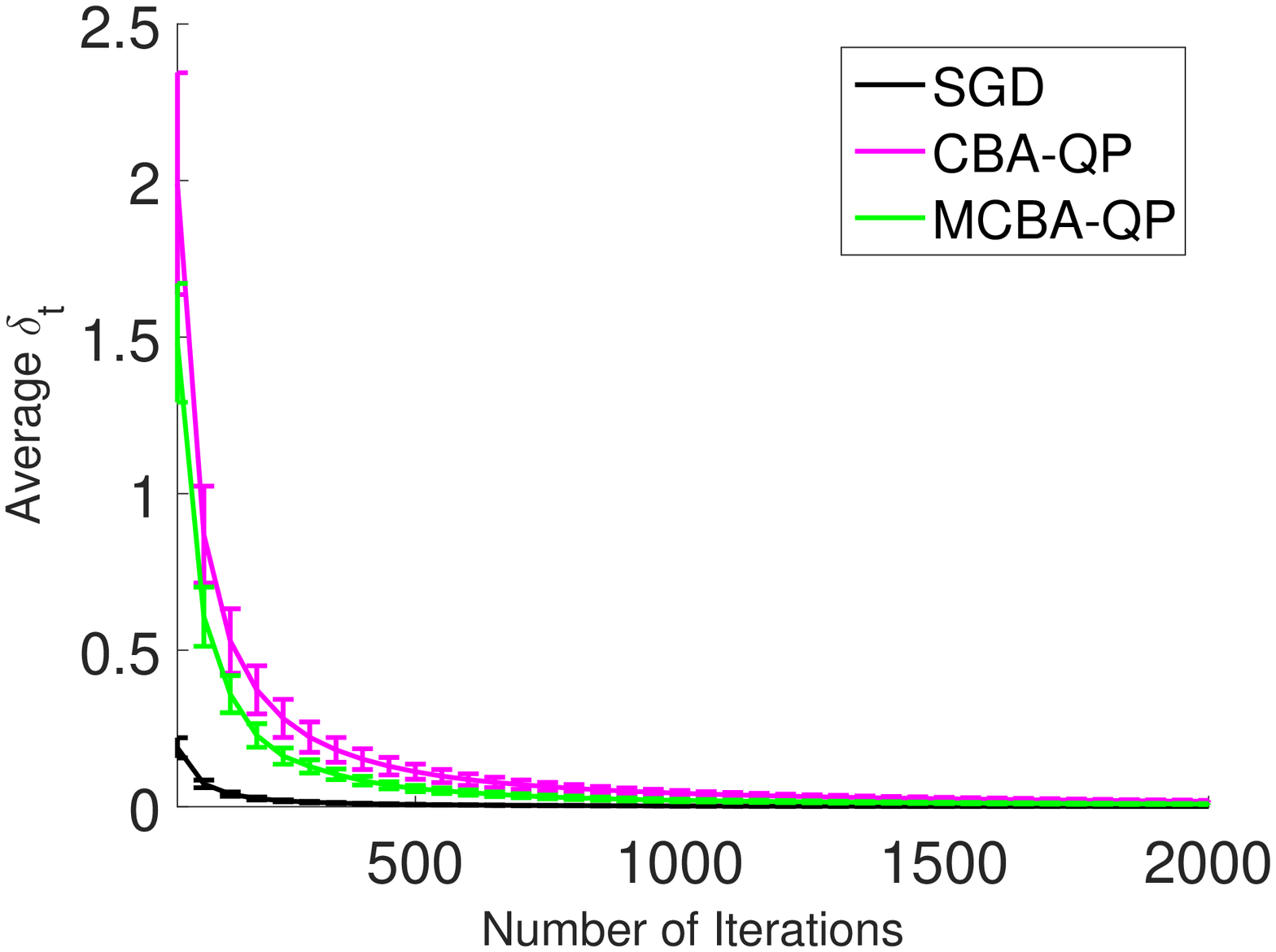}
        \label{fig:normal_QP_exp_d5}
    }
    \subfigure{
        \centering
        \includegraphics[width=0.45\columnwidth]{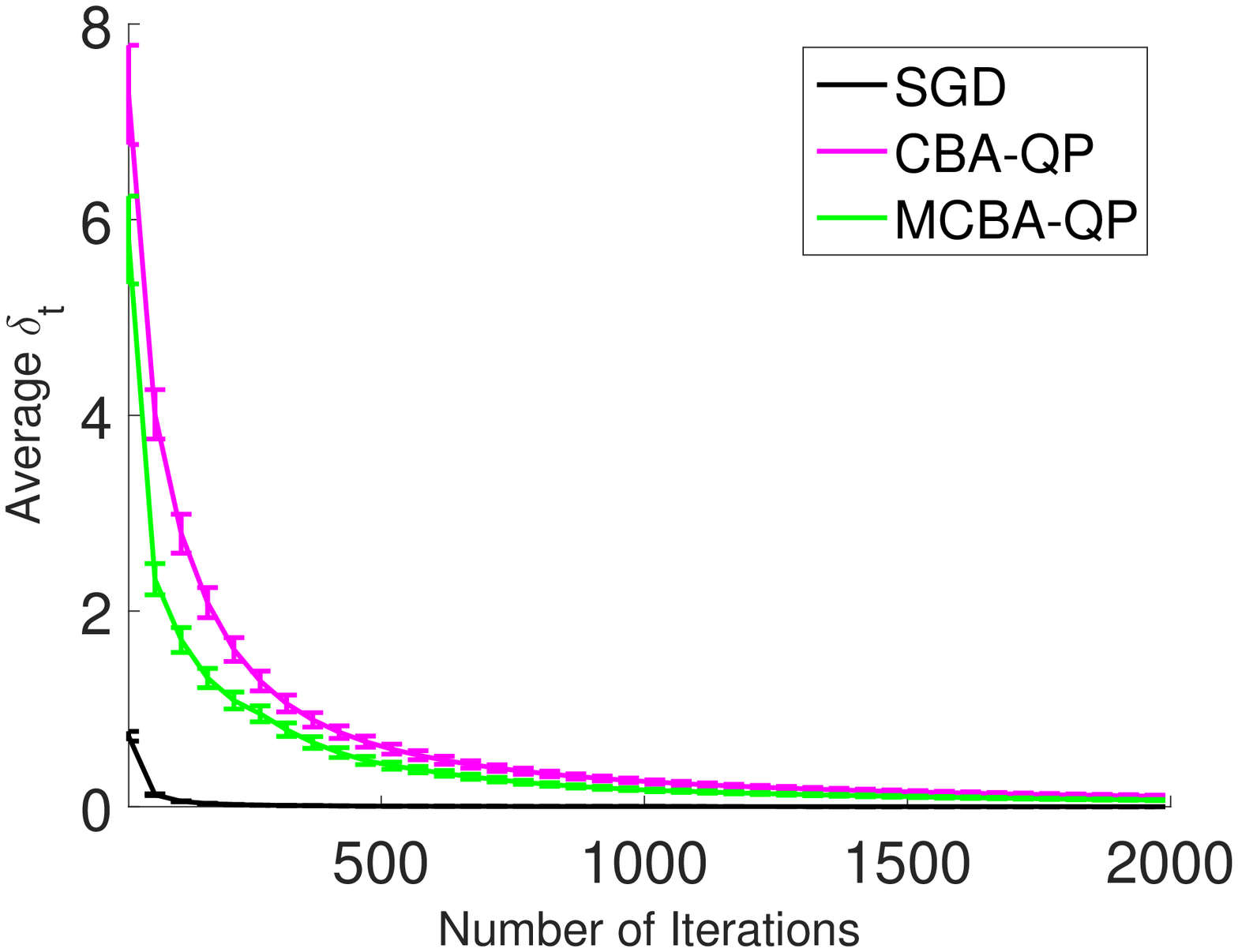}
        \label{fig:normal_QP_exp_d20}
    }
    \caption{The convergence of the average relative
        optimality gap for different algorithms for the multi-dimensional quadratic program~\eqref{eq:QPobj}.
        Left: $d=5$; Right: $d=20$.}
    \label{fig:synthetic-QP}
    \vspace{-0.3cm}
\end{figure*}

We randomly generate  $Q$ in the aforementioned method, start each algorithm at the same initial point $x_1$ that is uniformly randomly sampled in the box $\mathcal{X}$, and run each algorithm
for $T=2000$ iterations. We report the average relative optimality
gap $\delta_t = \frac{H(\bar x_t)-H(x^*)}{H(x^*)}$ and its standard error over $2000$
independent trails for each $t$. We run these experiments with the dimension $d=5$ and $d=20$.
The results are reported in Figure
\ref{fig:synthetic-QP}.

From the results shown in Figure \ref{fig:synthetic-QP}, we see that
even though CBA-QP and MCBA-QP use much less information than the
SGD method, they can eventually find a solution with an accuracy
comparable to SGD. As shown in Table~\ref{table1}, the computation
time for 2000 iterations is less than one second in each algorithm
and the time is not much different across different algorithms.
(This short runtime is because of the low dimensionality of the
problems.) Moreover, in our tests, MCBA-QP has a faster convergence
rate than CBA-QP which is consistent with Proposition~\ref{thm:CBA1-QP}
and~\ref{thm:restart-QP}.  However,  CBA-QP and  MCBA-QP converge
more slowly than  SGD. This is because  SGD utilizes the
information of $\xi$ with the full dimension while  CBA-QP and
 MCBA-QP can only exploit information along a random direction
$u$. This difference is more significant as the dimension increases,
which is confirmed by Figure 5.

\begin{table}[!h]
    \vspace{-0.3cm}
    \begin{center}
        \begin{tabular}{|c|c| c|}
            \hline
            Algorithm&\multicolumn{1}{c|}{$d=5$}&\multicolumn{1}{c|}{$d=20$}\\\hline
            SGD&0.236&0.344\\
            CBA-QP&0.391&0.490\\
            MCBA-QP&0.440&0.585\\
            \hline
        \end{tabular}
        \caption{The computation time (in seconds) of different algorithms for 2000 iterations for the four instances.
            \label{table1}
        }
        \vspace{-0.5cm}
    \end{center}
\end{table}

\section{Concluding Remarks}
\label{sec:discussions}

In this paper, we considered a stochastic optimization problem when
neither the underlying uncertain parameters nor the objective value
at the sampled point can be observed. Instead, the decision maker
can only access to comparative information between the sample point
and two chosen decision points in each iteration. We proposed an algorithm
that gives unbiased gradient estimates for this problem, which
achieves the same asymptotic convergence rate as standard stochastic
gradient methods. Numerical experiments demonstrate that our
proposed algorithm is efficient.

There is one remark we would like to make. In this paper,
we assumed that $\xi$ follows a continuous distribution. However, we
only need that in each iteration, the probability that $\xi_t = x_t$
is $0$. This can be guaranteed by only requiring ${\mathbb P}(\xi =
\ell) = {\mathbb P}(\xi = u) = 0$, and then in the CBA, adding a
small and decaying random perturbation to $g(x_t,\xi_t,z_t)$ (for
example, a uniform distribution on $[-1/2^t, 1/2^t]$ in iteration
$t$). By doing this, one can still use the same analysis, and the
same performance guarantee holds for the modified algorithm.

 There are several future directions of research.
First, for multi-dimensional problems, we only considered convex
quadratic problems in this paper (as mentioned in Section
\ref{sec:QP}, we can also generalize our method to separable
objective function cases). It would be of interest to see whether
the ideas and techniques can be generalized to more general
multi-dimensional settings. Second, in this paper, we assumed that
the distribution of $\xi$ is stationary over time, and the
comparative information is always reported accurately.
However, in practice, the distribution of $\xi$ may change over time
or the comparative information may be reported in a noisy
fashion. It would be interesting to see whether we could extend our
discussions to consider such situations. Finally, our paper only
considers a continuous decision setting, i.e., we assumed that all
the decision variables as well as the test variables can take any
continuous values. In many practical situations, the decision
variables may only be chosen from a finite set. It is worth further
research to see whether a similar idea can be applied in such
settings.


\section*{Acknowledgments}
The authors thank the editor-in-chief, the associated editor and three anonymous referees for many useful suggestions and feedback.
Xi Chen was supported by the NSF IIS-18454444, the Alibaba Innovation Research Award, and the Bloomberg Data Science Research Grant.

\bibliographystyle{ormsv080}
\bibliography{reference}

\newpage
\APPENDICES

\section{Proofs of Proposition~\ref{thm:CBA1}, \ref{thm:CBA1sc} and
\ref{thm:restart}} \label{sec:convproof}

{\noindent\bf Proof of Proposition~\ref{thm:CBA1}.}
\label{appendix:proof_of_theorem1-3}
We prove the proposition by considering the case when Assumption A5(a)
or A5(b) holds respectively. For the ease of notation, we shall use
$\mathbb{E}_t$ to denote the conditional expectation taken over
$\xi_t$ and $z_t$ conditioning on
$\xi_1,z_1,\xi_2,z_2,\dots,\xi_{t-1},z_{t-1}$.

\noindent 1. When Assumption A5(a) holds (based on
\citealt{Nemirovski:09,Duchi:09}):

According to \eqref{eq:proj}, we have
\begin{eqnarray}
\nonumber ( x_{t+1}-x^*)^2&=&(\text{Proj}_{[\ell,
    u]}(x_t-\eta_tg(x_t, \xi_t, z_t))-\text{Proj}_{[\ell,
    u]}(x^*))^2\\\nonumber &\leq&(x_t-\eta_tg(x_t, \xi_t,
z_t)-x^*)^2\\\label{sa_eq3} &=&(x_t-x^*)^2-2\eta_tg(x_t, \xi_t,
z_t)(x_t-x^*)+\eta_t^2(g(x_t, \xi_t, z_t))^2.
\end{eqnarray}
Taking expectation of \eqref{sa_eq3} over $\xi_t$ and $z_t$, we have
\begin{eqnarray}
\nonumber
\mathbb{E}_t(x_{t+1}-x^*)^2&=&(x_t-x^*)^2-2\eta_tH'(x_t)(x_t-x^*)+\eta_t^2\mathbb{E}_t(g(x_t,
\xi_t, z_t))^2\\\label{sa_eq4}
&\leq&(1-\eta_t\mu)(x_t-x^*)^2-2\eta_t(H(x_t)-H(x^*))+\eta_t^2\mathbb{E}_t(g(x_t,
\xi_t, z_t))^2,
\end{eqnarray}
where the first equality is because of Proposition
\ref{prop:unbiased} and the last inequality is because of
\eqref{eq:strconv} ($H(\cdot)$ is $\mu$-convex). According to the
third statement in Proposition \ref{prop:unbiased}, \eqref{sa_eq4}
implies that
\begin{eqnarray}
\label{sa_eq5}
H(x_t)-H(x^*)&\leq&\frac{1-\eta_t\mu}{2\eta_t}(x_t-x^*)^2-\frac{1}{2\eta_t}\mathbb{E}_t(x_{t+1}-x^*)^2+\frac{\eta_tG^2}{2}.
\end{eqnarray}

If $\mu=0$ and we choose $\eta_t=\frac{1}{\sqrt{T}}$, then summing
\eqref{sa_eq5} for $t=1,2,\dots, T$ and taking expectation give
\begin{eqnarray}
\nonumber T\mathbb{E}(H(\bar
x_T)-H(x^*))\leq\sum_{t=1}^T\mathbb{E}(H(x_t)-H(x^*))\leq\frac{\sqrt{T}}{2}(x_1-x^*)^2+\frac{\sqrt{T}G^2}{2}.
\end{eqnarray}
The desired result for this part is obtained by dividing this
inequality by $T$. In addition, if both $u$ and $\ell$ are finite
and we choose $\eta_t=\frac{1}{\sqrt{t}}$, then summing
\eqref{sa_eq5} for $t=1,2,\dots, T$ and taking expectation give
\begin{eqnarray}
\nonumber T\mathbb{E}(H(\bar
x_T)-H(x^*))\leq\sum_{t=1}^T\mathbb{E}(H(x_t)-H(x^*))\leq\sum_{t=1}^T\frac{\sqrt{t}-\sqrt{t-1}}{2}(x_t-x^*)^2+\sum_{t=1}^T\frac{G^2}{2\sqrt{t}}.
\end{eqnarray}
Note that $(x_t-x^*)^2\leq(u-\ell)^2$ and
$\sum_{t=1}^T\frac{1}{\sqrt{t}}\leq\int_0^T\frac{1}{\sqrt{x}}dx=2\sqrt{T}$
so that the above inequality implies
\begin{eqnarray}
\nonumber T\mathbb{E}(H(\bar
x_T)-H(x^*))\leq\sum_{t=1}^T\mathbb{E}(H(x_t)-H(x^*))\leq\frac{\sqrt{T}}{2}(u-\ell)^2+\sqrt{T}G^2.
\end{eqnarray}
The desired result for this part is obtained by dividing this
inequality by $T$.

\noindent 2. When Assumption A5(b) holds (based on
\citealt{Lan:10a}):

The $\mu$-convexity property \eqref{eq:strconv} of $H(\cdot)$
implies that
\begin{eqnarray}
\nonumber H(x^*)&\geq&
H(x_t)+H'(x_t)(x^*-x_t)+\frac{\mu}{2}(x^*-x_t)^2\\\nonumber &=&
H(x_t)+\mathbb{E}_t[g(x_t, \xi_t,
z_t)(x^*-x_t)]+\frac{\mu}{2}(x^*-x_t)^2\\\nonumber &=&
H(x_t)+\mathbb{E}_t[(g(x_t, \xi_t,
z_t)-H'(x_t))(x_{t+1}-x_t)]+H'(x_t)(x_{t+1}-x_t)\\\label{sa_eq-1}
&&+\mathbb{E}_t[g(x_t, \xi_t,
z_t)(x^*-x_{t+1})]+\frac{\mu}{2}(x^*-x_t)^2,
\end{eqnarray}
where the first equality is because of Proposition
\ref{prop:unbiased}. By Assumption A5(b), $H'(x)$ is $L$-Lipschitz
continuous, thus
\begin{eqnarray}
\nonumber
H(x_{t+1})&\leq&H(x_t)+H'(x_t)(x_{t+1}-x_t)+\frac{L}{2}(x_{t+1}-x_t)^2
\end{eqnarray}
which, together with \eqref{sa_eq-1}, implies that
\begin{eqnarray*}
    \nonumber H(x^*)&\geq& \mathbb{E}_tH(x_{t+1})+\mathbb{E}_t[(g(x_t,
    \xi_t,
    z_t)-H'(x_t))(x_{t+1}-x_t)]-\frac{L}{2}\mathbb{E}_t(x_{t+1}-x_t)^2\\\nonumber
    &&+\mathbb{E}_t[g(x_t, \xi_t,
    z_t)(x^*-x_{t+1})]+\frac{\mu}{2}(x^*-x_t)^2\\\nonumber
    &\geq&\mathbb{E}_tH(x_{t+1})-\frac{1}{2a_t}\mathbb{E}_t(g(x_t,
    \xi_t,
    z_t)-H'(x_t))^2-\frac{L+a_t}{2}\mathbb{E}_t(x_{t+1}-x_t)^2\\\nonumber
    &&+\mathbb{E}_t[g(x_t, \xi_t,
    z_t)(x^*-x_{t+1})]+\frac{\mu}{2}(x^*-x_t)^2\\\label{sa_eq-3}
    &\ge&\mathbb{E}_tH(x_{t+1})-\frac{\sigma^2}{2a_t}-\frac{L+a_t}{2}\mathbb{E}_t(x_{t+1}-x_t)^2+\mathbb{E}_t[g(x_t,
    \xi_t, z_t)(x^*-x_{t+1})]+\frac{\mu}{2}(x^*-x_t)^2,
\end{eqnarray*}
where $a_t$ is a positive constant, the second inequality is due to
Young's inequality, namely, $xy\leq\frac{x^2}{2a}+\frac{ay^2}{2}$
for any $a>0$, and the last equality is due to the third statement
of Proposition \ref{prop:unbiased}. We will determine the value of
$a_t$ later but we always ensure that $a_t$ and $\eta_t$ satisfy
\begin{eqnarray}
\label{sa_stepsize} \frac{L+a_t}{2}-\frac{1}{2\eta_t}\leq0.
\end{eqnarray}

By the optimality of $x_{t+1}$ as a solution to the projection
problem~\eqref{eq:proj}, we have
\begin{eqnarray*}
    \label{sa_eq-2} (x_{t+1}-x_t+\eta_tg(x_t, \xi_t,
    z_t))(x^*-x_{t+1})\geq0.
\end{eqnarray*}
Thus we have
\begin{eqnarray}
\nonumber H(x^*)&\geq&
\mathbb{E}_tH(x_{t+1})-\frac{\sigma^2}{2a_t}-\frac{L+a_t}{2}\mathbb{E}_t(x_{t+1}-x_t)^2+\frac{1}{\eta_t}\mathbb{E}_t[(x_{t+1}-x_t)(x_{t+1}-x^*)]+\frac{\mu}{2}(x^*-x_t)^2\\\nonumber
&=&\mathbb{E}_tH(x_{t+1})-\frac{\sigma^2}{2a_t}-\left(\frac{L+a_t}{2}-\frac{1}{2\eta_t}\right)\mathbb{E}_t(x_{t+1}-x_t)^2\\\nonumber
&&+\frac{1}{2\eta_t}\mathbb{E}_t(x_{t+1}-x^*)^2-\left(\frac{1}{2\eta_t}-\frac{\mu}{2}\right)\mathbb{E}_t(x_t-x^*)^2\\\label{sa_eq-4}
&\geq&\mathbb{E}_tH(x_{t+1})-\frac{\sigma^2}{2a_t}+\frac{1}{2\eta_t}\mathbb{E}_t(x_{t+1}-x^*)^2-\left(\frac{1}{2\eta_t}-\frac{\mu}{2}\right)\mathbb{E}_t(x_t-x^*)^2,
\end{eqnarray}
where the second inequality is from \eqref{sa_stepsize}.

If $\mu=0$ and we choose $\eta_t=\frac{1}{L+\sqrt{T}}$ and
$a_t=\sqrt{T}$ so that \eqref{sa_stepsize} is satisfied. Summing
\eqref{sa_eq-4} for $t=1,2,\dots, T-1$ and organizing terms give
\begin{eqnarray}
\nonumber T\mathbb{E}(H(\bar
x_T)-H(x^*))\leq\sum_{t=1}^T\mathbb{E}(H(x_t)-H(x^*))\leq
\frac{L+\sqrt{T}}{2}(x_1-x^*)^2+H(x_1)-H(x^*)+\frac{\sqrt{T}\sigma^2}{2}.
\end{eqnarray}
The desired result for this part is obtained by dividing this
inequality by $T$.  In addition, if both $u$ and $\ell$ are finite
and we choose $\eta_t=\frac{1}{L+\sqrt{t}}$ and $a_t=\sqrt{t}$, then
summing \eqref{sa_eq-4} for $t=1,2,\dots, T-1$ and taking
expectation give
\begin{eqnarray}
\nonumber T\mathbb{E}(H(\bar x_T)-H(x^*))& \leq
&\sum_{t=1}^T\mathbb{E}(H(x_t)-H(x^*))\\\nonumber
&\leq&\sum_{t=1}^{T-1}\frac{\sqrt{t}-\sqrt{t-1}}{2}(x_t-x^*)^2+
H(x_1)-H(x^*)+\frac{L}{2}(x_1-x^*)^2+\sum_{t=1}^{T-1}\frac{\sigma^2}{2\sqrt{t}}.
\end{eqnarray}
Note that $(x_t-x^*)^2\leq(u-\ell)^2$ and
$\sum_{t=1}^{T-1}\frac{1}{\sqrt{t}}\leq\int_0^T\frac{1}{\sqrt{x}}dx=2\sqrt{T}$
so that the above inequality implies
\begin{eqnarray}
\nonumber T\mathbb{E}(H(\bar
x_T)-H(x^*))\leq\sum_{t=1}^T\mathbb{E}(H(x_t)-H(x^*))\leq\frac{\sqrt{T}}{2}(u-\ell)^2+H(x_1)-H(x^*)+\frac{L}{2}(x_1-x^*)^2+\sqrt{T}\sigma^2.
\end{eqnarray}
The desired result for this part is obtained by dividing this
inequality by $T$. $\hfill\Box$\\

{\noindent\bf Proof of Proposition \ref{thm:CBA1sc}.}
Similar as to the proof of Proposition \ref{thm:CBA1}, we consider the
case when Assumption A5(a) or A5(b) holds respectively, and use
$\mathbb{E}_t$ to denote the conditional expectation taken over
$\xi_t$ and $z_t$ conditioning on
$\xi_1,z_1,\xi_2,z_2,\dots,\xi_{t-1},z_{t-1}$.

\noindent 1. When Assumption A5(a) holds:
We can still show \eqref{sa_eq5} using exactly
the same argument in the proof of Proposition \ref{thm:CBA1}.

If $\mu>0$ and we choose $\eta_t=\frac{1}{\mu t}$, summing
\eqref{sa_eq5} for $t=1,2,\dots, T$ gives
\begin{eqnarray*}
    T\mathbb{E}(H(\bar
    x_T)-H(x^*))\leq\sum_{t=1}^T\mathbb{E}(H(x_t)-H(x^*))\leq\sum_{t=1}^T\frac{G^2}{2\mu
        t}.
\end{eqnarray*}
The desired result for this part is obtained by dividing this
inequality by $T$ and using the fact that
$\sum_{t=1}^T\frac{1}{t} \le \log T+1$.\\

\noindent 2. When Assumption A5(b) holds: Using exactly
the same argument in the proof of Proposition \ref{thm:CBA1},
we can show that \eqref{sa_eq-4} holds for any positive constant $a_t$ that satisfies \eqref{sa_stepsize}.

If $\mu>0$ and we choose $\eta_t=\frac{1}{\mu t+L}$ and $a_t=\mu t$
so that \eqref{sa_stepsize} is satisfied. Summing \eqref{sa_eq-4}
for $t=1,2,\dots, T-1$ gives
\begin{eqnarray}
\nonumber T\mathbb{E}(H(\bar
x_T)-H(x^*))\leq\sum_{t=1}^T\mathbb{E}(H(x_t)-H(x^*))\leq\sum_{t=1}^T\frac{\sigma^2}{2\mu
    t}+H(x_1)-H(x^*)+\frac{L(x_1-x^*)^2}{2}.
\end{eqnarray}
The desired result for this part is obtained by dividing this
inequality by $T$ and using the fact that
$\sum_{t=1}^T\frac{1}{t}\le\log T+1$. $\hfill\Box$\\

Before we prove Proposition \ref{thm:restart}, we first introduce the
following lemma:

\begin{lemma}
    \label{lemma1} If Assumption A5(a) holds, then by choosing
    $\eta_t=\eta\in(0,+\infty)$, the CBA ensures that
    \begin{eqnarray}
    \nonumber \mathbb{E}(H(\bar x_T)-H(x^*))\leq\frac{(x_1-x^*)^2}{2\eta
        T}+\frac{\eta G^2}{2}.
    \end{eqnarray}
    If Assumption A5(b) holds, then by choosing
    $\eta_t=\eta\in(0,\frac{1}{L})$ and $a_t = \frac{1}{\eta_t}-L$, the CBA
    ensures that
    \begin{eqnarray}
    \nonumber
    \mathbb{E}(H(\bar x_T)-H(x^*))\leq
    \frac{(x_1-x^*)^2}{2\eta T}+\frac{H(x_1)-H(x^*)}{T}+\frac{\sigma^2}{1/\eta-L}.
    \end{eqnarray}
\end{lemma}

{\noindent\bf Proof of Lemma \ref{lemma1}.} When Assumption A5(a)
holds, by choosing $\eta_t = \eta\in (0,
\infty)$ and summing \eqref{sa_eq5} over $t=1,2,\dots, T$, we have
\begin{eqnarray*}
    T\mathbb{E}(H(\bar
    x_T)-H(x^*))\leq\sum_{t=1}^T\mathbb{E}(H(x_t)-H(x^*))\leq\frac{(x_1-x^*)^2}{2\eta}+\frac{T\eta
        G^2}{2}.
\end{eqnarray*}
The first conclusion is obtained by dividing this inequality by $T$.

When Assumption A5(b) holds, by setting
$\eta_t=\eta\in(0,\frac{1}{L})$, $a_t=\frac{1}{\eta}-L$ and summing
\eqref{sa_eq-4} over $t=1,2,\dots, T-1$, we have
\begin{eqnarray*}
    T\mathbb{E}(H(\bar x_T)-H(x^*))\leq
    \frac{(x_1-x^*)^2}{2\eta}+H(x_1)-H(x^*)+\frac{T\sigma^2}{1/\eta-L}.
\end{eqnarray*}
The second conclusion is obtained by dividing this inequality by
$T$. $\hfill\Box$\\

{\noindent\bf Proof of Proposition \ref{thm:restart}.} The optimality of
$x^*$ and the $\mu$-convexity property \eqref{eq:strconv} of
$H(\cdot)$ imply
\begin{eqnarray}
\label{eq:scv_lb} \frac{\mu}{2}(\hat x^k-x^*)^2\leq H(\hat
x_k)-H(x^*).
\end{eqnarray}
With a slight abuse of notation, let $\mathbb{E}_k$ be the
conditional expectation conditioning on $\hat x^1,\hat
x^2,\dots,\hat x^k$.

\noindent 1. When Assumption A5(a) holds (based on
\citealt{Hazan:14}):

Define $\Delta_k=H(\hat x^k)-H(x^*)$ for $k\geq 1$. In the
following, we use induction to show $\mathbb{E}\Delta_k \leq
\frac{\Delta_1+G^2/\mu}{2^{k-1}}$ for $k\geq 1$. Note that this
statement holds trivially when $k=1$. Suppose
$\mathbb{E}\Delta_k\leq \frac{\Delta_1+G^2/\mu}{2^{k-1}}$. Now we
consider $\mathbb{E}\Delta_{k+1}$.

By Lemma~\ref{lemma1} and $\eta_t^k=\frac{1}{2^{k+1}\mu}$, we have
\begin{eqnarray*}
    \mathbb{E}\Delta_{k+1} = \mathbb{E}_k(H(\hat x^{k+1})-H(x^*))\leq
    \frac{2^k\mu(\hat{x}^k-x^*)^2}{ T_k}+\frac{G^2}{2^{k+2}\mu}
    \leq\frac{2^{k+1}\Delta_k}{T_k}+\frac{G^2}{2^{k+2}\mu},
\end{eqnarray*}
where the second inequality is due to \eqref{eq:scv_lb}. Taking
expectation over $\hat x^1,\hat x^2,\dots,\hat x^k$ and applying the
induction assumption
$\mathbb{E}\Delta_k=\frac{\Delta_1+G^2/\mu}{2^{k-1}}$, we have
\begin{eqnarray}\label{sa_eq19}
\mathbb{E}\Delta_{k+1} \leq\frac{2^{k+1}(\Delta_1+G^2/\mu)}{
    T_k2^{k-1}}+\frac{G^2}{2^{k+2}\mu}
\leq\frac{\Delta_1}{2^{k+1}}+\frac{G^2}{2^{k+1}\mu}+\frac{G^2}{2^{k+2}\mu}
\leq\frac{\Delta_1+G^2/\mu}{2^k},
\end{eqnarray}
where we use the facts that $T_k=2^{k+3}$. By induction assumption,
$\mathbb{E}\Delta_k\leq\frac{\Delta_1+G^2/\mu}{2^{k-1}}$ for $k\geq
1$. Since $T=\sum_{k=1}^KT_k=\sum_{k=1}^K2^{k+3}\leq2^{K+4}$, we
have $K\geq\log_2(T/16)$. Let $k=K$ in \eqref{sa_eq19}, we have
\begin{eqnarray*}
    \mathbb{E}(H(\hat x^{K+1})-H(x^*))=\mathbb{E}\Delta_{K+1}
    \leq\frac{\Delta_1+G^2/\mu}{2^K} \le \frac{16(\Delta_1+G^2/\mu)}{T}.
\end{eqnarray*}

\noindent 2. When Assumption A5(b) holds:

Define $\Delta_k=H(\hat x^k)-H(x^*)+\frac{L}{2}(\hat x^k-x^*)^2$ for
$k\geq 1$. In the following, we use induction to show
$\mathbb{E}\Delta_k \leq \frac{\Delta_1+\sigma^2/\mu}{2^{k-1}}$ for
$k\geq 1$. Note that this statement holds trivially when $k=1$.
Suppose $\mathbb{E}\Delta_k\leq
\frac{\Delta_1+\sigma^2/\mu}{2^{k-1}}$. Now we consider
$\mathbb{E}\Delta_{k+1}$.

By Lemma~\ref{lemma1} and
$\eta_t^k=\frac{1}{2^{k+1}\mu+L}\in(0,\frac{1}{L})$, we have
\begin{eqnarray*}
    \mathbb{E}\Delta_{k+1} = \mathbb{E}_k(H(\hat x^{k+1})-H(x^*))\leq
    \frac{(2^{k+1}\mu+L)(\hat
        x^k-x^*)^2}{2T_k}+\frac{\Delta_k}{T_k}+\frac{\sigma^2}{2^{k+1}\mu}
    \leq\left(2^{k+1}+1\right)\frac{\Delta_k}{T_k}+\frac{\sigma^2}{2^{k+1}\mu},
\end{eqnarray*}
where the second inequality is due to \eqref{eq:scv_lb}. Taking
expectation over $\hat x^1,\hat x^2,\dots,\hat x^k$ and applying the
induction assumption
$\mathbb{E}\Delta_k=\frac{\Delta_1+\sigma^2/\mu}{2^{k-1}}$, we have
\begin{eqnarray}\label{sa_eq16}
\mathbb{E}\Delta_{k+1}
\leq\left(2^{k+1}+1\right)\frac{\Delta_1+\sigma^2/\mu}{T_k2^{k-1}}+\frac{\sigma^2}{2^{k+1}\mu}
\leq\frac{\Delta_1}{2^{k+1}}+\frac{\sigma^2}{2^{k+1}\mu}+\frac{\sigma^2}{2^{k+1}\mu}
\leq\frac{\Delta_1+\sigma^2/\mu}{2^k},
\end{eqnarray}
where we use the facts that $T_k=2^{k+3}+4$. Thus,
$\mathbb{E}\Delta_k\leq\frac{\Delta_1+\sigma^2/\mu}{2^{k-1}}$ for
$k\geq 1$. Since
$T=\sum_{k=1}^KT_k=\sum_{k=1}^K\left(2^{k+3}+4\right)\leq2^{K+5}$,
we have $K\geq\log_2(T/32)$. Let $k=K$ in \eqref{sa_eq16}, we have
\begin{eqnarray*}
    \mathbb{E}(H(\hat x^{K+1})-H(x^*))=\mathbb{E}\Delta_{K+1}
    \leq\frac{\Delta_1+\sigma^2/\mu}{2^K}=\frac{32(\Delta_1+\sigma^2/\mu)}{T}.
\end{eqnarray*}
Thus the proposition is proved. $\hfill\Box$\\

\section{Proofs of Proposition \ref{thm:CBA1-QP} and
\ref{thm:restart-QP}} \label{sec:convproofQP}

{\noindent\bf Proof of Proposition~\ref{thm:CBA1-QP}.}
\label{appendix:proofoftheorem4-5} For the ease of notation, we
shall use $\mathbb{E}_t$ to denote the conditional expectation taken
over $\xi_t$, $u_t$ and $z_t$ conditioning on
$\xi_1,z_1,u_t,\xi_2,z_2,u_2\dots,\xi_{t-1},z_{t-1},u_{t-1}$. The
proof is similar to the proof of Proposition~\ref{thm:CBA1} when
Assumption A5(b) holds.

Since $\mu$ is the smallest eigenvalues of $Q$ and it is positive, $H(\cdot)$ defined in \eqref{eq:QPobj} is $\mu$-convex, which implies \eqref{sa_eq-1} using Proposition~\ref{prop:unbiased_multi_dimensional}.
Moreover, since $L$ is the largest eigenvalues of $Q$, $H'(x)$ is $L$-Lipschitz
continuous so that
\begin{eqnarray}
\nonumber
H(x_{t+1})&\leq&H(x_t)+H'(x_t)(x_{t+1}-x_t)+\frac{L}{2}(x_{t+1}-x_t)^2
\end{eqnarray}
which, by Proposition~\ref{prop:unbiased_multi_dimensional} and the same argument in the proof of Proposition~\ref{thm:CBA1}, implies that
\begin{eqnarray*}
    \nonumber H(x^*)&\geq& \mathbb{E}_tH(x_{t+1})-\frac{\sigma^2}{2a_t}-\frac{L+a_t}{2}\mathbb{E}_t(x_{t+1}-x_t)^2+\mathbb{E}_t[g(x_t,
    \xi_t, z_t)(x^*-x_{t+1})]+\frac{\mu}{2}(x^*-x_t)^2,
\end{eqnarray*}
where $a_t$ is a positive constant chosen to satisfy \eqref{sa_stepsize}.

By the optimality of $x_{t+1}$ as a solution to the projection
problem~\eqref{eq:projqp}, we have
\begin{eqnarray*}
    \label{sa_eq-2} (x_{t+1}-x_t+\eta_tg(x_t, \xi_t,
    z_t))(x^*-x_{t+1})\geq0.
\end{eqnarray*}
Using this inequality and the same argument in the proof of Proposition~\ref{thm:CBA1} (under Assumption A5(b)), we can show that \eqref{sa_eq-4} holds.

Since $\mu>0$ by assumption,  we can choose $\eta_t=\frac{1}{\mu t+L}$ and $a_t=\mu t$
so that \eqref{sa_stepsize} is satisfied. Summing \eqref{sa_eq-4}
for $t=1,2,\dots, T-1$ gives
\begin{eqnarray}
\nonumber T\mathbb{E}(H(\bar
x_T)-H(x^*))\leq\sum_{t=1}^T\mathbb{E}(H(x_t)-H(x^*))\leq\sum_{t=1}^T\frac{\sigma^2}{2\mu
    t}+H(x_1)-H(x^*)+\frac{L(x_1-x^*)^2}{2}.
\end{eqnarray}
The conclusion of this proposition is obtained by dividing this
inequality by $T$ and using the fact that
$\sum_{t=1}^T\frac{1}{t}\le\log T+1$.  $\hfill\Box$\\

{\noindent\bf Proof of Proposition~\ref{thm:restart-QP}.}
The optimality of
$x^*$ and the $\mu$-convexity  of
$H(\cdot)$ imply \eqref{eq:scv_lb}.  Let $\mathbb{E}_k$ be the
conditional expectation
conditioning on $\hat x^1,\hat x^2,\dots,\hat x^k$.\\

Following a similar argument to the proof of Proposition~\ref{thm:restart} when Assumption A5(b) holds, we
define $\Delta_k=H(\hat x^k)-H(x^*)+\frac{L}{2}(\hat x^k-x^*)^2$ for
$k\geq 1$. We then use induction to show
$\mathbb{E}\Delta_k \leq \frac{\Delta_1+\sigma^2/\mu}{2^{k-1}}$ for
$k\geq 1$. Note that this statement holds trivially when $k=1$.
Suppose $\mathbb{E}\Delta_k\leq
\frac{\Delta_1+\sigma^2/\mu}{2^{k-1}}$. Now we consider
$\mathbb{E}\Delta_{k+1}$.

By setting $\eta_t^k=\frac{1}{2^{k+1}\mu+L}\in(0,\frac{1}{L})$ and summing
\eqref{sa_eq-4} over $t=1,2,\dots, T-1$, we have
\begin{eqnarray*}
    \mathbb{E}\Delta_{k+1} = \mathbb{E}_k(H(\hat x^{k+1})-H(x^*))\leq
    \frac{(2^{k+1}\mu+L)(\hat
        x^k-x^*)^2}{2T_k}+\frac{\Delta_k}{T_k}+\frac{\sigma^2}{2^{k+1}\mu}
    \leq\left(2^{k+1}+1\right)\frac{\Delta_k}{T_k}+\frac{\sigma^2}{2^{k+1}\mu},
\end{eqnarray*}
 Taking
expectation over $\hat x^1,\hat x^2,\dots,\hat x^k$ and applying the
induction assumption
$\mathbb{E}\Delta_k=\frac{\Delta_1+\sigma^2/\mu}{2^{k-1}}$, we can obtain \eqref{sa_eq16}. Thus,
$\mathbb{E}\Delta_k\leq\frac{\Delta_1+\sigma^2/\mu}{2^{k-1}}$ for
$k\geq 1$. Since
$T=\sum_{k=1}^KT_k=\sum_{k=1}^K\left(2^{k+3}+4\right)\leq2^{K+5}$,
we have $K\geq\log_2(T/32)$. Let $k=K$ in \eqref{sa_eq16}, we have
\begin{eqnarray*}
    \mathbb{E}(H(\hat x^{K+1})-H(x^*))=\mathbb{E}\Delta_{K+1}
    \leq\frac{\Delta_1+\sigma^2/\mu}{2^K}=\frac{32(\Delta_1+\sigma^2/\mu)}{T}.
\end{eqnarray*}
Thus the proposition is proved. $\hfill\Box$
$\hfill\Box$\\

\section{CBA with Categorical Results from Comparison}
\label{appendix:multilevelCBA}

    In CBA given in Algorithm~\ref{alg:cba1}, the results of the comparison between $x$ and $\xi$ is binary, i.e., either $\xi>x$ or $\xi<x$. In this section, we extend CBA to allow the comparison result to be categorical and depend on the gap between $\xi$ and $x$. In particular, we assume there exist $m+1$ non-negative quantities $\theta_0,\theta_1,\dots,\theta_{m-1}$ and $\theta_m$ satisfying $0=:\theta_0<\theta_1<\dots<\theta_{m-2}<\theta_{m-1}<\theta_m:=+\infty$  such that, after presenting a solution $x$, we know whether
    $\xi\in (x+\theta_i,x+\theta_{i+1}]$ or $\xi\in [x-\theta_{i+1},x-\theta_i)$ for all $i=0,1,\dots,m-1$. Using Example~\ref{example:mean} as an example, this type of comparison result corresponds to the case where the customer reports a coarse level of the difference between his/her preferred value of
    the feature and the actual value of the feature of the product
    presented to him/her (e.g., whether the size of the product is too small, a little small, a little large or too large).  Note that the binary comparison result is a special case of the categorical result with $m=1$.

    To facilitate the development of algorithm, we need to replace  Assumption
    A2 with the following assumption:
    \begin{itemize}
        \item[](A2') For each $\xi$, $h(x,\xi)$ is continuously differentiable with respect to $x$ on $[\ell, \xi)$ and
        $(\xi, u]$ with the derivative denoted by $h'_x(x,\xi)$.
        Furthermore, for any $x\in [\ell, u]$, $h'_{-}(x,x-\theta_i) :=
        \underset{z\rightarrow x-\theta_i}{\lim} h'_x(x, z)$ and $h'_{+}(x,x+\theta_i) :=
        \underset{z\rightarrow x+\theta_i}{\lim} h'_x(x, z)$ exist and are finite for any $i=0,1,\dots,m-1$.
    \end{itemize}

    With this comparative   information,  we propose a comparison-based algorithm with categorical comparison result (CBA-C) to
    solve (\ref{formulation}).
    The algorithm requires specification of $2m$ functions,
    $f_{-}^i(x, z)$ and $f_{+}^i(x, z)$ for $i=0,1,\dots,m-1$, which need to satisfy the following
    conditions.
    \begin{itemize}
        \item (C1') $f_{-}^i(x,z) = 0$ for all $z \notin [x-\theta_{i+1},x-\theta_i)$ and $f_{-}^i(x,z) > 0$ for all $\max\{\underline{s},x-\theta_{i+1}\}\leq z\leq x-\theta_i$. In addition,  for all $x$, we have $\int_{x-\theta_{i+1}}^{x-\theta_i}
        f_{-}^i(x,z) dz = 1$.
        \item (C2') $f_{+}^i(x,z) =0$ for all $z \notin (x+\theta_i,x+\theta_{i+1}]$ and $f_{+}^i(x,z) > 0$ for all $\min\{\bar{s},x+\theta_{i+1}\} \ge z > x +\theta_i$. In addition, for all $x$, we have $\int_{x+\theta_i}^{x+\theta_{i+1}}
        f_{+}^i(x,z) dz = 1$.
        \item (C3') There exists a constant $K_3^i$ such that $\int_{x-\theta_{i+1}}^{x-\theta_i}\frac{(F(z)-F(x-\theta_{i+1}))(h''_{x,z}(x,z))^2}{f_{-}^i(x,z
            )}dz\leq K_3$ and
        $\int_{x-\theta_{i+1}}^{x-\theta_i}\frac{(F(x+\theta_{i+1})-F(z))(h''_{x,z}(x,z))^2}{f_{+}^i(x,z)}dz\leq
        K_3$ for all $x\in[\ell, u]$, where $F(\cdot)$ is the c.d.f. of $\xi$.
    \end{itemize}

    Note that $f_{-}^i(x,z)$ and $f_{+}^i(x,z)$ essentially define the
    density function on $[x-\theta_{i+1},x-\theta_i)$ and $(x+\theta_i,x+\theta_{i+1}]$, respectively, for any given
    $x\in[\ell, u]$. The functions $f_{-}^i(x,z)$ and $f_{+}^i(x,z)$ satisfying C1'-C3' can be constructed in a similar way as in
    Example~\ref{example:funiform} and \ref{example:fexponential} and their optimal choices are similar to \eqref{eq:optfneg} and \eqref{eq:optfpos}.
    Next, in Algorithm \ref{alg:cba1categorical}, we describe the detailed procedure of
    the CBA-C.

    \begin{algorithm}[h] \caption{Comparison-Based Algorithm with Categorical Comparison Result
            (CBA-C):\label{alg:cba1categorical}}
        \begin{enumerate}
            \item {\bf Initialization.} Set $t = 1$, $x_1 \in[\ell, u]$. Define $\eta_t$
            for all $t\ge 1$. Set the maximum number of iterations $T$.
            Choose functions $f_{-}^i(x,z)$ and $f_{+}^i(x,z)$ for $i=0,1,\dots,m-1$ that satisfy
            (C1')-(C3').
            \item {\bf Main iteration.} Sample $\xi_t$ from the distribution of $\xi$. If $\xi_t = x_t$, then resample
            $\xi_t$ until it does not equal $x_t$. (This step will always
            terminate in a finite number of steps as long as $\xi$ is not
            deterministic.)
            \begin{enumerate}
                \item If $\xi_t\in [x_t-\theta_{i+1},x_t-\theta_i)$, then generate $z_t$ from a
                distribution on $[x_t-\theta_{i+1},x_t-\theta_i)$ with p.d.f. $f_{-}^i(x_t,z_t)$. Set
                \begin{eqnarray}\label{gradient_case1-c}
                g(x_t, \xi_t, z_t) = \left\{\begin{array}{ll} h'_{-}(x_t,x_t-\theta_i), & \mbox{
                    if } z_t < \xi_t,
                \\
                h'_{-}(x_t,x_t-\theta_i) - \frac{h''_{x,z}(x_t,z_t)}{f_{-}(x_t,z_t)}, & \mbox{ if
                } z_t \ge \xi_t. \end{array}\right.
                \end{eqnarray}
                \item If $\xi_t\in (x_t+\theta_i,x_t+\theta_{i+1}]$, then generate $z_t$ from a
                distribution on $(x_t+\theta_i,x_t+\theta_{i+1}]$ with p.d.f. $f_{+}^i(x_t, z_t)$. Set
                \begin{eqnarray}\label{gradient_case2-c}
                g(x_t, \xi_t, z_t) = \left\{\begin{array}{ll} h'_{+}(x_t,x_t+\theta_i), & \mbox{
                    if } z_t > \xi_t,
                \\
                h'_{+}(x_t,x_t+\theta_i) + \frac{h''_{x,z}(x_t,z_t)}{f_{+}(x_t,z_t)}, & \mbox{ if
                } z_t \le \xi_t.\end{array} \right.
                \end{eqnarray}
            \end{enumerate}
            Let
            \begin{eqnarray} \label{eq:proj-c}
            x_{t+1}&=&\text{Proj}_{[\ell, u]}\left(x_t-\eta_tg(x_t, \xi_t,
            z_t)\right) = \max\left(\ell, \min\left(u, x_t - \eta_tg(x_t, \xi_t,
            z_t)\right)\right).
            \end{eqnarray}
            \item {\bf Termination.} Stop when $t\ge T$. Otherwise, let $t\leftarrow t+ 1$ and go
            back to Step 2.
            \item {\bf Output.} $\textbf{CBA-C}(x_1,T,\{\eta_t\}_{t=1}^T) = \bar x_T=\frac{1}{T}\sum_{t=1}^Tx_t$.
        \end{enumerate}
    \end{algorithm}

    We have the following proposition about the
    comparison-based algorithm with categorical result.

    \begin{proposition}\label{prop:unbiased-c}
        Suppose $f_{-}^i(x,z)$ and $f_{+}^i(x,z)$ satisfy (C1')-(C3') and
        Assumption \ref{assumption:onedimensional} holds with (A2) replaced by (A2'). Then
        \begin{enumerate}
            \item
            ${\mathbb E}_{z }g(x , \xi , z ) = h'_x(x ,\xi )$, for all $x\in
            [\ell, u]$, $x \neq \xi $.
            \item ${\mathbb E}_{z , \xi } g(x , \xi , z ) =
            H'(x )$, for all $x \in [\ell, u]$.
            \item If Assumption A5(a) holds, then $\mathbb{E}_{z,\xi} (g(x, \xi, z))^2\leq G'^2:= K_1^2+2\sum_{i=0}^{m-1}K_3^i$. If Assumption A5(b) holds, then ${\mathbb E}_{z,\xi} (g(x, \xi, z)-H'(x))^2\leq \sigma'^2:= K_2^2+2\sum_{i=0}^{m-1}K_3^i$.
        \end{enumerate}
    \end{proposition}

    {\bf\noindent Proof of Proposition \ref{prop:unbiased-c}.} First, we
    consider the case when $\xi  \in [x-\theta_{i+1},x-\theta_i) $. We have
    \begin{eqnarray*}
        {\mathbb E}_{z } g(x , \xi , z ) = h_{-}'(x,x -\theta_i) - \int_{\xi }^{x -\theta_i}
        h''_{x,z}(x , z)dz = h'_x(x , \xi).
    \end{eqnarray*}
    Similarly, when $\xi \in (x+\theta_i,x+\theta_{i+1}]$,
    \begin{eqnarray*} {\mathbb
            E}_{z } g(x , \xi , z ) = h_{+}'(x,x +\theta_i) + \int_{x +\theta_i}^{\xi } h''_{x,z}(x
        , z)dz =  h'_x(x , \xi).
    \end{eqnarray*}
    Thus the first conclusion of the proposition is proved. The second
    conclusion of the proposition follows from Assumption A1 (which
    ensures $\xi=x$ is a zero-measure event) and Assumption A4.

    Next, we show the first part of the third conclusion when Assumption
    A5(a) is true. If $\xi  \in [x-\theta_{i+1},x-\theta_i) $, then we have
    \small
    \begin{eqnarray}
    \label{eq:A1-c} {\mathbb E}_{z} (g(x, \xi, z))^2&=&(h'_{-}(x,x-\theta_i))^2 + \int_{\xi }^{x -\theta_i} \left(-
    2h'_{-}(x-\theta_i)\frac{h''_{x,z}(x,z)}{f_{-}^i(x ,z
        )}+\left(\frac{h''_{x,z}(x,z)}{f_{-}^i(x ,z
        )}\right)^2\right)f_{-}^i(x ,z )dz \nonumber \\
    &=&(h'_{-}(x,x-\theta_i))^2-2h'_{-}(x,x-\theta_i)(h'_{-}(x,x-\theta_i)-h'_x(x , \xi
    ))+ \int_{\xi }^{x -\theta_i}\frac{(h''_{x,z}(x,z))^2}{f_{-}^i(x ,z )}dz\\
    &\leq&(h'_x(x , \xi))^2+\int_{\xi }^{x
        -\theta_i}\frac{(h''_{x,z}(x,z))^2}{f_{-}^i(x ,z )}dz.\nonumber
    \end{eqnarray}
    \normalsize
    where the last inequality is because $a^2 + b^2 \ge 2ab$ for any
    $a$, $b$. By similar arguments, if $\xi\in (x+\theta_i,x+\theta_{i+1}]$, then
    \small
    \begin{eqnarray*}
        \label{eq:A2-c} {\mathbb E}_{z} (g(x, \xi, z))^2 &\leq&(h'_x(x ,
        \xi))^2+ \int_{x +\theta_i }^{\xi}\frac{(h''_{x,z}(x,z))^2}{f_{+}^i(x ,z )}dz.
    \end{eqnarray*}
    \normalsize
    These two inequalities and Assumption A5(a) further imply
    \small
    \begin{eqnarray}
    &&{\mathbb E}_{z,\xi} (g(x, \xi,
    z))^2\nonumber\\
    &\leq&K_1^2+\sum_{i=0}^{m-1}\int_{x-\theta_{i+1}}^{x-\theta_i}\left(\int_{\xi }^{x
        -\theta_i}\frac{(h''_{x,z}(x,z))^2}{f_{-}^i(x ,z )}dz\right)dF(\xi)+\sum_{i=0}^{m-1}\int_{x +\theta_i
    }^{x+\theta_{i+1}}\left(\int_{x + \theta_i}^{\xi}\frac{(h''_{x,z}(x,z))^2}{f_{+}^i(x
        ,z )}dz\right)dF(\xi)\nonumber \\
    &=&K_1^2+\sum_{i=0}^{m-1}\int_{x-\theta_{i+1}}^{x-\theta_i}\frac{(F(z)-F(x-\theta_{i+1}))(h''_{x,z}(x,z))^2}{f_{-}^i(x
        ,z )}dz+\sum_{i=0}^{m-1}\int_{x +\theta_i
    }^{x+\theta_{i+1}}\frac{(F(x+\theta_{i+1})-F(z))(h''_{x,z}(x,z))^2}{f_{+}^i(x ,z
        )}dz\nonumber \\ &\leq&K_1^2+2\sum_{i=0}^{m-1}K_3^i, \label{eq:A1A2-c}
    \end{eqnarray}
    \normalsize
    where the interchanging of integrals in the equality is justified by
    Tonelli's theorem and the last inequality is due to (C3').

    Next, we show the second part of the third conclusion when
    Assumption A5(b) is true. If $\xi  \in [x-\theta_{i+1},x-\theta_i) $, then following the similar
    analysis as in \eqref{eq:A1-c}, we have
    \begin{eqnarray*}
        {\mathbb E}_{z} (g(x, \xi, z)- h'_x(x ,\xi ))^2 ={\mathbb E}_{z}
        (g(x, \xi, z))^2- (h'_x(x ,\xi ))^2 \leq \int_{\xi }^{x
            -\theta_i}\frac{(h''_{x,z}(x,z))^2}{f_{-}^i(x ,z )}dz.
    \end{eqnarray*}
    Similarly, if $\xi\in (x+\theta_i,x+\theta_{i+1}]$, then
    \begin{eqnarray*}
        {\mathbb E}_{z} (g(x, \xi, z)- h'_x(x ,\xi ))^2&\leq& \int_{x +\theta_i }^{\xi}\frac{(h''_{x,z}(x,z))^2}{f_{+}^i(x ,z )}dz.
    \end{eqnarray*}
    By using the same argument as in \eqref{eq:A1A2-c}, we have
    \begin{eqnarray*}
        {\mathbb E}_{z, \xi} (g(x, \xi, z)- h'_x(x ,\xi ))^2 \le 2\sum_{i=0}^{m-1}K_3^i.
    \end{eqnarray*}
    Finally, we note that,
    \begin{eqnarray*}
        {\mathbb E}_{z,\xi} (g(x, \xi, z)-H'(x))^2&=&{\mathbb E}_{z,\xi}
        (g(x, \xi, z)- h'_x(x ,\xi ))^2+
        \mathbb{E}_{\xi}(h_x'(x,\xi)-H'(x))^2.
    \end{eqnarray*}
    Therefore, when Assumption A5(b) holds, we have ${\mathbb E}_{z,\xi}
    (g(x, \xi, z)-H'(x))^2 \le K_2^2 + 2\sum_{i=0}^{m-1}K_3^i$. Thus the proposition
    holds. $\hfill\Box$\\

    Proposition \ref{prop:unbiased-c} shows that in the CBA-C, the gradient
    estimate $g(x, \xi, z)$ is an unbiased estimate of the true gradient
    at $x$ and can be utilized as a stochastic gradient of $H(x)$. As a result, we can also
    prove that the convergence rates for CBA-C is exactly the same as those of CBA given in
    Proposition~\ref{thm:CBA1} and~\ref{thm:CBA1sc} except that $G^2$ and $\sigma^2$ are replaced by $G'^2$ and $\sigma'^2$.

   We also conduct some numerical experiments using CBA-C and the results are given in Appendix
    ~\ref{sec:expcategorical}.\\

    \section{Additional Numerical Tests}
\label{appendix:additional_numerical}

    \subsection{Additional Numerical Results for Mini-Batch Methods}
     \label{subsec:numerical_minibatch_add}
    
    In this section, we present more numerical results on the performance of CBA
    depends on the sample size $S$ in the mini-batch technique described
    in Section~\ref{sec:MB}. The results here are supplementary to the results presented in Section~\ref{subsec:numerical_minibatch}.

\begin{figure*}[!h]
        \centering
        \subfigure{
            \centering
            \includegraphics[width=0.45\columnwidth]{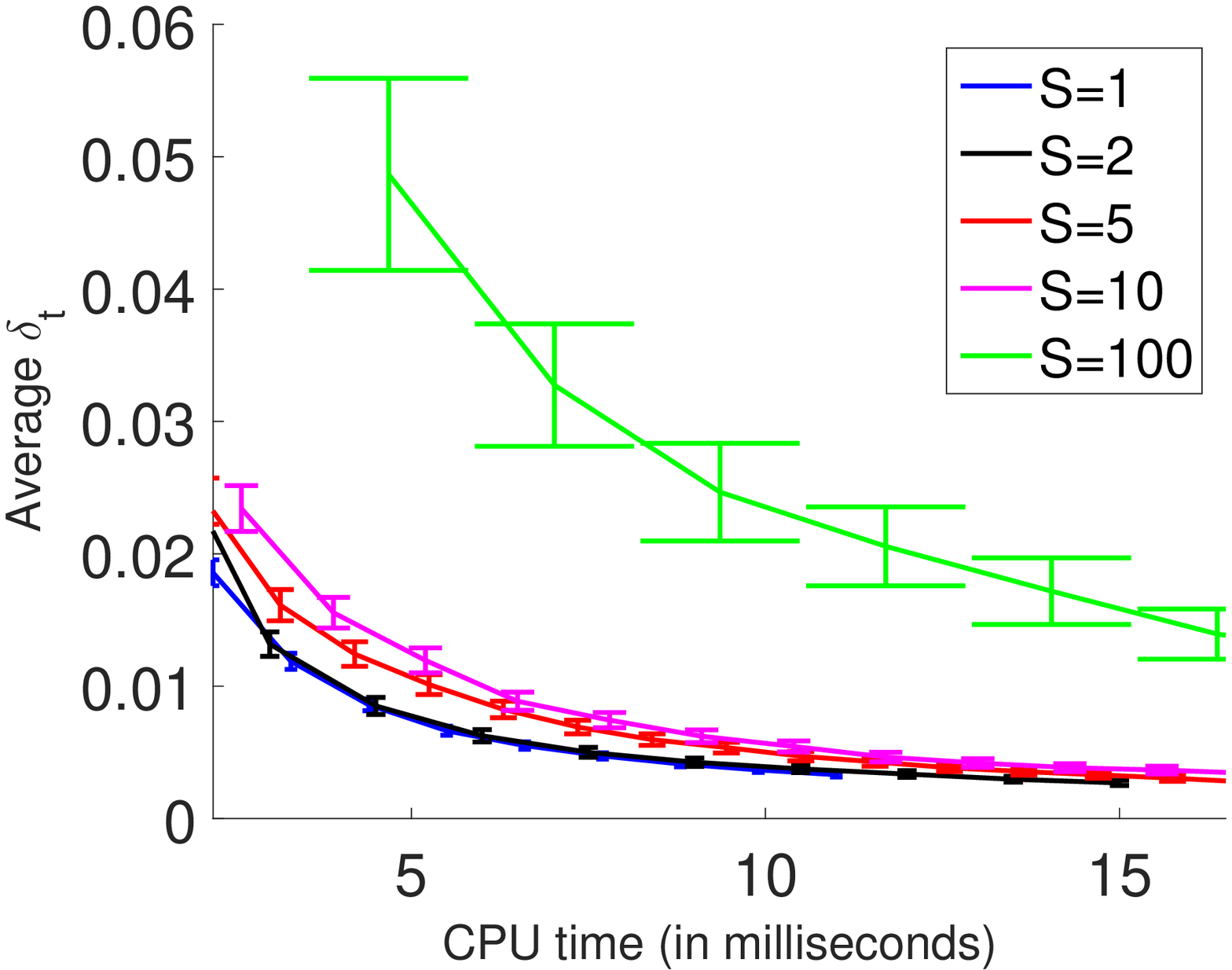}
            \label{fig:uniform_quad_uniform}
        }
        \subfigure{
            \centering
            \includegraphics[width=0.45\columnwidth]{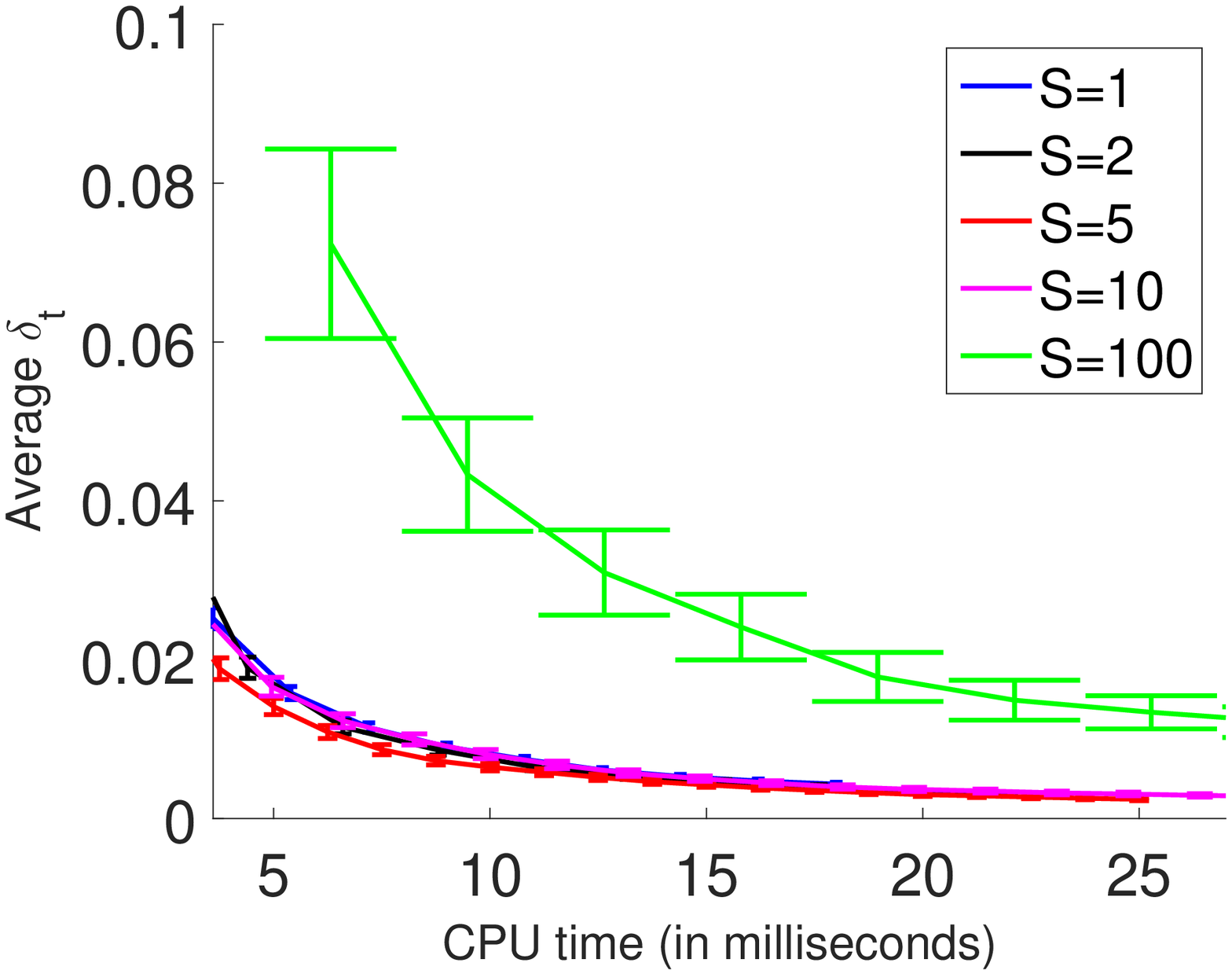}
            \label{fig:normal_quad_exp}
        }
        \subfigure{
            \centering
            \includegraphics[width=0.45\columnwidth]{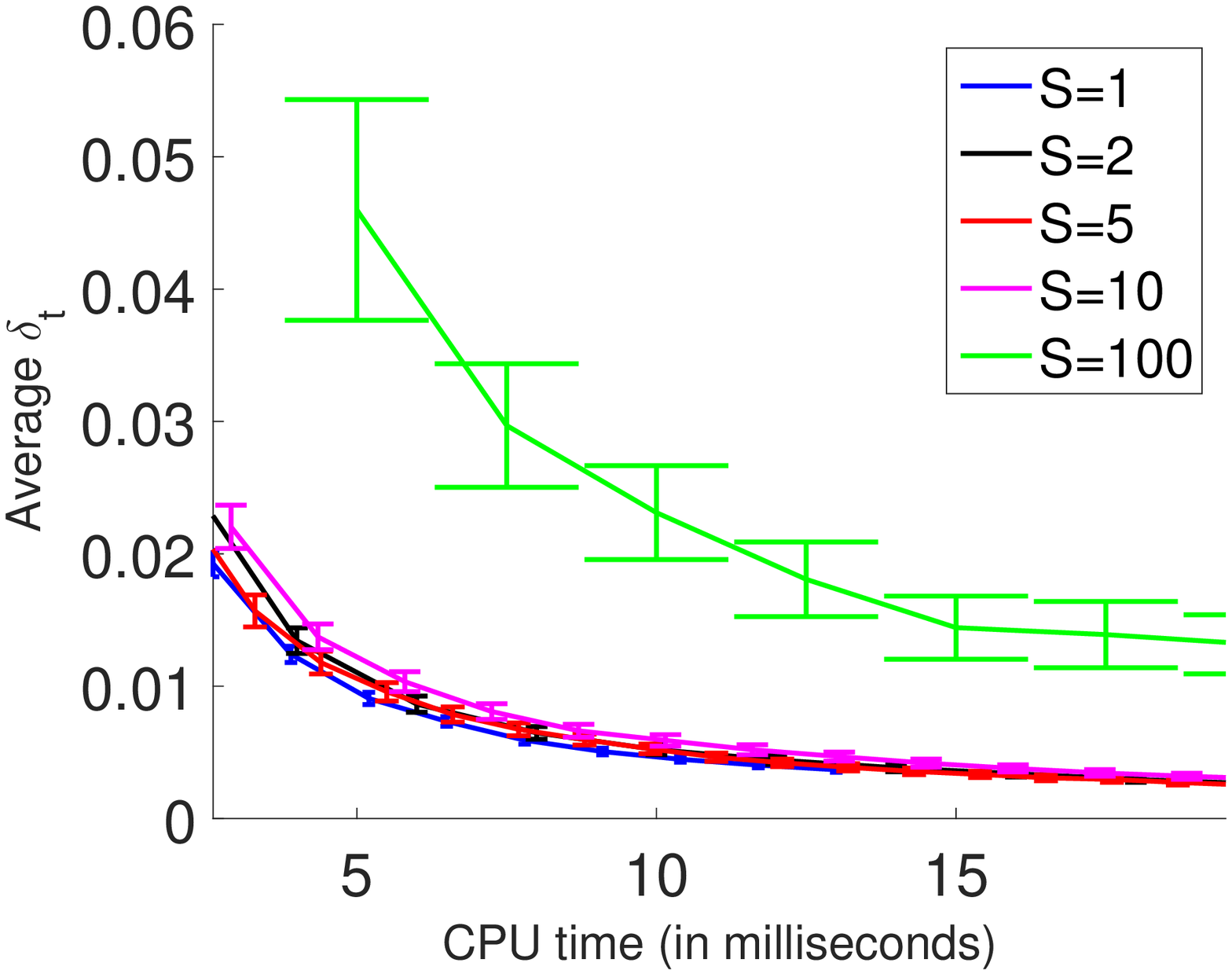}
            \label{fig:uniform_pquad_uniform}
        }
        \subfigure{
            \centering
            \includegraphics[width=0.45\columnwidth]{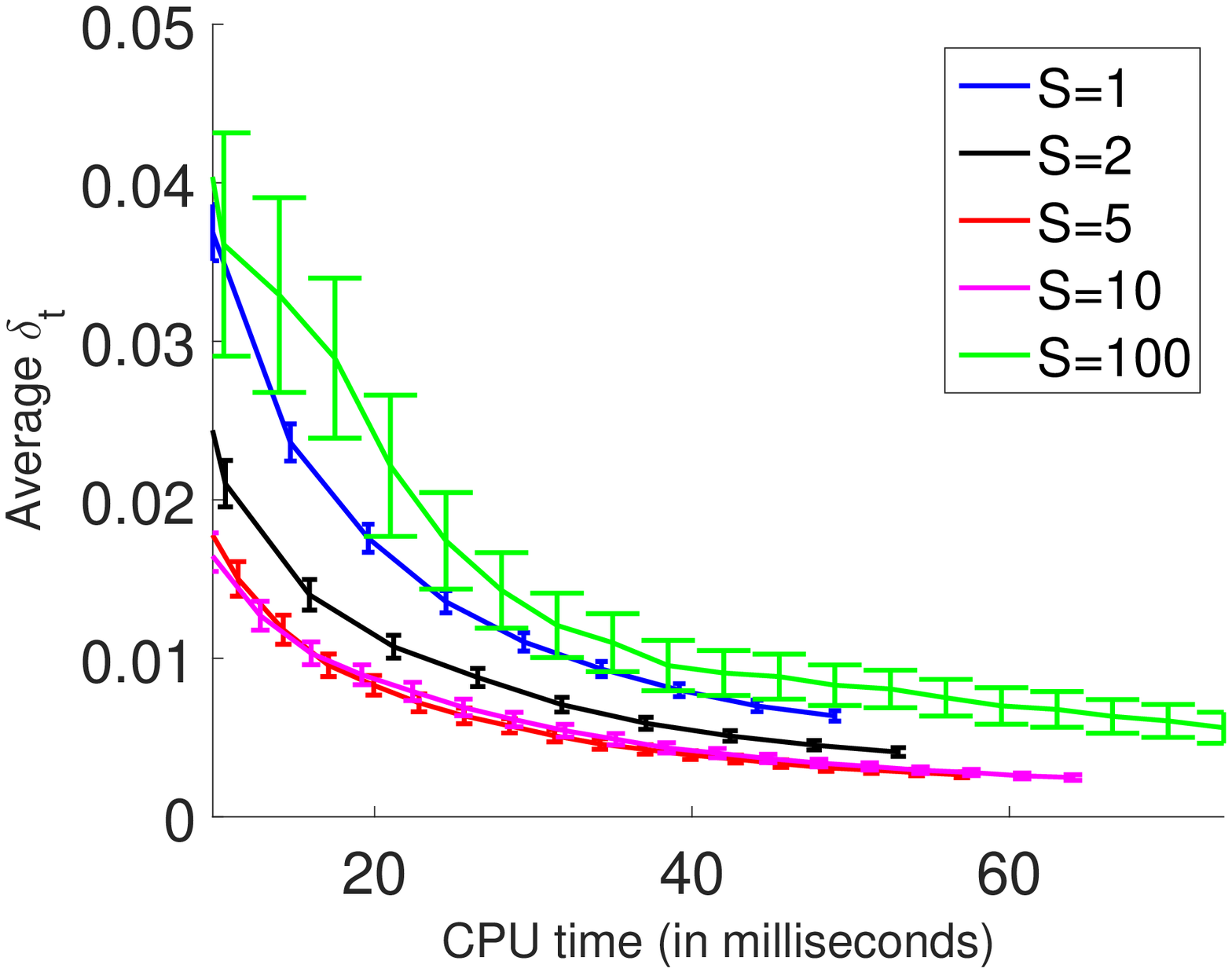}
            \label{fig:normal_pquad_exp}
        }
        \caption{The convergence of the average relative
            optimality gap in CBA when using different numbers of comparisons ($S$) in the mini-batch method.
            First row: $h_1(x,\xi)$; Second row: $h_2(x,\xi)$. First
            column: $\xi\sim\mathcal{U}[50,150]$;  Second column:
            $\xi\sim\mathcal{N}(100,100)$.}
        \label{fig:synthetic_MB_time}
    \end{figure*}

    With the same instances and setting as in Section~\ref{subsec:numerical_minibatch},
    we present the convergence of CBA with $S=1,2,5,10,100$ in
    Figure~\ref{fig:synthetic_MB_time}. Different from Figure~\ref{fig:synthetic_MB},
    the horizontal axis in Figure~\ref{fig:synthetic_MB_time} represents the CPU time elapsed.
    According to the figures, additional comparisons with $z$ do not always improve the time efficiency of CBA.
    In fact, the CBAs with a small batch size ($S=1$ or $2$) converges faster as time increases in the first column of Figure~\ref{fig:synthetic_MB} while
    the CBAs with a medium batch size ($S=5$ or $10$) converges faster in the second column.
    This does not contradict with the results in Figure~\ref{fig:synthetic_MB} because, within the same amount of time, the CBAs with a small $S$
    can run more iterations than the CBAs with a large $S$. Hence, the CBA with the largest batch size (i.e. $S=100$) converges most slowly due to its long runtime per iteration. In practice, one can experiment with difference values of $S$ to achieve the best time efficiency in CBA.

\begin{figure*}[!h]
    \centering
    \subfigure{
        \centering
        \includegraphics[width=0.45\columnwidth]{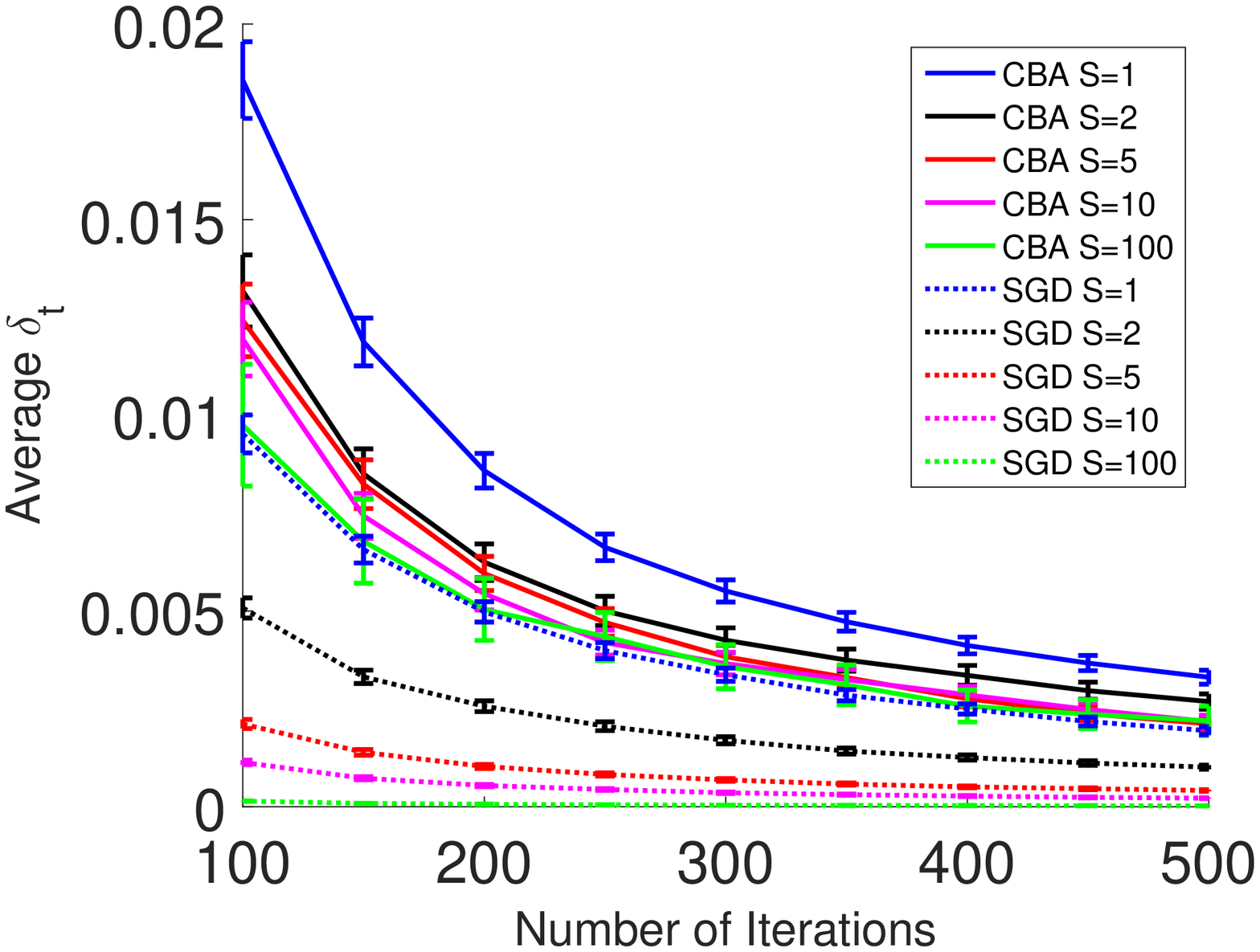}
        \label{fig:uniform_quad_uniform}
    }
    \subfigure{
        \centering
        \includegraphics[width=0.45\columnwidth]{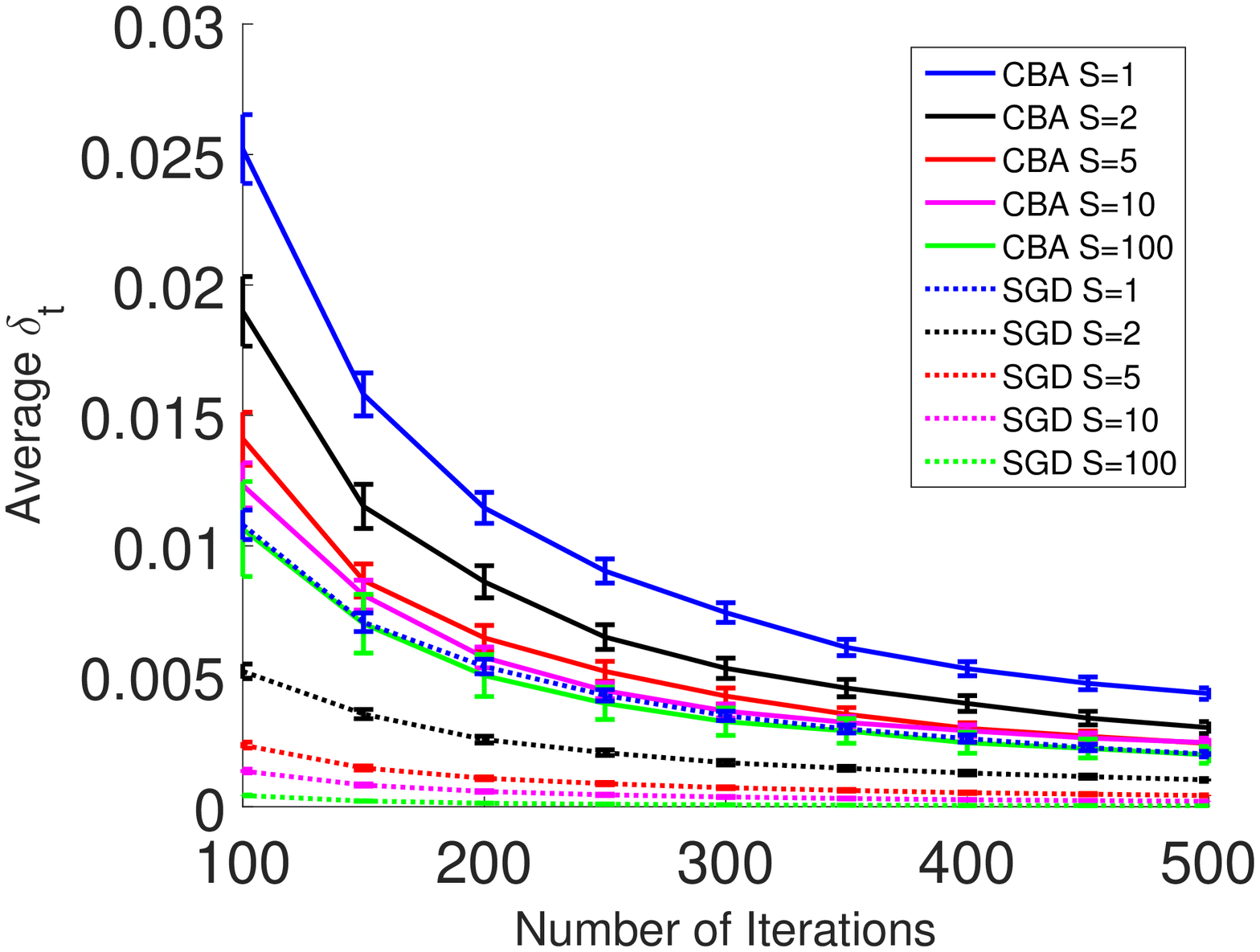}
        \label{fig:normal_quad_exp}
    }
    \subfigure{
        \centering
        \includegraphics[width=0.45\columnwidth]{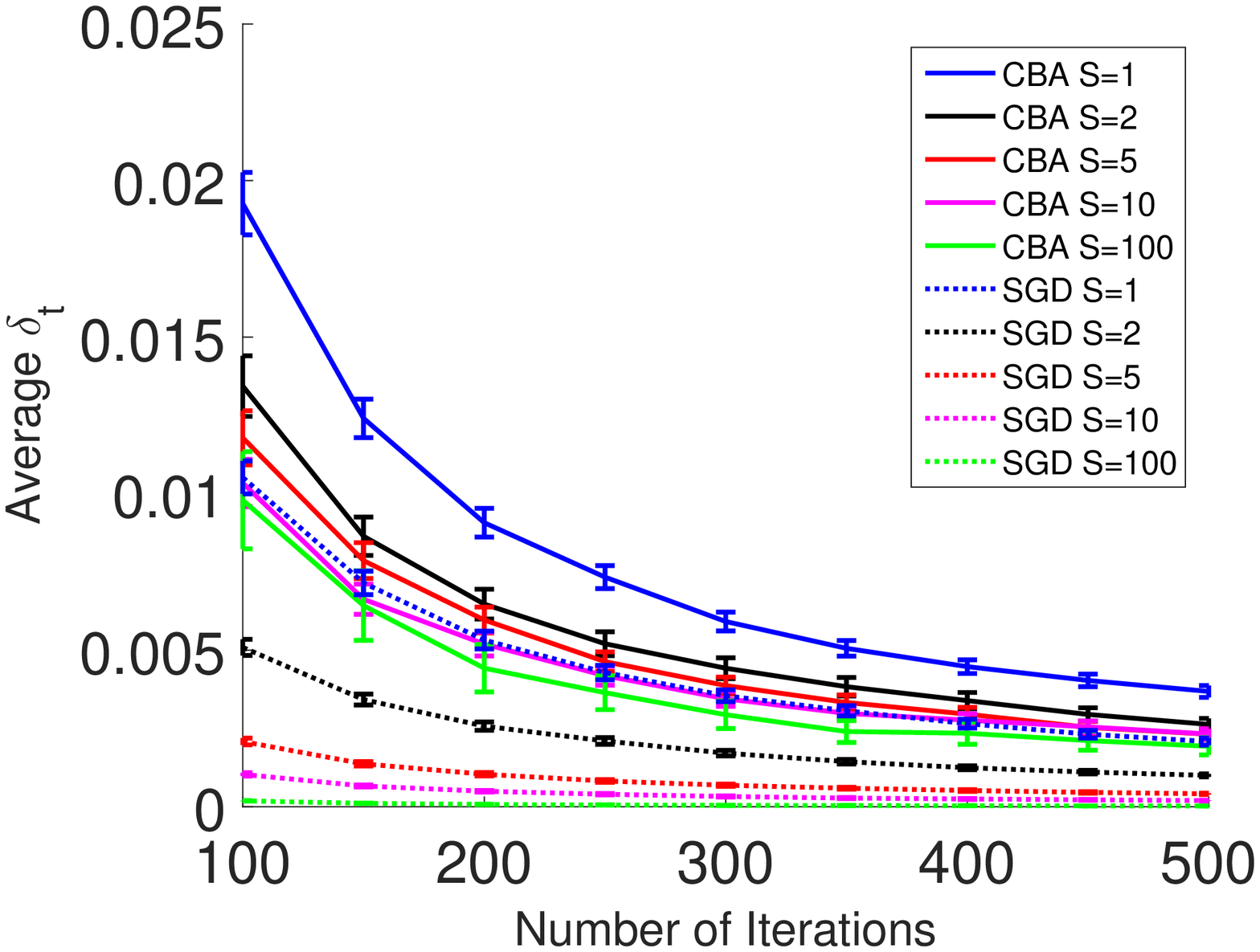}
        \label{fig:uniform_pquad_uniform}
    }
    \subfigure{
        \centering
        \includegraphics[width=0.45\columnwidth]{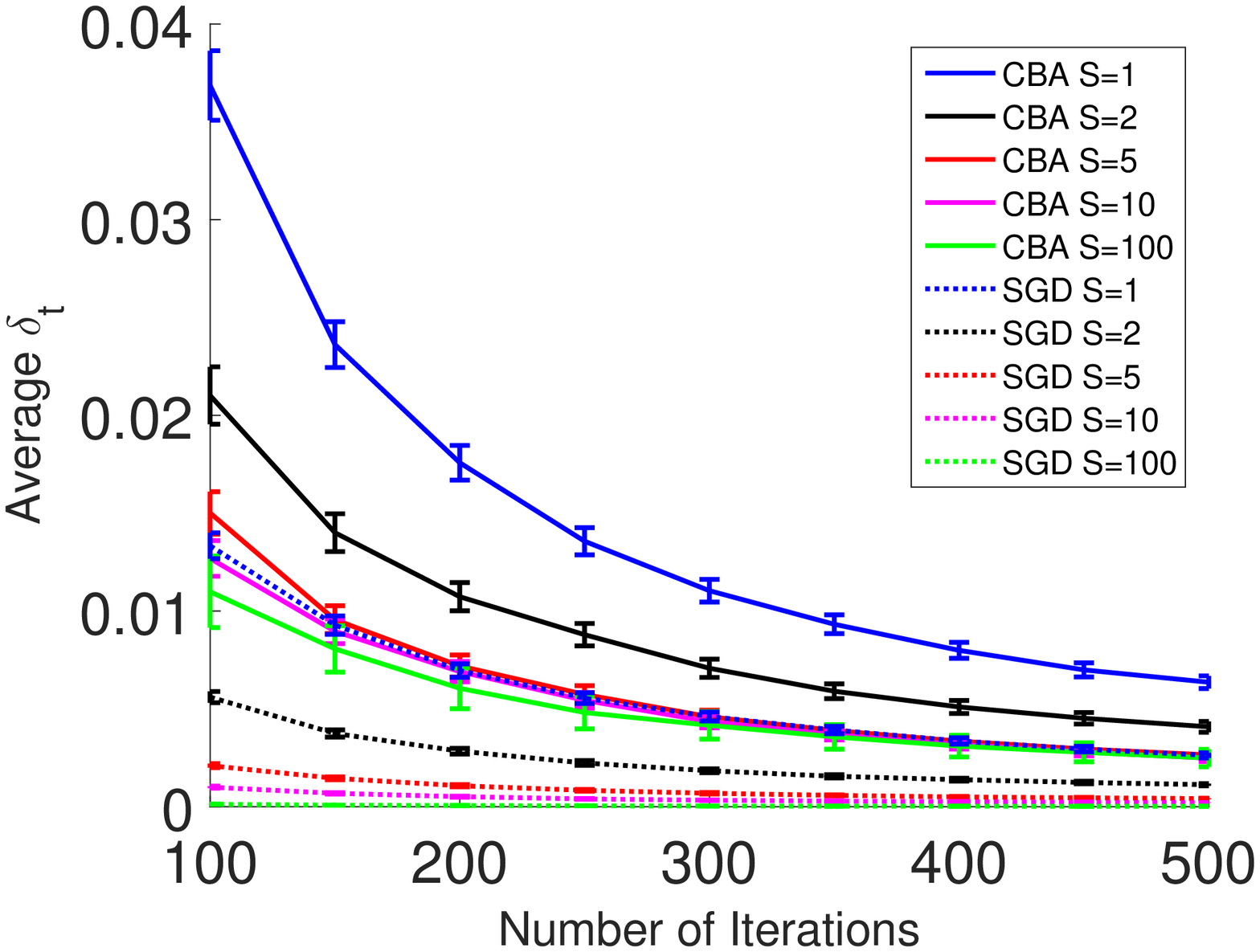}
        \label{fig:normal_pquad_exp}
    }
    \caption{The convergence of the average relative
        optimality gap in CBA (solid lines) and SGD (dot lines) when using different batch sizes ($S$) per iteration.
        First row: $h_1(x,\xi)$; Second row: $h_2(x,\xi)$. First
        column: $\xi\sim\mathcal{U}[50,150]$;  Second column:
        $\xi\sim\mathcal{N}(100,100)$.}
    \label{fig:synthetic_MB_SGD}
\end{figure*}

    After studying the impact of mini-batch method to CBA, we also show how mini-batch affects SGD and whether the performance of CBA relative to SGD we found in Figure~\ref{fig:synthetic}   still holds when using different batch sizes.
    Therefore, we also implement SGD by generating $S$ samples of $\xi$ and using the average of the gradients of $h$ at each sample to perform the stochastic gradient descent step. We also present the convergence of both and CBA and SGD in
    Figure~\ref{fig:synthetic_MB_SGD}  with $S=1,2,5,10,100$ in order to compare their performances with different batch sizes.
    These figures reveal the similar phenomenon as we found in Section~\ref{subsec:numerical:convergence} that CBA converges more slowly than SGD with the same batch size because CBA uses less information. However, CBA does not require knowing $\xi$ exactly so that it can be applied in some scenarios where SGD cannot be implemented.

    \subsection{Categorical Results from Comparison}
    \label{sec:expcategorical}
        We conduct numerical experiments to test the
        performance of the CBA-C (Algorithm~\ref{alg:cba1categorical}) through a comparison with SGD and CBA. We consider the same two objective
        functions $h_1$ and $h_2$ and the same two distributions of $\xi$ as in Section~\ref{subsec:numerical:convergence}. The implementation of SGD and CBA, including the choice of $\eta_t$, $f_+$ and $f_-$, is also the same as in Section~\ref{subsec:numerical:convergence}. The number of the categorical results of the comparisons used in CBA-C is chosen to be $m=3$ or $m=5$. To implement CBA-C, we need to specify $2m$ distributions $f_+^i$ and $f_-^i$ for $i=0,1,\dots,m-1$. When $\theta_{i+1}<+\infty$ so that $ [x_t-\theta_{i+1},x_t-\theta_i)$ and $(x_t+\theta_i,x_t+\theta_{i+1}]$ are bounded intervals, we choose $f_{-}^i$ and $f_{+}^i$ to be uniform distributions on the corresponding intervals. For the cases
        where $\theta_{i+1}=+\infty$ (i.e., $i=m-1$) so that $ [x_t-\theta_{i+1},x_t-\theta_i)$ and $(x_t+\theta_i,x_t+\theta_{i+1}]$ are unbounded in one end, we choose $f_{-}^i$ and $f_{+}^i$   to be exponential distributions as in Example
        \ref{example:fexponential} with $\lambda_{+}=\lambda_{-}=2^{-4}=0.0625$. For the same reasons given in Examples~\ref{example:funiform} and \ref{example:fexponential}, the distribution function $f_+^i$ and $f_-^i$ chosen in this way satisfy (C1')-(C3').

\begin{figure*}[!h]
            \centering
            \subfigure{
                \centering
                \includegraphics[width=0.45\columnwidth]{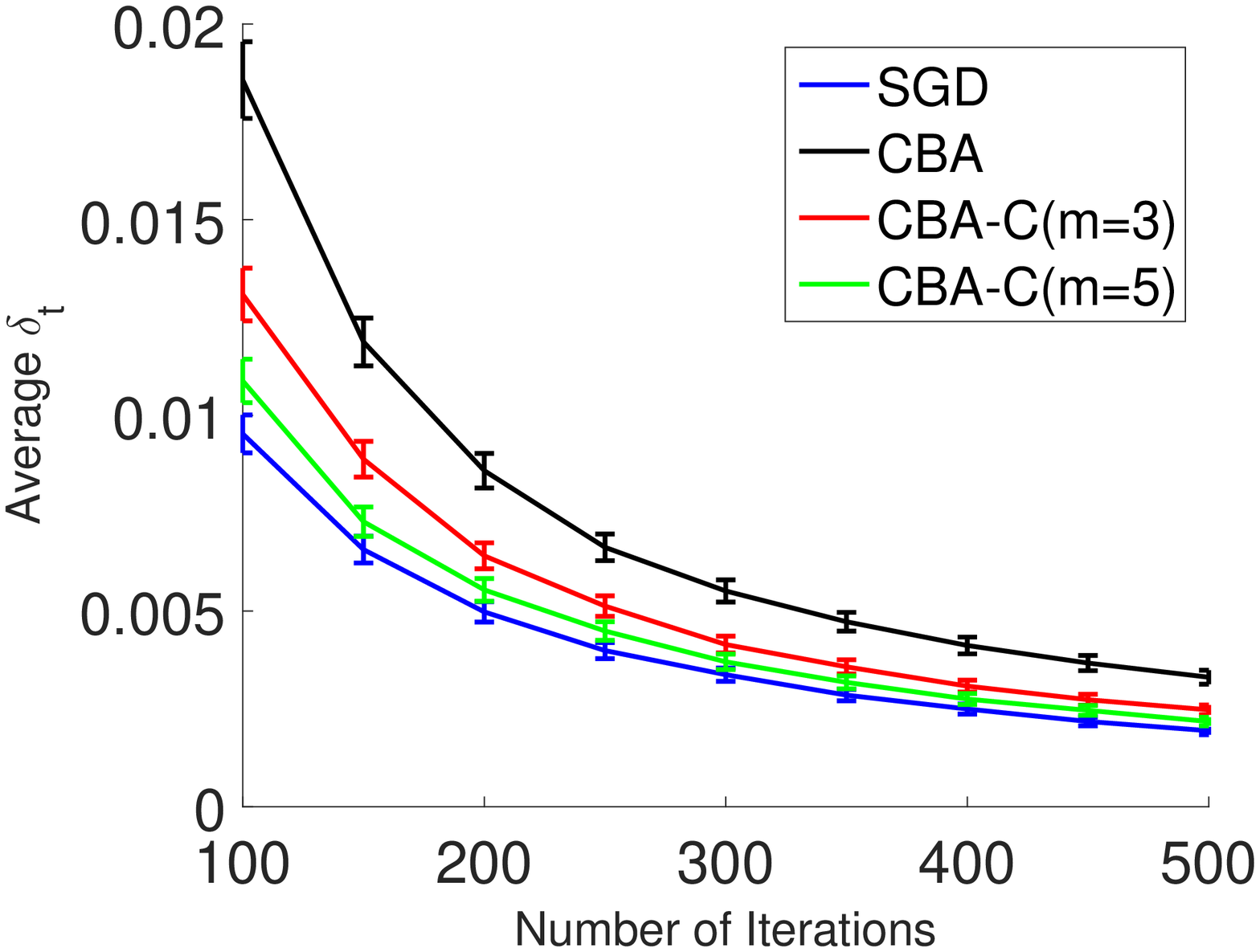}
                \label{fig:uniform_quad_uniform_discrete}
            }
            \subfigure{
                \centering
                \includegraphics[width=0.45\columnwidth]{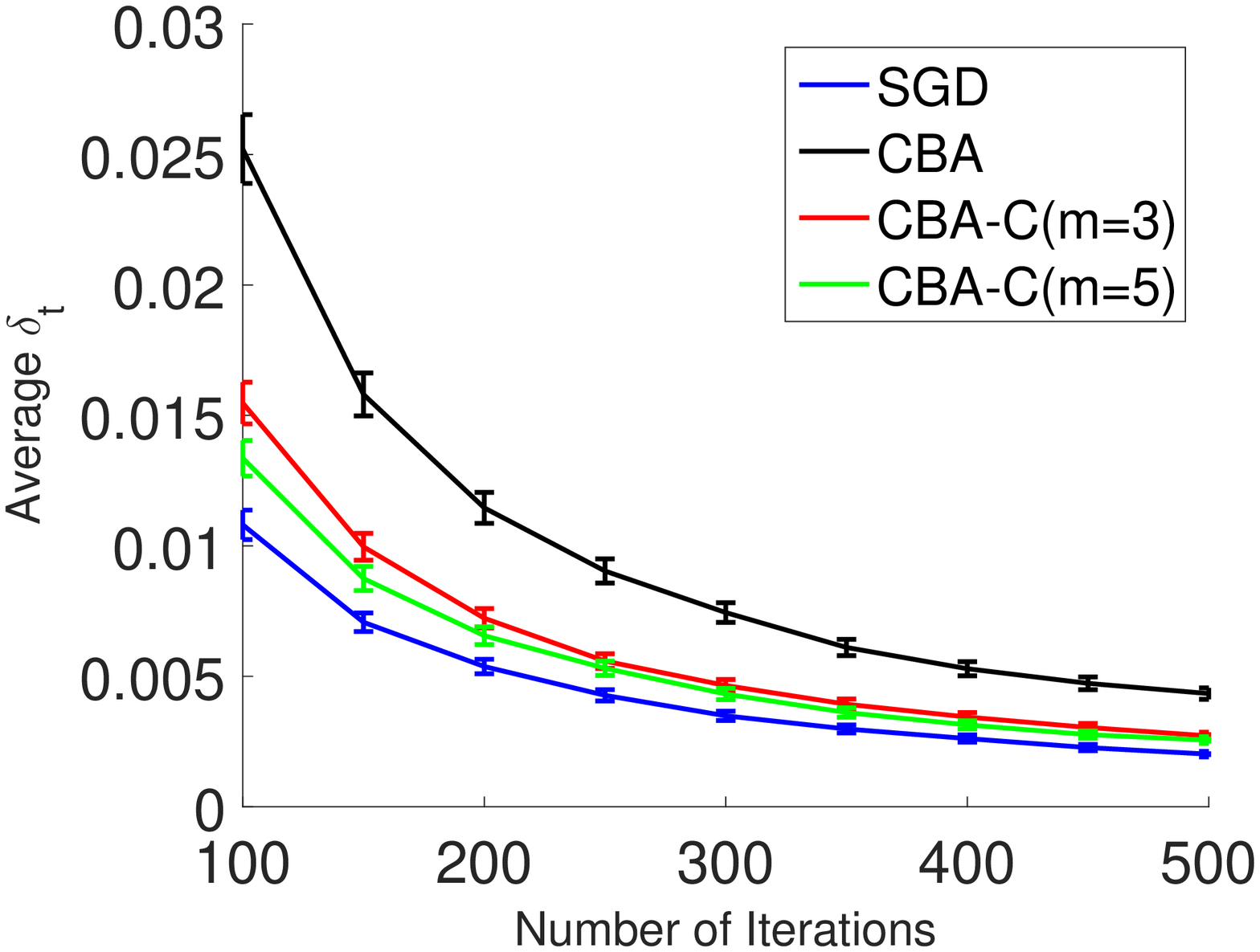}
                \label{fig:normal_quad_exp_discrete}
            }
            \subfigure{
                \centering
                \includegraphics[width=0.45\columnwidth]{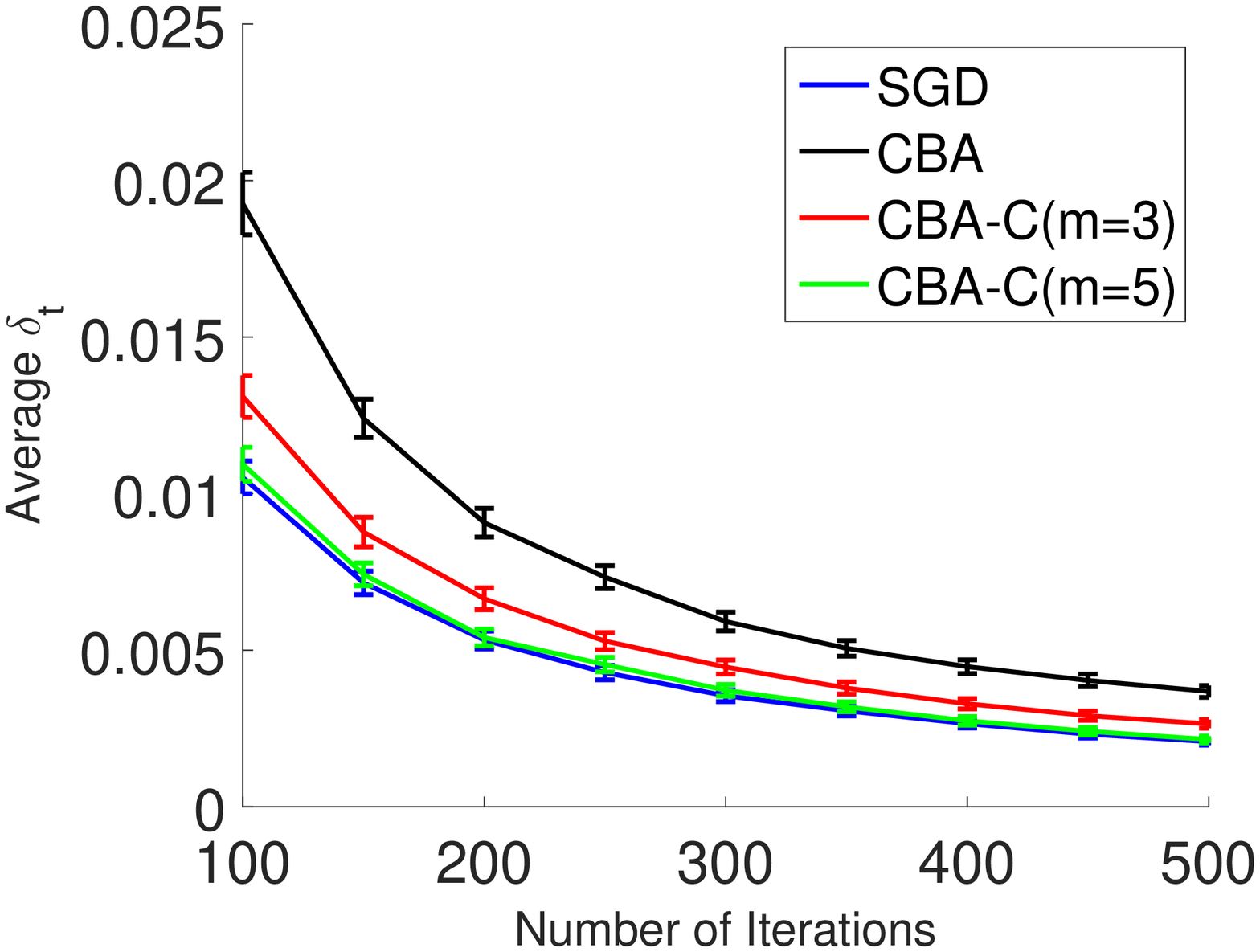}
                \label{fig:uniform_pquad_uniform_discrete}
            }
            \subfigure{
                \centering
                \includegraphics[width=0.45\columnwidth]{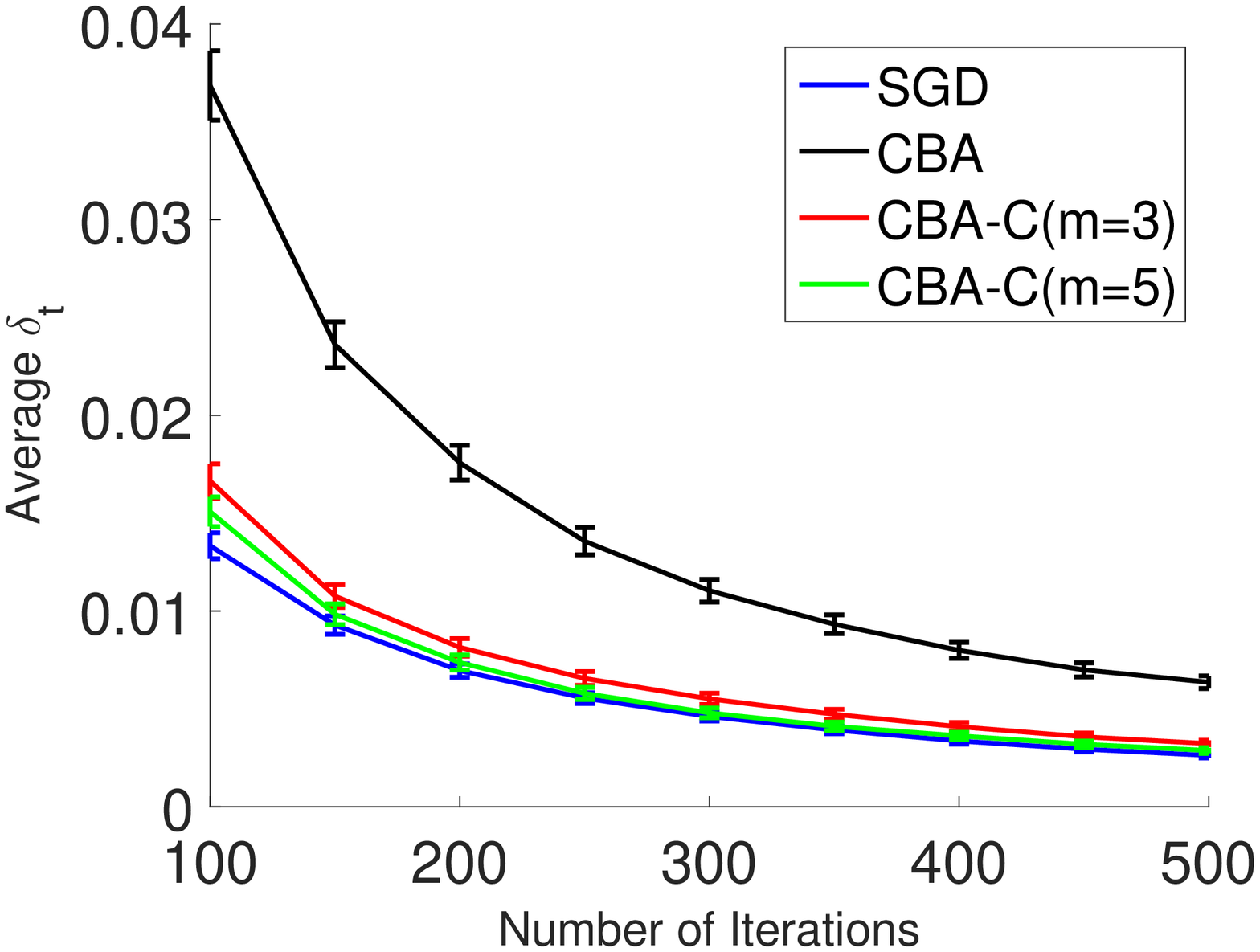}
                \label{fig:normal_pquad_exp_discrete}
            }
            \caption{The convergence of the average relative
                optimality gap for different algorithms for the four instances.
                First row: $h_1(x,\xi)$; Second row: $h_2(x,\xi)$. First
                column: $\xi\sim\mathcal{U}[50,150]$;  Second column:
                $\xi\sim\mathcal{N}(100,100)$.}
            \label{fig:synthetic_discrete}
            \vspace{-0.3cm}
        \end{figure*}

        Similar to Section~\ref{subsec:numerical:convergence}, in all tests, we start from a random
        initial point $x_1 \sim\mathcal{U}[50,150]$ and run each algorithm
        for $T=500$ iterations and we report the average relative optimality
        gap $\delta_t = \frac{H(\bar x_t)-H(x^*)}{H(x^*)}$ over $2000$
        independent trails for each $t$. The results are reported in Figure
        \ref{fig:synthetic_discrete}. The bar at each point represents the standard error of the
        corresponding $\delta_t$.

        \begin{table}[!h]
            \vspace{-0.3cm}
            \begin{center}
                \begin{tabular}{|c|c|c | c | c|}
                    \hline
                    &\multicolumn{2}{c|}{$h_1(x,\xi)$}&\multicolumn{2}{c|}{$h_2(x,\xi)$}\\\hline
                    Algorithm&$\xi\sim\mathcal{U}[50,150]$& $\xi\sim\mathcal{N}(100,100)$&$\xi\sim\mathcal{U}[50,150]$& $\xi\sim\mathcal{N}(100,100)$\\\hline
                    SGD&0.008&0.007&0.010&0.040\\
                    CBA&0.011&0.018&0.013&0.049\\
                    CBA-C$(m=3)$&0.012&0.023&0.023&0.065\\
                    CBA-C$(m=5)$&0.012&0.026&0.025&0.064\\
                    \hline
                \end{tabular}
                \caption{The computation time (in seconds) of different algorithms for 500 iterations for the four instances.
                    \label{table1_discrete}
                }
                \vspace{-0.5cm}
            \end{center}
        \end{table}

        From the results shown in Figure \ref{fig:synthetic_discrete}, we can see
        that CBA-C converges faster than CBA in these
        problems. Moreover, the convergence speed of CBA-C is higher for a larger $m$. This is because the comparison results in CBA-C are categorical instead of binary so that they provide more information about the gap between $x$ and $\xi$ which helps to construct a more precision stochastic gradient than CBA. This information becomes more precise as $m$ increases. Also, as shown in
        Table~\ref{table1_discrete}, the computation time of CBA-C for 500 iterations is less
        than $0.1$ seconds and is similar to CBA.

\end{document}